\documentclass[11pt]{amsart}
\usepackage{amssymb,amsthm}
\usepackage{setspace}

\newtheorem*{124Theorem}{1.24 Theorem}
\newtheorem*{213Lemma}{2.13 Lemma}
\newtheorem*{233Lemma}{2.33 Lemma}
\newtheorem*{320Lemma}{3.20 Lemma}
\newtheorem*{324Corollary}{3.24 Corollary}
\newtheorem*{328Lemma}{3.28 Lemma}
\newtheorem*{334Lemma}{3.34 Lemma}
\newtheorem*{338Theorem}{3.38 Theorem}
\newtheorem*{342Lemma}{3.42 Lemma}
\newtheorem*{349Lemma}{3.49 Lemma}
\newtheorem*{353Lemma}{3.53 Lemma}
\newtheorem*{357Lemma}{3.57 Lemma}
\newtheorem*{360Lemma}{3.60 Lemma}
\newtheorem*{365Lemma}{3.65 Lemma}

\newtheorem*{380Lemma}{3.80 Lemma}
\newtheorem*{44Lemma}{4.4 Lemma}
\newtheorem*{48Corollary}{4.8 Corollary}
\newtheorem*{412Lemma}{4.12 Lemma}
\newtheorem*{415Lemma}{4.15 Lemma}
\newtheorem*{423Theorem}{4.23 Theorem}
\newtheorem*{427Lemma}{4.27 Lemma}
\newtheorem*{439Proposition}{4.39 Proposition}
\newtheorem*{441Lemma}{4.41 Lemma}
\newtheorem*{444Lemma}{4.44 Lemma}
\newtheorem*{447Lemma}{4.47 Lemma}
\newtheorem*{453Lemma}{4.53 Lemma}
\newtheorem*{454Lemma}{4.54 Lemma}
\newtheorem*{459Lemma}{4.59 Lemma}
\newtheorem*{465Lemma}{4.65 Lemma}
\newtheorem*{469Lemma}{4.69 Lemma}
\newtheorem*{476Lemma}{4.76 Lemma}

\newtheorem*{478Lemma}{4.78 Lemma}
\newtheorem*{484Corollary}{4.84 Corollary}
\newtheorem*{487Lemma}{4.87 Lemma}
\newtheorem*{489Lemma}{4.89 Lemma}
\newtheorem*{495Theorem}{4.95 Theorem}
\newtheorem*{57Lemma}{5.7 Lemma}
\newtheorem*{511Corollary}{5.11 Corollary}
\newtheorem*{514Lemma}{5.14 Lemma}
\newtheorem*{524Lemma}{5.24 Lemma}
\newtheorem*{528Corollary}{5.28 Corollary}
\newtheorem*{531Proposition}{5.31 Proposition}
\newtheorem*{533Lemma}{5.33 Lemma}
\newtheorem*{535Proposition}{5.35 Proposition}
\newtheorem*{536Lemma}{5.36 Lemma}
\newtheorem*{538Definition}{5.38 Definition}
\newtheorem*{540Lemma}{5.40 Lemma}
\newtheorem*{542Proposition}{5.42 Proposition}
\newtheorem*{549Lemma}{5.49 Lemma}
\newtheorem*{552Lemma}{5.52 Lemma}
\newtheorem*{558Lemma}{5.58 Lemma}
\newtheorem*{564Corollary}{5.64 Corollary}
\newtheorem*{566Lemma}{5.66 Lemma}
\newtheorem*{572Lemma}{5.72 Lemma}
\newtheorem*{574Lemma}{5.74 Lemma}
\newtheorem*{576Proposition}{5.76 Proposition}
\newtheorem*{577Corollary}{5.77 Corollary}
\newtheorem*{68Theorem}{6.8 Theorem}
\newtheorem*{614Theorem}{6.14 Theorem}
\newtheorem*{616Theorem}{6.16 Theorem}
\newtheorem*{652Theorem}{6.52 Theorem}
\newtheorem*{655Theorem}{6.55 Theorem}

\theoremstyle{remark}

\numberwithin{equation}{section}

\DeclareMathOperator{\Ree}{Re}
\DeclareMathOperator{\Imm}{Im}
\DeclareMathOperator{\sgn}{sgn}

\newcommand{\al}{\alpha}
\newcommand{\BC}{\mathbb C}
\newcommand{\credit}{{\rm cr}}
\newcommand{\de}{\delta}
\newcommand{\De}{\Delta}
\newcommand{\ve}{\varepsilon}
\newcommand{\g}{\gamma}
\newcommand{\G}{\Gamma}
\newcommand{\ka}{\kappa}
\newcommand{\la}{\lambda}
\newcommand{\cL}{\mathcal L}
\newcommand{\om}{\omega}
\newcommand{\Om}{\Omega}
\newcommand{\vp}{\varphi}
\newcommand{\BR}{\mathbb R}
\newcommand{\s}{\sigma}
\newcommand{\vt}{\vartheta}
\newcommand{\p}{\partial}
\newcommand{\BZ}{\mathbb Z}

\newcommand{\dvp}{\dot\vp}
\newcommand{\dvt}{\dot\vt}
\newcommand{\dz}{\dot z}

\newcommand{\sqrtonealtwo}{\sqrt{1-\al^2}}
\newcommand{\sqrtonealplustwo}{\sqrt{1-\al_+^2}}
\newcommand{\sqrtonealminustwo}{\sqrt{1-\al_-^2}}

\newcommand{\oN}{\overline N}

\newcommand{\oU}{\overline U}
\newcommand{\ow}{\overline w}
\newcommand{\oW}{\overline W}
\newcommand{\oz}{\overline z}
\newcommand{\oZ}{\overline Z}

\newcommand{\tgamma}{\tilde\gamma}

\title[A Hamiltonian approach to the heat kernel of a subLaplacian]{A Hamiltonian approach to the heat kernel of a subLaplacian on $S^{2n+1}$}

\author{Peter C. Greiner}
\date{} 
\thanks{ Research partially supported by NSERC Grant OGP0003017.}

\begin{document}

\maketitle

\begin{abstract}
The heat kernel for the Cauchy-Riemann subLaplacian on $S^{2n+1}$ is derived in a manner which is
completely analogous to the classical derivation of elliptic heat kernels.
This suggests that the classical hamiltonian construction
of elliptic heat kernels, with appropriate modifications, will yield heat kernels for subelliptic operators.
\end{abstract}

\section{Introduction}
\label{sec1}

Let $\sqrt 2 Z_1,\ldots, \sqrt 2 Z_{n+1}$ denote an orthonormal basis of holomorphic vector fields on $\BC^{n+1}$
with respect to the Euclidean metric,  $Z_j = \sum\limits_{k=1}^{n+1} a_{j_k} \p / \p z_k$, $j=1,\ldots, n+1$,
and let $Z_j^*$ denote the operator adjoint to $Z_j$ with respect to the Euclidean volume form
\begin{equation}
\label{eq1.1}
dx = dx_1 \wedge \ldots \wedge dx_{2n+2} = \frac{(-1)^{n^2-1}}{2^{n+1}} dz_1 \wedge d\oz_1\wedge \ldots
\wedge dz_{n+1} \wedge d\oz_{n+1} ,
\end{equation}
$z_j = x_j + ix_{j+n+1}$, $j=1,\ldots, n+1$.
The Laplacian on $\BC^{n+1}$ is given by
\begin{equation}
\label{eq1.2}
\De = 2\sum_{j=1}^{n+1} \frac{\p}{\p z_j} \frac{\p}{\p \oz_j}= -\sum_{j=1}^{n+1}
(Z_j^* Z_j + \oZ_j^* \oZ_j) .
\end{equation}
Similarly, the CR-subLaplacian $S^{2n+1}$ is
\begin{equation}
\label{eq1.3}
\De_C = -\sum_{j=1}^n (W_j^* W_j + \oW_j^* \oW_j) \mid_{S^{2n+1}} ,
\end{equation}
where $\sqrt 2 W_1 , \ldots, \sqrt 2 W_n$ represents a local orthonormal basis of those holomorphic
vector fields on $\BC^{n+1}$ which are orthogonal to the holomorphic radial vector field
\begin{equation}
\label{eq1.4}
N = \sum_{j=1}^{n+1} \frac{z_j}{r} \frac{\p}{\p z_j} ,\qquad r^2 = z\cdot\oz ,
\end{equation}
and therefore tangent to spheres.
Then $\sqrt 2 N$, $\sqrt 2 W_1 , \ldots, \sqrt 2 W_n$ is a local orthonormal basis of $T^{(1,0)} \BC^{n+1}$ and one has
\begin{equation}
\label{eq1.5}
\De = -2 \Ree N^* N - 2 \Ree \sum_{j=1}^n W_j^* W_j .
\end{equation}
Let $\De_S$ denote the restriction of $\De$ to $C^\infty (S^{2n+1})$.
Then \eqref{eq1.3} and \eqref{eq1.5} yield
\begin{equation}
\label{eq1.6}
\De_C = \De_S + 2\Ree N^* N \mid_{C^\infty (S^{2n+1})} .
\end{equation}
Thus $\De_C$ is globally defined on $S^{2n+1}$ and is fully explicit in view of \eqref{eq1.2} and \eqref{eq1.4}.
We note that $\De_S$ is elliptic but $\De_C$ is not.

Let $P_C (t,w,z)$ denote the heat kernel of $\De_C$ with pole at $w$.
$P_C$ is characterized by
\begin{equation}
\label{eq1.7}
\frac{\p P_C}{\p t} = \De_C P_C ,
\end{equation}
\begin{equation}
\label{eq1.8}
\lim_{t\to 0} \, P_C = \de (z-w) .
\end{equation}
$\De_C$ is invariant with respect to complex (unitary) rotations of $\BC^{n+1}$, hence $P_C (t)$ is a function
of one complex variable
\begin{equation}
\label{eq1.9}
z\cdot \ow = z_1 \ow_1 + \cdots + z_{n+1} \ow_{n+1} = \cos\vt_1 e^{i\vp_1} ,
\end{equation}
and its complex conjugate $\oz\cdot w$, or of 2 real variables $\vt_1 $, $\vp_1$, $0\leq \vt_1 \leq \frac{\pi}{2}$,
$-\pi < \vp_1 \leq \pi$.
Thus, with a slight abuse of notation, we may write
\begin{equation}
\label{eq1.10}
P_C (t,w,z) = P_C (t, |z\cdot w|, \arg z\cdot w) = P_C (t,\cos\vt_1 , \vp_1) .
\end{equation}
In particular we may assume that $w=(1,0,\ldots, 0)$, and then $z_1 = z\cdot \ow = \cos\vt_1 e^{i\vp_1}$.
One extends $\vt_1 , \vp_1$ to a complete system of spherical coordinates $(\vt,\vp) = (\vt_1 ,\ldots, \vt_n ,
\vp_1 , \ldots, \vp_{n+1})$ on $S^{2n+1}$, see \eqref{eq2.1}, and then $z$ is represented by $(\vt,\vp)$.
$\De_C$ restricted to functions of $\vt_1,\vp_1$ yields the reduced operator $\cL_C$,
\begin{equation}
\label{eq1.11}
\cL_C = \frac12 \frac{\p^2}{\p\vt_1^2} + \big( (n-1) \cot\vt_1 + \cot 2\vt_1\big) \frac{\p}{\p\vt_1} +\frac12
\tan^2\vt_1 \frac{\p^2}{\p\vp_1^2} .
\end{equation}
The heat kernel of $\cL_C$ is also the heat kernel of $\De_C$ modulo a normalizing factor. $\cL_C$ is not elliptic
since $(\tan^2\vt_1)\p^2 / \p\vp_1^2$ vanishes at $\vt_1 =0$.
On the other hand one has
\begin{equation}
\label{eq1.12}
\left[ \frac{\p}{\p\vt_1} , (\tan\vt_1)\frac{\p}{\p\vp_1}\right] = \frac{1}{\cos^2\vt_1} \frac{\p}{\p\vp_1} ,
\end{equation}
which does not vanish at $\vt_1=0$, hence $\cL_C$ is step 2 subelliptic.
In \S\ref{sec4} the heat kernel $P_C$ of $\De_C$ is obtained in the following form:
{\it 
\begin{equation}
\label{eq1.13}
P_C = \frac{e^{\frac{n^2}{2}t}}{(2\pi t)^{n+1}} \sum_{k=-\infty}^{\infty} \int_{-\infty}^\infty e^
{-\frac{f(u,\kappa +i2k\pi)}{2t}} V_n (\kappa+i2k\pi, t) du ,
\end{equation}
\begin{equation}
\label{eq1.14}
f(u,\kappa) = u^2 - \kappa^2 = u^2 - \big( \cosh^{-1} (\cos\vt_1 \cosh (u+i\vp_1) ) \big)^2 ,
\end{equation}
\begin{equation}
\label{eq1.15}
V_n (\kappa,t) = V^n W_n = \left( \frac{\kappa}{\sinh \kappa}\right)^n \left(\sum_{\ell=0}^{n-1} W_{n,\ell} (\kappa) t^\ell\right) ,
\end{equation}
$W_{n,0}=1$ and $W_{n,\ell}$, $0 < \ell \leq n-1$, are found by iteration,
\begin{equation}
\label{eq1.16}
W_{n,\ell} (\kappa) = \frac{1}{\kappa^\ell} \int_0^\kappa V^{-n} x^{\ell-1} \Big(\cL_S - \frac{n^2}{2}\Big)
V_{n,\ell-1} dx ,
\end{equation}
where we set $V_{n,\ell} (\kappa) = V(\kappa)^n W_{n,\ell} (\kappa)$ and $\cL_S$ is given by
\begin{equation}
\label{eq1.17}
\cL_C = \cL_S - \frac12 \frac{\p^2}{\p\vp_1^2} .
\end{equation}
}

\eqref{eq1.13} has geometric meaning.
To facilitate its description let me recall the heat kernel $P_S$ of $\De_S$:
\begin{equation}
\label{eq1.18}
P_S = \frac{e^{\frac{n^2}{2}t}}{(2\pi t)^{n+1/2}} \sum_{k=-\infty}^{\infty} e^{-\frac{(\gamma + 2k\pi)^2}{2t}}
v_n (\gamma +2k\pi,t) = \sum_{k=-\infty}^\infty P_{S,k},
\end{equation}
where the angle $\gamma$, $0\leq \gamma < \pi$ subtends the
points $z,w\in S^{2n+1}$, $\cos\g = x\cdot y = \Ree z \cdot w$
and $V_n (i\gamma,t) = v_n(\g,t)$.
$P_S$ is a solution of $\De_S \, P_S = \p P_S / \p t$,  hence so are $P_{S,k}$,
$k\in\BZ$. 
Write
\begin{equation}
\label{eq1.19}
P_{S,k} = \frac{1}{(2\pi t)^{n+1/2}} e^{-\frac{(\g+2k\pi)^2}{2t}} (a_0 + a_1 t + \cdots ) ,
\end{equation}
and note that the coefficient of the principal singularity in $t$ in $\De_S P_{S,k} - \p P_{S,k} / \p t$ is
\begin{equation}
\label{eq1.20}
e^{-\frac{(\g+2k\pi)^2}{2t}} \left\{ \frac12 \left( \frac{\p (-(\g+2k\pi)^2/2t)}{\p\g}\right)^2 - \frac{\p (-(\g+2k\pi)^2/2t)}{\p t}\right\} = 0 ,
\end{equation}
that is, $-(\g+2k\pi)^2/2t$ is a solution of the Hamilton-Jacobi equation
 induced by \eqref{eq3.43}.
$\g = \cos^{-1} (\cos\vt_1 \cos\vp_1)$ is the
Riemannian distance between $z$ and $w$ in $S^{2n+1}$.
Vice-versa, given a heat kernel one may obtain the Riemannian distance from its small time asymptotic.
For example, in \S\ref{sec3} we derive
\begin{equation}
\label{eq1.21}
P_S (t\to 0) \sim \frac{1}{(2\pi t)^{n+1/2}} e^{-\frac{\g^2}{2t}} \left( \frac{\g}{\sin \g}\right)^n .
\end{equation}

\eqref{eq1.13} is obtained by a completely analogous calculation.
In particular, an action integral yields $f=u^2 - (\kappa + i2k\pi)^2$, $k=0,\pm 1,\ldots$, where $u$
is essentially the dual of the missing direction $\vp_1$.
In view of \eqref{eq4.28} the coefficient of the principal singularity in $t$ in
\begin{equation}
\label{eq1.22}
\left( \De_C - \frac{\p}{\p t}\right) \frac{1}{(2\pi t)^{n+1}} e^{-\frac{f}{2t}}
\end{equation}
is given by
\begin{equation}
\label{eq1.23}
-e^{-\frac{f}{2t}} \left\{\frac12 \left( \frac{\p (-f/2t)}{\p\kappa}\right)^2 + \frac12 \left( \frac{\p (-f/2t)}{\p\vp_1}\right)^2 -
\frac{\p (-f/2t)}{\p t}\right\} .
\end{equation}
This turns out to be a $u$-derivative whose integral over $-\infty < u < \infty$ vanishes.
Note the difference.
One needs \eqref{eq1.20} to vanish but we need only the integral of
\eqref{eq1.23} to vanish and that is exactly what we get.
$V_{n,0}$ is a solution of a first order differential equation and the rest of the $V_{n,\ell}$'s are obtained by quadrature.
$f(u) = f\big( u,\kappa(u)\big)$ is the square of the length of a CR-geodesic at the critical points
$u_{\credit}$ of $f$, i.e. where $f'(u_\credit)=0$.
This yields the small time asymptotic for $P_C$ by the method of stationary phase and the length of the
shortest geodesic, or the Carnot-Caratheodory distance.
The relevant results are (all CR-geodesics are local, i.e.
$-\pi < \vp_1 < \pi$):

\begin{124Theorem}
{\rm (i)} All the critical points of $f(u)$ are on the imaginary axis.

{\rm (ii)} Given $z=(\vt,\vp) \in S^{2n+1}$, $\vt_1 \ne 0$, $f(u)$ has at least one and at most a finite number of 
critical points.

{\rm (iii)} When $\vt_1=0$, $0< |\vp_1| < \pi$, $f(u)$ has a discrete infinity of critical points.

{\rm (iv)} There is a 1-1 correspondence between the set of critical points of $f$ and the set of the lengths
of CR-geodesics joining $(1,0,\ldots,0)$ and $z$.

{\rm (v)} When $\vp_1 > 0$ ($\vp_1<0$), let $u_1(z), u_2(z),\ldots$ denote the critical points of $f(u)$
ordered by increasing modulus.
Then $f\big( u_j (z)\big)$ represents the square of the length of the geodesic attached to $u_j (z)$, and
\addtocounter{equation}{1}
\begin{equation}
\label{eq1.25}
f\big( u_1(z)\big) < f\big( u_2(z)\big) < \cdots
\end{equation}
The shortest geodesic length is the Carnot-Caratheodory distance $d_c(z)$ of the point $z$ from
$(1,0,\ldots, 0)$ where $d_c^2 (z) = f\big( u_c (z)\big)$, $u_c(z) = u_1 (z)$.

{\rm (vi)} Let $\eta$ denote the unique solution of
\begin{equation}
\label{eq1.26}
\tan \left( \vp_1 + \eta \sqrt{1 - \frac{\sin^2\vt_1}{\sin^2 \eta}}\right) = \sqrt{1 - \frac{\sin^2\vt_1}{\sin^2\eta}}\,
\tan\eta
\end{equation}
in $\vt_1 \leq \eta \leq \pi - \vt_1$ when $\vp_1 \geq 0$ and in $-\pi + \vt_1 < \eta < -\vt_1$
when $\vp_1 < 0$.
Then
\begin{equation}
\label{eq1.27}
d_c (z) = \sin\vt_1 \frac{\eta}{\sin\eta} .
\end{equation}

{\rm (vii)}
\eqref{eq1.13} and the method of the stationary phase yield the small time asymptotic of $P_C$, see \eqref{eq6.53},
\eqref{eq6.56}.

{\rm (viii)}
The canonical curve through $(1,0,\ldots,0)$ is represented by $\vt_1 =0$, $-\pi < \vp_1 \le \pi$; this is the
set of points which can be joined to $(1,0,\ldots,0)$ by a discrete infinity of 
CR-geodesics of different lengths.
\end{124Theorem}

The concept of a canonical curve, or, in general, of a canonical submanifold through a point was suggested in
\cite{9}, \cite{6}, \cite{7}, where the canonical curve in a step 3 case was carefully worked out.
\cite{1}, \cite{2} note that $\cL_C$ is a sum of 2 commuting elliptic operators and a 
convolution of their
heat kernels, obtained using special functions, yields $P_C$ in the form \eqref{eq4.40}.
This observation was already used by Staubach in \cite{15} to find the heat kernel on the Heisenberg group;
he composed the two known heat kernels via the Trotter product formula and path integration.
Our motivation for finding $P_C$ came from Yum-Tong Siu's interest in CR-zeta functions;
note that the Mellin transform of the heat kernel is the zeta function.
Furthermore, heat kernels yield all functions of the operator.
The CR-zeta function on $S^{2n+1}$ has been studied by Der-Chen Chang and Song-Ying Li in \cite{4} and \cite{14} studies CR-geodesics.

Explicit expressions for $\De_S$, $\De_C$, $\cL_S$ and $\cL_C$ in spherical coordinates are derived in \S\ref{sec2}.
As preparation for the construction of $P_C$ we carefully derive $P_S$ via a Hamiltonian formalism in \S\ref{sec3}.
$P_C$ is obtained in \S\ref{sec4} in complete analogy with the work done for $P_S$ in \S\ref{sec3}.
\S\ref{sec5} is devoted to finding the subRiemannian (CR-) geodesics, and their connection with the complex
action $f$ is explored in \S\ref{sec6}; this also yields the small time asymptotic of $P_C$.
We mostly work with $0\leq \vp_1 \leq \pi$.
The calculations with $-\pi < \vp_1 < 0$ are completely analogous, here $\eta<0$, and we shall not redo them
since we do not wish to burden this already lengthy article with unenlightening repetition.

The construction of $P_C$ is a word-for-word analogue of the construction of $P_S$.
This suggests a pattern and introduces the possibility that the Hamiltonian formalism will yield heat kernels for
subelliptic operators in general.

\section{Laplacians in spherical coordinates}
\label{sec2}
Let
\begin{align}
\label{eq2.1}
&z_1 = r \cos\vartheta_1 e^{i\vp_1} ,\\
&\cdots\nonumber\\
&z_k = r \sin\vartheta_1\ldots \sin\vartheta_{k-1}\cos\vartheta_k e^{i\vp_k} ,\qquad k=2,\ldots, n,\nonumber\\
&\cdots\nonumber\cr
&z_{n+1} = r\sin\vt_1\ldots\sin\vt_{n-1}\sin\vt_n e^{i\vp_{n+1}},\nonumber
\end{align}
$0\leq \vt_j \leq \pi/2$, $j=1,\ldots,n$, $-\pi < \vp_j\leq \pi$, $j=1,\ldots,n+1$, and $0\leq r < \infty$ denote spherical
coordinates in $\BC^{n+1}$.
Then
\begin{equation}
\label{eq2.2}
\frac{\p r}{\p z_k} = \frac{\oz_k}{2 r} ,\qquad k=1,\ldots, n+1 ,
\end{equation}
and
\begin{equation}
\label{eq2.3}
|z_k|^2 \tan^2\vt_k = |z_{k+1}|^2 + \cdots + |z_{n+1}|^2 ,\qquad k=1,\ldots,n
\end{equation}
imply
\begin{align}
\label{eq2.4}
\frac{\p\vt_k}{\p z_j} &= \frac12 \cos^2\vt_k \cot\vt_k \frac{\oz_j}{|z_k|^2} ,\qquad k<j , \\
\label{eq2.5}
\frac{\p\vt_k}{\p z_k} &= - \frac{\sin 2 \vt_k}{4 z_k} .
\end{align}
Also, $z_k/\oz_k = e^{2i \vp_k}$, so
\begin{equation}
\label{eq2.6}
\frac{\p\vp_k}{\p z_j} = \frac{\delta_{jk}}{2i z_k} ,
\end{equation}
and one has
\begin{align}
\label{eq2.7}
&\frac{\p}{\p z_k} = \frac{\oz_k}{2r} \frac{\p}{\p r} + \frac{1}{2i z_k} \frac{\p}{\p\vp_k} +
\sum_{j=1}^{\min (k,n)} \frac{\p\vt_j}{\p z_k} \frac{\p}{\p\vt_j}\\
&= \frac{\oz_k}{2} \left\{ \frac{1}{r} \frac{\p}{\p r} + \sum_{j=1}^{k-1} \frac{\cos^2 \vt_j \cot\vt_j}{|z_j|^2} \frac{\p}{\p\vt_j} 
- \frac{\sin 2\vt_k}{2|z_k|^2} \frac{\p}{\p\vt_k} - \frac{i}{|z_k|^2} \frac{\p}{\p\vp_k}\right\} ,\nonumber
\end{align}
$k=1,\ldots,n+1$; $\sum\limits_{k=1}^0$ and $\p/\p\vt_{n+1}$ vanish.
We write $\De=\De^{(2)} + \De^{(1)}$, where $\De^{(j)}$ denotes the homogeneous part of $\De$ with the $j$-th order
derivatives.
In particular, $\De^{(2)}$ is the second order part of
\begin{equation*}
\frac12 \sum_{k=1}^{n+1} \left\{ |z_k|^2 \left( \frac1r \frac{\p}{\p r} + \sum_{j=1}^{k-1} \frac{\cos^2 \vt_j \cot\vt_j}{|z_j|^2}
\frac{\p}{\p\vt_j}
-\frac{\sin 2\vt_k}{2|z_k|^2}\frac{\p}{\p\vt_k}\right)^2 + \frac{1}{|z_k|^2} \frac{\p^2}{\p\vp_k^2}\right\} .
\end{equation*}
Squaring the round bracket the cross terms in the second order part vanish, so
\[
2\De^{(2)} = \frac{\p^2}{\p r^2} + \sum_{j=1}^n \frac{\cos^2 \vt_j}{|z_j|^2} \frac{\p^2}{\p\vt_j^2} +
\sum_{j=1}^{n+1} \frac{1}{|z_j|^2}\frac{\p^2}{\p\vp_j^2} .
\]
To find $\De^{(1)}$ we first differentiate the coefficients of $4\p / \p\oz_k$ with respect to $z_k$ to find
\begin{align*}
&\frac1r \left( 2 - \frac{|z_k|^2}{r^2}\right) \frac{\p}{\p r} + \frac{\sin 4\vt_k}{4|z_k|^2} \frac{\p}{\p\vt_k} \\
&+ \sum_{j=1}^{k-1} \frac{\cos^2 \vt_j \cot\vt_j}{|z_j|^2} \left( 2 - \frac{|z_k|^2}{|z_j|^2} \cot^2 \vt_j
(1+2\sin^2 \vt_j)\right) \frac{\p}{\p\vt_j} ,\nonumber
\end{align*}
which summed over $k$ yields
\begin{align*}
2\De^{(1)} &= \frac{2n+1}{r} \frac{\p}{\p r} + \sum_{j=1}^n \left\{ \frac{\sin 4\vt_j}{4|z_j|^2}\right. \\
&\quad  + \frac{\cos^2\vt_j \cot\vt_j}{|z_j|^2} \sum_{k=j+1}^{n+1}\left. \left( 2 - \frac{|z_k|^2}{|z_j|^2}
\cot^2\vt_j (1+ 2\sin^2 \vt_j)\right)\right\} \frac{\p}{\p \vt_j} \nonumber \\
&= \frac{2n+1}{r} \frac{\p}{\p r} + \sum_{j=1}^n \left\{ \frac{(2\cos^2 \vt_j -1)\sin 2\vt_j}{2|z_j|^2}\right. \nonumber\\
&\quad \left. + \frac{\cos^2 \vt_j \cot\vt_j}{|z_j|^2}  (2n+1 - 2j - 2\sin^2 \vt_j) \right\} \frac{\p}{\p \vt_j} \nonumber \\
&= \frac{2n+1}{r} \frac{\p}{\p r} + 2\sum_{j=1}^n \frac{\cos^2 \vt_j}{|z_j|^2} \big( (n-j)\cot\vt_j + \cot 2\vt_j\big)
\frac{\p}{\p\vt_j} . \nonumber
\end{align*}
Thus we have derived
\begin{align}
\label{eq2.8}
\De &= \frac12 \left( \frac{\p^2}{\p r^2} + \frac{2n+1}{r} \frac{\p}{\p r}\right) + \frac{1}{r^2} \De_S \\
&= \frac12 \frac{1}{r^{2n+1}} \frac{\p}{\p r} r^{2n+1} \frac{\p}{\p r} + \frac{1}{r^2} \De_S \nonumber
\end{align}
with
\begin{align}
\label{eq2.9}
\De_S = \De_{S,1} &+ \sum_{j=2}^n \frac{1}{\sin^2\vt_1\ldots \sin^2\vt_{j-1}}
\De_{S,j} \\
&+ \frac12 \frac{1}{\sin^2\vt_1 \ldots \sin^2 \vt_n} \frac{\p^2}{\p\vp_{n+1}^2} ,\nonumber
\end{align}
where we set
\begin{equation}
\label{eq2.10}
\De_{S,j}=\frac12 \frac{\p^2}{\p\vt_j^2} + \big( (n-j) \cot\vt_j + \cot 2\vt_j\big)
\frac{\p}{\p\vt_j} + \frac12 \frac{1}{\cos^2\vt_j} \frac{\p^2}{\p\vp_j^2} .
\end{equation}
$\De$ is symmetric  on $\BC^{n+1}$ with respect to the Euclidean volume element.
So is $\De_S$ on $S^{2n+1}$ with respect to the induced volume element.
Indeed,
\begin{align}
\label{eq2.11}
\De_S = &\tfrac12 \sum_{j=1}^n \frac{r^2\cos^2 \vt_j}{|z_j|^2} \cdot \\
&\qquad \cdot \left\{ \frac{1}{\sin^{2(n-j)} \vt_j \sin 2\vt_j} \frac{\p}{\p\vt_j} \sin^{2(n-j)} \vt_j \sin 2\vt_j
\frac{\p}{\p\vt_j}\right. \nonumber\\
&\qquad\quad \left.+ \frac{1}{\cos^2 \vt_j} \frac{\p^2}{\p\vp_j^2}\right\} + \frac{r^2}{|z_{n+1}|^2}
\frac{\p^2}{\p\vp_{n+1}^2} ; \nonumber
\end{align}
note that the coefficient of the curly bracket is independent of $\vt_j$, and the symmetry of $\De_S$
follows from (\ref{eq2.11}) if the induced volume element $dS$ on $S^{2n+1}$ is given by
\begin{equation}
\label{eq2.12}
dS = \frac{1}{2^n} \prod_{j=1}^n \big( \sin^{2(n-j)} \vt_j \sin 2\vt_j\big) d\vt d\vp ,
\end{equation}
$d\vt = d\vt_1 \ldots d\vt_n$, $d\vp = d\vp_1 \ldots d\vp_{n+1} $.
Now
\[
\prod_{j=1}^{n+1} dz_j \wedge d\oz_j = 2^{n+1} i^{n^2-1} dx_1\wedge\ldots\wedge dx_{2n+2} 
= 2^{n+1} i^{n^2-1} r^{2n+1} dr dS ,
\]
and then (\ref{eq2.12}) is a consequence of

\begin{213Lemma}
One has
\addtocounter{equation}{1}
\begin{equation}
\label{eq2.14}
\prod_{j=1}^{n+1} dz_j \wedge d\oz_j 
= 2i^{n^2-1} r^{2n+1} \prod_{j=1}^n \big( \sin^{2(n-j)} \vt_j \sin 2\vt_j\big) dr d\vt d\vp .
\end{equation}
\end{213Lemma}

\begin{proof}
An elementary calculation gives the result when $n=1$, i.e. in $\BC^2$.
Assuming it holds in $\BC^n$ we write $\BC^{n+1} = \BC\times \BC^n$.
Then
\begin{align}
\label{eq2.15}
&\prod_{j=1}^{n+1} dz_j \wedge d\oz _j \\
&= dz_1 \wedge d\oz_1 
\wedge (r\sin\vt_1)^{2n-1} 2i^{(n-1)^2-1} \prod_{j=2}^n \big( \sin^{2(n-j)} \vt_j \sin 2\vt_j\big) \cdot
\nonumber \\
&\quad \cdot d(r\sin \vt_1) d\vt_2 \wedge \ldots \wedge d\vt_n \wedge d\vp_2 \wedge \ldots \wedge d\vp_{n+1}
\nonumber
\end{align}
which implies (\ref{eq2.14}).
\end{proof}

Let $\sqrt 2 Z_1 ,\ldots, \sqrt 2 Z_n$ denote an orthonormal set of holomorphic vector fields in some domain
in $\BC^{n+1}$ which are orthogonal to $N$ of (\ref{eq1.4}) and let $\sqrt 2 Z_1^* ,\ldots, \sqrt 2 Z_n^*$ denote
the adjoint operators.
We set
\begin{equation}
\label{eq2.16}
\De^{(C)} = - \sum_{j=1}^n (Z_j^* Z_j + \oZ_j^* \oZ_j) = -2\Ree \sum_{j=1}^n Z_j^* Z_j .
\end{equation}
$\De^{(C)}$ is independent of the choice of the defining orthonormal set.
Indeed, let $U = (u_{jk})$ represent a smooth $U(n)$ valued function in a domain in $\BC^{n+1}$,
and set
\begin{align}
\label{eq2.17}
W_j &= \sum_{k=1}^n u_{jk} Z_k ,\qquad j=1,\ldots, n .\\
\noalign {\noindent Then}
\label{eq2.18}
W_j^* &= \sum_{k=1}^n Z_k^* \overline{u_{jk}} ,
\end{align}
and
\begin{align}
\label{eq2.19}
\sum_{j=1}^n W_j^* W_j &= \sum_{j=1}^n \sum_{k=1}^n Z_k^* \overline{u_{jk}} \sum_{\ell=1}^n u_{j\ell} Z_\ell \\
&= \sum_{k,\ell=1}^n Z_k^* \left( \sum_{j=1}^n \overline {u_{jk}} u_{j\ell}\right) Z_\ell \nonumber \\
&=\sum_{k=1}^n Z_k^* Z_k .\nonumber
\end{align}
The same argument also implies that
\begin{equation}
\label{eq2.20}
\De = - \sum_{j=1}^{n+1} (Z_j^* Z_j + \overline{Z_j^*}\, \oZ_j ) ,
\end{equation}
where $\sqrt 2 Z_1,\ldots, \sqrt 2 Z_{n+1}$ represents an arbitrary orthonormal basis of the holomorphic
vector fields in $\BC^{n+1}$; recall that $\De$ was defined by (\ref{eq1.2}) with $Z_j = \p/\p z_j$.
Let $Z_{n+1}=N$. Then
\begin{equation}
\label{eq2.21}
\De = -2 \Ree \sum_{j=1}^n Z_j^* Z_j - 2\Ree N^* N ,
\end{equation}
or,
\begin{equation}
\label{eq2.22}
\De^{(C)} = \De + 2\Ree N^* N .
\end{equation}
Now,
\begin{equation}
\label{eq2.23}
\sum_{j=1}^{n+1} a_j \frac{\p}{\p z_j} \perp N \ \Rightarrow \ \sum_{j=1}^{n+1} a_j \oz_j = 0.
\end{equation}
Consequently,
\begin{equation}
\label{eq2.24}
\sum_{j=1}^{n+1} a_j \frac{\p r}{\p z_j} = \frac{1}{2r} \sum_{j=1}^{n+1} a_j \oz_j = 0,
\end{equation}
that is, a holomorphic vector field orthogonal to $N$ in complex Euclidean metric is tangent
to spheres centred at the origin.
Let $Z$ denote a holomorphic vector field in $\BC^{n+1}$ which is orthogonal to $N$ and let $Z_S$
denote its restriction to $S^{2n+1}$; $Z_S\in T (S^{2n+1})$.
Then an
integration by parts shows that
\begin{equation}
\label{eq2.25}
(Z^*)_S = (Z_S)^* ,
\end{equation}
where $(Z_S)^*$ is the adjoint operator with respect to the measure $dS$ on $S^{2n+1}$.
In view of this we define the Cauchy-Riemann subLaplacian $\De_C$ to be the restriction
of $\De^{(C)}$ to the unit sphere:
\begin{equation}
\label{eq2.26}
\De_C = \De^{(C)}\mid_{S^{2n+1}} = - \sum_{j=1}^n \left\{ (Z_{j,S})^* Z_{j,S} + \overline{(Z_{j,S})^*}\,
\overline{Z_{j,S}}\right\} .
\end{equation}
The invariance argument of (\ref{eq2.19}) implies the following:
given an orthonormal set of holomorphic vector fields $\sqrt 2 Z_1,\ldots, \sqrt 2 Z_n\in T(S^{2n+1})$,
one has
\begin{equation}
\label{eq2.27}
\De_C = - 2 \Ree \sum_{j=1}^n Z_j^* Z_j .
\end{equation}
Note that
\begin{equation}
\label{eq2.28}
N = \frac12 \frac{\p}{\p r} - \frac{i}{2r} \sum_{j=1}^{n+1} \frac{\p}{\p \vp_j}
\end{equation}
in spherical coordinates.
Integrating by parts
\begin{equation}
\label{eq2.29}
-N^* = \oN + \frac{2n+1}{2r} = \frac12 \frac{\p}{\p r} + \frac{2n+1}{2r} + \frac{i}{2r} \sum_{j=1}^{n+1} \frac{\p}{\p\vp_j} ,
\end{equation}
\begin{align}
\label{eq2.30}
-2N^*N = &\frac12 \left( \frac{\p^2}{\p r^2} + \frac{2n+1}{r} \frac{\p}{\p r}\right) + \frac{1}{2r^2}
\left( \sum_{j=1}^{n+1} \frac{\p}{\p\vp_j}\right)^2 \\
&- i \frac{n}{r^2} \sum_{j=1}^{n+1} \frac{\p}{\p\vp_j} .\nonumber
\end{align}
In view of (\ref{eq2.8}), (\ref{eq2.22}) and \eqref{eq2.30}, one has
\begin{equation}
\label{eq2.31}
\De^{(C)} = \frac{1}{r^2} \De_S - \frac{1}{2r^2} \left( \sum_{j=1}^{n+1} \frac{\p}{\p\vp_j}\right)^2 ,
\end{equation}
and therefore
\begin{equation}
\label{eq2.32}
\De_C = \De_S -\frac12 \left( \sum_{j=1}^{n+1} \frac{\p}{\p\vp_j}\right)^2 .
\end{equation}

\begin{233Lemma}
$N$ is invariant under complex rotations of the coordinates.
\end{233Lemma}

\begin{proof}
$U=(u_{jk}) \in U(n+1)$, that is  $U\oU^t = I$, and set $w=Uz$.
Then
\addtocounter{equation}{1}
\begin{equation}
\label{eq2.34}
w_m = \sum_{\ell=1}^{n+1} u_{m\ell} z_\ell \ \Rightarrow\ \frac{\p}{\p z_\ell} = \sum_{m=1}^{n+1} u_{m\ell}
\frac{\p}{\p w_m} ,
\end{equation}
and
\begin{equation}
\label{eq2.35}
z^t \p_z = w^t (U^t)^{-1} U^t \p_w  = w^t \p_w .
\qedhere
\end{equation}
\end{proof}

Consequently, $\oN$, $N^*$ and $\oN^*$ are also invariant under such rotations and so are $\Ree N^* N$
and $\Imm N^* N$.
Also, $\p_z\cdot \p_{\oz} = (\p_w U ) \cdot (\oU\vphantom{U}^t \p_{\ow}) = \p_w \cdot \p_{\ow}$, so $\De$
of (\ref{eq1.2}) is also invariant under complex rotations of the coordinates.
In view of (\ref{eq2.22}) so is $\De^C$, hence its restriction
to $S^{2n+1}$, $\De_C$, is also invariant under complex rotations of the
coordinates.
Note that
\begin{equation}
\label{eq2.36}
\sum_{j=1}^{n+1} \frac{\p}{\p \vp_j} = -2r \Imm N ,
\end{equation}
see (\ref{eq2.28}).
The invariance of $\De_C$ under complex rotations implies a similar invariance for its heat kernel
$P_C (t,w,z)$, $w,z\in S^{2n+1}$.
We choose coordinates in which $w$ is represented by $(1,0,\ldots,0)$.
The complex rotation which achieves this is given by any unitary matrix $U$ whose first
row is $\ow = (\ow_1,\ldots,\ow_{n+1})$.
Such a $U$ sends the point $z$ to
\begin{equation}
\label{eq2.37}
Uz = (z\cdot \ow , z'_2 ,\ldots, z'_{n+1}),\qquad z\in S^{2n+1} ,
\end{equation}
with $|z'_2|^2 + \ldots + |z'_{n+1}|^2 = 1 - |z\cdot \ow|^2$.
Again let $U'=U(n)$ denote a complex rotation in $\BC^n$ which sends
$z'=(z'_2,\ldots, z'_{n+1})$ to $(\sqrt{1-|z\cdot \ow|^2} , 0,\ldots,0)$ and
let $1$ represent the identity operator on $\BC$.
Then $1\oplus U' \in U(n+1)$ and the composition $(1\oplus U')\circ U$ sends $w$ to $(1,0,\ldots,0)$ and $z$
to $(z\cdot\ow , \sqrt{1-|z\cdot \ow|^2},0,\ldots,0)$.
Consequently, the heat kernel $P_C(t)$ of $\De_C$ is a function of the real variables $|z\cdot \ow|$ and $\arg (z\cdot \ow)$
only, or, in spherical coordinates $P_C(t)=P_C (t,\cos\vt_1 , \vp_1)$ is a function of $\vt_1$ and $\vp_1$ only when $w=(1,0,\ldots,0)$.
This can also be seen from (\ref{eq2.32}) in view of \eqref{eq2.9}, \eqref{eq2.10}.
Let $\cL_C$ denote $\De_C$ acting on functions of $\vt_1$ and $\vp_1$ only,
\begin{equation}
\label{eq2.38}
\cL_C = \frac12 \frac{\p^2}{\p\vt_1^2} + \big( (n-1) \cot\vt_1 + \cot 2\vt_1\big) \frac{\p}{\p\vt_1} +
\frac12 (\tan^2 \vt_1) \frac{\p^2}{\p\vp_1^2} .
\end{equation}
Then (\ref{eq2.32}) implies that
\begin{equation}
\label{eq2.39}
\De_C = \cL_C + \ldots
\end{equation}
where $\ldots$ denotes terms with a $\p/\p\vt_j$, $j=2,\ldots,n$, or with a $\p/\p\vp_k$, $k=2,\ldots, n+1$.
Therefore the fundamental solution for $\cL_C - \p/\p t$ with base point $(0,0)$ is also a solution
of $\De_C - \p/\p t$ and normalized on $S^{2n+1}$ yields the heat kernel for $\De_C$.
Analogously we write
\begin{equation}
\label{eq2.40}
\cL_S = \frac12 \frac{\p^2}{\p\vt_1^2} + \big( (n-1) \cot\vt_1 + \cot 2\vt_1\big) \frac{\p}{\p\vt_1} +
\frac12 \frac{1}{\cos^2 \vt_1} \frac{\p^2}{\p\vp_1^2} ,
\end{equation}
which agrees with \eqref{eq1.17}.
$\cL_S=\De_{S,1}$ when acting on functions of $\vt_1$ and $\vp_1$ only (see \eqref{eq2.9}).
In particular, the heat kernel of $\cL_S$ with base point $(0,0)$ is also the heat kernel of $\De_S$
when normalized on $S^{2n+1}$.
$\De_S$ is invariant with respect to the full orthogonal group of rotations $O(2n+2)$ hence
so is its heat kernel $P_S (t,y,x)$, $x,y\in \BR^{2n+2}$, $|x|=|y|=1$. 
Rotating $y$ to $(1,0,\ldots,0)$ shows that $P_S (t,y,x)$ is a function of $x_1y_1 +\ldots+x_{2n+2}
y_{2n+2}$ only.
Assuming $y=(1,0,\ldots,0)$ we see that $P_S(t)$ is a function of $x_1 = \cos\vt_1 \cos\vp_1$ only and
\begin{equation}
\label{eq2.41}
\cL_S = \frac12 (1-x_1^2) \frac{d^2}{dx_1^2} - \Big( n + \frac12\Big) x_1 \frac{d}{dx_1} .
\end{equation}
We shall need to integrate functions $f(x_1)=f(\cos\vt_1\cos\vp_1)$ on $S^{2n+1}$.
Set $v=\cos\vt_1 \sin\vp_1$. 
Then (\ref{eq2.12}) yields
\begin{align}
\label{eq2.42}
\int_{S^{2n+1}} f(x_1) dS
&= \pi^n \int f(x_1) \left( \prod_{j=1}^n \sin^{2(n-j)} \vt_j \sin 2\vt_j d\vt_j\right) d\vp \\
&= \frac{\pi^n}{\Gamma (n)} \int f(x_1) \sin^{2(n-1)} \vt_1 \sin 2\vt_1 d\vt_1 d\vp_1 \nonumber \\
&= \frac{2\pi^n}{\Gamma (n)} \int_{x_1^2 + v^2 \leq 1} f(x_1) (1-x_1^2 - v^2)^{n-1} dx_1 dv \nonumber \\
&= \frac{2\pi^{n+1/2}}{\Gamma(n+1/2)} \int_{-1}^1 f(x_1) (1-x_1^2)^{n-1/2} dx_1 .\nonumber
\end{align}

\section{The heat kernel of $\De_S$}
\label{sec3}

A second order elliptic differential operator
\begin{equation}
\label{eq3.1}
\De = \frac12 \sum_{j=1}^n X_j^2 + \cdots
\end{equation}
has a local heat kernel $P$ which can be expressed in the following form:
\begin{equation}
\label{eq3.2}
P = \frac{1}{(2\pi t)^{n/2}} e^{-\frac{d(x,x')^2}{2 t}} (a_0 + a_1 t + \cdots ) ;
\end{equation}
here $d(x,x')$ denotes the Riemannian distance between $x$ and $x'$ when the metric is
induced by the orthonormal basis $X_1,\ldots, X_n$ of $TM_n$, and the $a_j$'s are smooth
functions of $x,x'$ near the diagonal.
The coefficient of the term with the highest singularity in $t$ in $\p P /\p t - \De P$ is
\begin{equation}
\label{eq3.3}
- \frac{\p}{\p t} \frac{d(x,x')^2}{2 t} - \frac12 \sum_{j=1}^n \left( X_j \frac{d(x,x')^2}{2 t}\right)^2 = 0 .
\end{equation}
Thus $d^2 / 2 t$ is a solution of the Hamilton-Jacobi
equation.

We shall carry out this procedure in the case of $\De_S$ of (\ref{eq2.9}); the result is known,
but the technique suggests an approach to the construction of the heat kernel $P_C$ of the subelliptic
Laplacian $\De_C$ of (\ref{eq2.32}) which we shall carry out in \S\ref{sec4}.
Let $S$ denote a solution of the Hamilton-Jacobi equation induced by $\De_S$,
\begin{align}
\label{eq3.4}
0 = \frac{\p S}{\p t} &+ \frac12 \left\{ \left( \frac{\p S}{\p \vt_1}\right)^2 +  \frac{(\p S / \p\vp_1)^2}{\cos^2 \vt_1} + \sum_{j=2}^n
\frac{(\p S / \p\vt_j)^2}{\sin^2\vt_1\ldots \sin^2\vt_{j-1}} \right . \\
&\left. + \sum_{j=2}^n
\frac{(\p S / \p\vp_j)^2}{\sin^2\vt_1\ldots \sin^2\vt_{j-1}\cos^2\vt_j} +
\frac{(\p S / \p\vp_{n+1})^2}{\sin^2\vt_1\ldots \sin^2 \vt_n}\right\} .\nonumber
\end{align}
This is a first order partial differential equation whose solution can be found in the form
of an action integral $S$ along the bicharacteristics.
To this end we set
\begin{equation}
\label{eq3.5}
\omega = \nabla_\vt S ,\qquad \tau = \nabla_\vp S , \qquad \zeta =\frac {\p S}{\p t} ,
\end{equation}
\begin{align}
\label{eq3.6}
&H(\vt,\omega,\tau) = \frac12 \left\{ \omega_1^2 + \frac{\tau_1^2}{\cos^2\vt_1} + \sum_{j=2}^n
\frac{\omega_j^2}{\sin^2\vt_1 \ldots \sin^2 \vt_{j-1}} \right.\\
&\left.  + \sum_{j=2}^n \frac{\tau_j^2}{\sin^2\vt_1 \ldots \sin^2 \vt_{j-1}
\cos^2 \vt_j} + \frac{\tau_{n+1}^2}{\sin^2\vt_1\ldots \sin^2 \vt_n}\right\} . \nonumber 
\end{align}
$H$ is the Hamiltonian and (\ref{eq3.4}) is $\zeta + H(\vt,\om,\tau) = 0$.
The bicharacteristic curves $\big( \vt(s), \vp(s), t(s), S(s), \om(s),\tau(s),\zeta(s)\big)$, $0\leq s\leq t$, are solutions of
\begin{align}
&\dot\vt = \nabla_\om H ,\qquad \dot\vp = \nabla_\tau H ,
\label{eq3.7} \\
&\dot\om = - \nabla_\vt H ,\qquad \dot\tau = -\nabla_\vp H = 0 ,
\label{eq3.8} \\
&\dot t = 1,\qquad \dot\zeta = 0 ,
\label{eq3.9} \\
&\dot S = \om \cdot \dot\vt + \tau\cdot \dot\vp - H ,\qquad \cdot \Doteq d/ds ,
\label{eq3.10}
\end{align}
since $H=H(s) = -\zeta =$ constant along the bicharacteristic.
(3.8) yields 
\begin{equation}
\label{eq3.11}
\tau_j (s) = \tau_j = {\mathrm {constant}},\qquad j=1,2,\ldots, n+1 .
\end{equation}
(3.7) and (3.8) represent a system of $4n+2$ first order ordinary differential
equations for $4n+2$ unknown functions.
To fix the solution we need $4n+2$ conditions.
Start with
\begin{equation}
\label{eq3.12}
\vt_1(0) = 0 ,\qquad \vt(t) = \vt ,\qquad \vp (t) = \vp .
\end{equation}
Substituting this into \eqref{eq3.6} implies
\begin{equation}
\label{eq3.13}
\omega_j (0) = 0,\quad j=2,\ldots, n ,\qquad
\tau_j (0) = 0 ,\quad j=2,\ldots, n+1 .
\end{equation}
\eqref{eq3.12}, \eqref{eq3.13} give $4n+1$ conditions.
For the last one needed to fix the bicharacteristic we keep $\tau_1$ as an arbitrary parameter,
at least for now.
Also,
\begin{equation}
\label{eq3.14}
\tau_j (s) = \tau_j = 0 ,\qquad j = 2,\ldots, n+1 ,
\end{equation}
in view of \eqref{eq3.11}, \eqref{eq3.13}, and therefore \eqref{eq3.6}, (3.7) yield
\begin{equation}
\label{eq3.15}
\dot\vp_j = 0 \ \Rightarrow\ \vp_j (s) = \vp_j ,\qquad
j = 2 , \ldots, n+1.
\end{equation}
 From \eqref{eq3.14} one has
\begin{equation}
\label{eq3.16}
\dot\om_n = -H_{\vt_n} = a_n \tau_n^2 + a_{n+1} \tau_{n+1}^2 = 0 ,
\end{equation}
and \eqref{eq3.13} implies
\begin{equation}
\label{eq3.17}
\om_n (s) = \om_n (0) = 0 .
\end{equation}
When $k < n$ one finds
\begin{equation}
\label{eq3.18}
\dot\om_k = -H_{\vt_k} = \sum_{j=k}^{n+1} a_j \tau_j^2 + \sum_{j=k+1}^n b_j \om_j^2 .
\end{equation}
In particular $\dot\om_{n-1} =0$. 
Continuing in this manner one obtains $\dot\om_k=0$, $k=2,\ldots,n$, and therefore
$\om_k (s) = \om_k (0) = 0$, $k=2,\ldots, n$, so one has
\begin{equation}
\label{eq3.19}
\dot\vt_k = H_{\om_k} = 0 \ \Rightarrow\ \vt_k (s) = \vt_k ,\qquad k=2,\ldots, n .
\end{equation}
We have derived

\begin{320Lemma}
The boundary conditions \eqref{eq3.12}, \eqref{eq3.13} plus a choice of the parameter
$\tau_1\in\BR$ fix the following parts of the bicharacteristic curve:
\addtocounter{equation}{1}
\begin{equation}
\label{eq3.21}
\tau_1(s) = \tau_1 ,
\end{equation}
\begin{equation}
\label{eq3.22}
\vt_j (s) = \vt_j ,\qquad \om_j (s) = 0 ,\qquad j=2,\ldots, n,
\end{equation}
\begin{equation}
\label{eq3.23}
\vp_j(s) = \vp_j ,\qquad \tau_j(s) = 0 , \qquad j=2,\ldots, n+1 .
\end{equation}
\end{320Lemma}

Substituting \eqref{eq3.21}--\eqref{eq3.23} into \eqref{eq3.6} yields

\begin{324Corollary}
$\vt_1(s)$, $\vp_1(s)$, $\om_1(s)$ and $\tau_1$ represent the bicharacteristic
curve induced by the Hamiltonian
\addtocounter{equation}{1}
\begin{equation}
\label{eq3.25}
H = \frac12 \left( \om_1^2 + \frac{\tau_1^2}{\cos^2\vt_1}\right)
\end{equation}
of $\cL_S$; in particular,
\begin{equation}
\label{eq3.26}
\dot\vt_1 = \om_1 = \pm\sqrt{2H - \frac{\tau_1^2}{\cos^2\vt_1}} ,
\end{equation}
\begin{equation}
\label{eq3.27}
\dot\vp_1 = \frac{\tau_1}{\cos^2\vt_1} .
\end{equation}
\end{324Corollary}

\begin{328Lemma}
Let $E^2 = 2H$.
Then the $\vt_1(s)$ component of the bicharacteristic is given by
\addtocounter{equation}{1}
\begin{equation}
\label{eq3.29}
\sin^2 \vt_1 (s) = \left( 1 - \frac{\tau_1^2}{E^2}\right) \sin^2 Es,\qquad
0\leq s \leq t .
\end{equation}
\end{328Lemma}

\begin{proof}
We set $k^2=2H/\tau_1^2$ and write \eqref{eq3.26} as
\begin{align}
\pm \frac{d\vt_1}{ds} &= \sqrt{2H - \frac{2\tau_1^2}{1+\cos 2\vt_1}} \\
&= \frac{|\tau_1|k}{\sin 2\vt_1} \frac{k^2-1}{k^2} 
\sqrt{1 - \left( \frac{k^2}{k^2-1} \cos 2\vt_1 - \frac{1}{k^2-1}\right)^2} ,\nonumber
\end{align}
or,
\begin{equation*}
\frac{d\left( \frac{k^2}{k^2-1} \cos 2\vt_1 - \frac{1}{k^2-1}\right)}
{\sqrt {1 - \left( \frac{k^2}{k^2-1} \cos 2\vt_1 - \frac{1}{k^2-1}\right)^2} }= \pm 2 |\tau_1|k \, ds ,
\end{equation*}
which yields
\begin{equation}
\label{eq3.31}
\frac{k^2}{k^2-1} \cos 2\vt_1 - \frac{1}{k^2-1} = \sin \big( \pm 2 |\tau_1| k (s-s_0)\big) ,
\end{equation}
with some $s_0$ still to be determined.
At $s=0$ \eqref{eq3.31} is $1=\sin (\pm 2 |\tau_1| k s_0)$, and expanding the sine function
in \eqref{eq3.31} via the addition formula gives
\begin{equation*}
\frac{k^2}{k^2-1} \cos 2\vt_1 - \frac{1}{k^2-1} = \cos \big( 2|\tau_1| k s\big) .
\end{equation*}
Subtracting $1$ from each side leads to \eqref{eq3.29}.
\end{proof}

Next, \eqref{eq3.27} and
\begin{equation}
\label{eq3.32}
\cos^2 \vt_1 (s) = \cos^2 E s + \frac{\tau_1^2}{E^2} \sin^2 Es = \frac
{1 + \frac{\tau_1^2}{E^2} \tan^2 Es}{1+\tan^2 Es}
\end{equation}
imply
\begin{equation}
\label{eq3.33}
d\vp_1 = \frac{d\left( \frac{\tau_1}{E} \tan Es\right)}{1+\left( \frac{\tau_1}{E}\tan Es\right)^2} ,
\end{equation}
and we have

\begin{334Lemma}
Assuming $\vp_1 (0) = 0$, the $\vp_1$-component of the bicharacteristic of $\cL_S (\De_S)$ is given by
\addtocounter{equation}{1}
\begin{equation}
\label{eq3.35}
\vp_1 (s) = \tan^{-1} \left( \frac{\tau_1}{E} \tan Es\right)  ,
\end{equation}
which can be continued for all $s>0$.
\end{334Lemma}

The parameters $E$ and $\tau_1$ are to be chosen so that $\vt_1(t) = \vt_1$, $\vp_1 (t) = \vp_1$.
In particular, \eqref{eq3.29} and \eqref{eq3.35} imply
\begin{equation*}
\frac{\sin^2\vt_1}{\sin^2 Et} = 1 - \left( \frac{\tau_1}{E}\right)^2 = 1 - \frac{\tan^2\vp_1}{\tan^2 Et} ,
\end{equation*}
\begin{equation*}
\sin^2 Et - (\cos^2 Et) \tan^2 \vp_1 = \sin^2 \vt_1 .
\end{equation*}
Subtracting $1$ from each side yields $\cos Et  = \cos\vt_1 \cos\vp_1$, hence
\begin{equation}
\label{eq3.36}
Et = \cos^{-1} (\cos\vt_1 \cos\vp_1) + 2k \pi = \g + 2k\pi ,\quad k\in Z , 
\end{equation}
where we set $\cos \g = \cos\vt_1 \cos \vp_1$; the angle $\g$ subtends the arc from $(1,0,\ldots, 0)$
to $z$. (3.10) gives
$\dot S = H = E^2 / 2$, so one has
\begin{equation}
\label{eq3.37}
S = Ht = \frac
{\big( \cos^{-1} (\cos\vt_1 \cos\vp_1)+2k\pi\big)^2}{ 2t} = \frac{(\g + 2k\pi)^2}{2t} , \quad k\in Z .
\end{equation}
$S$ satisfies the Hamilton-Jacobi equation \eqref{eq3.4}.
Then \eqref{eq3.2}, \eqref{eq3.3} and \eqref{eq3.37} suggest the following result:

\begin{338Theorem}
Given $z,w\in S^{2n+1}$, let $\gamma$ denote the angle which subtends the arc that joins
$z$ and $w$ on a great circle, $0\leq\gamma < \pi$.
Then the heat kernel $P_S$ of $\De_S$ on $S^{2n+1}$ is given by
\addtocounter{equation}{1}
\begin{equation}
\label{eq3.39}
P_S = \frac{e^{\frac{n^2}{2}t}}{(2\pi t)^{n+1/2}} \sum_{k=-\infty}^\infty e^{-\frac{(\g+2k\pi)^2}{2t}}
v_n (\g + 2k\pi , t) ,
\end{equation}
where
\begin{equation}
\label{eq3.40}
v_n (\g , t) = v^n w_n = \left( \frac{\g}{\sin\g}\right)^n \left(\sum_{\ell=0}^{n-1} w_{n,\ell} (\g) t^\ell
\right) ,
\end{equation}
with $w_{n,0}=1$ and the $w_{n,\ell}$, $\ell=1,2,\ldots, n-1$ are found by iteration,
\begin{equation}
\label{eq3.41}
w_{n,\ell}(\g) = \frac{1}{\g^\ell} \int_0^\g v^{-n} \sigma^{\ell-1} \Big( \cL_S - \frac{n^2}{2}\Big) v_{n,\ell-1} d\sigma ,
\end{equation}
$v_{n,\ell} (\g) = v(\g)^n w_{n,\ell} (\g)$.
\end{338Theorem}

The proof will be given in several steps.

\begin{342Lemma}
\eqref{eq3.40} and \eqref{eq3.41} follow from \eqref{eq3.39}.
\end{342Lemma}

\addtocounter{equation}{1}

\begin{proof}
Set $x_1 = \cos\g$.
Then \eqref{eq2.41} yields 
\begin{equation}
\label{eq3.43}
\cL_S^{(n)} = \frac12 \frac{d^2}{d\g^2} + n(\cot\g)\frac{d}{d\g} ;
\end{equation}
we let $\De_S ^{(n)}$ represent $\De_S$ acting on $C^\infty (S^{2n+1})$ and let $\cL_S^{(n)}$
denote the reduced operator.
With $v'= \p v / \p \g $, one has
\begin{align}
\label{eq3.44}
0 = &\left( \frac{\p}{\p t} - \cL_S^{(n)}\right) \frac{e^{A_nt} e^{-\frac{\g^2}{2t}} v_n} {t^{n+1/2}} \cdot \\
&\cdot \left\{ \frac{ 2\g v'_n - 2n\g \Big( \frac{1}{\g} -\cot \g\Big) v_n}
{2t} - (\cL_S^{(n)} - A_n) v_n + \frac{\p v_n}{\p t}\right\} .\nonumber
\end{align}
Expanding in powers of $t$ gives
\begin{equation}
\label{eq3.45}
v'_{n,0} = n \left( \frac{1}{\g} - \cot\g \right) v_{n,0} ,\quad \text{i.e.}\ 
v_{n,0} = \left( \frac{\g}{\sin\g}\right)^n = v^n ,
\end{equation}
and then
\begin{equation}
\label{eq3.46}
w_{n,\ell} = \frac{1}{\g^\ell} \int_0^\g \frac{\sigma^{\ell-1} (\cL_S^{(n)} - A_n) v_{n,\ell-1}}{v_{n,0}} d\sigma .
\end{equation}
To fix $A_n$ we calculate $w_{n,1}$:
\begin{equation}
\label{eq3.47}
w_{n,1} = \frac{1}{\g} \int_0^\g \frac{\cL_S^{(n)} v^n}{v^n} - A_n ,
\end{equation}
\begin{align*}
\int_0^\g \frac{\cL_S^{(n)} v^n}{v^n} &= \frac12 \int_0^\g v^{-n}
\left( \frac{d}{d\g} + 2n \cot \g\right) \frac{dv^n}{d\g} \\
&= v^{-n} \frac{dv^n}{2d\g} + \frac{n}{2} \int_0^\g v^{-n} \left( \frac{d\log v}{d\g} + 2\cot\g\right) \frac{dv^n}{d\g}\\
&= n \frac{d\log v}{2d\g} + \frac{n^2}{2} \int_0^\g \left( \frac{1}{\g^2} - \cot^2 \g\right) \\
&= -n (n-1) \frac{d\log v}{2d\g} + \frac{n^2}{2} \g ,
\end{align*}
hence
\begin{equation}
\label{eq3.48}
w_{n,1} = - \frac{n(n-1)}{2} \frac{d\log v}{d(\g^2/2)} ,
\end{equation}
if we choose $A_n = n^2/2$, the largest positive eigenvalue of $\cL_S^{(n)}$; one eigenfunction is
$(1-x_1^2)^{1/2-n} P_{n-1}^{(1/2-n,1/2-n)} (x_1)$, where $P_k^{(\al,\beta)}$ is a Jacobi polynomial.
In particular, \eqref{eq3.48} and \eqref{eq3.46} imply that $v_{1,\ell}=0$, $\ell=1,2,\ldots$
\end{proof}

\noindent That $v_{n,\ell}=0$, $\ell=n,n+1,\ldots$ for all $n$ will be shown later.

\begin{349Lemma}
\eqref{eq3.39} yields the heat kernel $P_S^{(1)}$ of $\De_S$ on $S^3$; i.e. Theorem 3.38 holds when $n=1$.
\end{349Lemma}
\addtocounter{equation}{1}

\begin{proof}
We claim that 1) $P_S^{(1)}$ is a smooth function when $t>0$, and 2) $P_S^{(1)} \to 0$ as $t\to 0$ uniformly
on sets which are bounded away from $\g=0$.
Assume $\g\ne \pi$.
Adding the $k$-th and $(-k)$-th term in \eqref{eq3.39} one obtains
\begin{align}
\label{eq3.50}
P_S^{(1)} &=\frac{2e^{t/2}}{(2\pi t)^{3/2}} \frac{\g}{\sin \g} e^{-\frac{\g^2}{2t}}\cdot\\
&\qquad  \cdot \sum_{k=0}^\infty e^
{-\frac{4k^2\pi^2}{2t}}  \left( \cosh \frac{4k\pi\g}{2t} - \frac{8k^2\pi^2}{2t} \frac{\sinh\frac{4k\pi\g}{2t}}{\frac{4k\pi\g}{2t}}
\right)\nonumber\\
&= O \left( \frac{e^{-\frac{\g^2}{2t}}}{(2\pi t)^{3/2}} \frac{\g}{\sin\g} \left\{ 1 + \sum_{k=1}^\infty
\frac{4k^2\pi^2}{2t} e^{-\frac{4k^2\pi^2}{2t}}\right\}\right) ;  
\nonumber
\end{align}
the second equality makes sense if $\g$ is bounded away from $\pi$.
Note our use of $1\leq (\sinh x) / x \leq \cosh x \leq e^x$.
The last sum is $O(1)$ by comparison with an integral.
The first equality shows that $P_S^{(1)}$ is a smooth function when $t>0$ and $\g\ne \pi$, and the last estimate
implies that $P_S^{(1)}\to 0$ when $t\to 0$, uniformly on sets which are bounded away from $\g=0,\pi$.
Near $\g=\pi$ we replace $\g$ by $\pi-\ve$:
\begin{align}
\label{eq3.51}
P_S^{(1)} &=
\frac{e^{t/2}}{(2\pi t)^{3/2} \sin(\pi-\ve)} \sum_{k=-\infty}^{\infty} (\pi-\ve + 2k\pi) e^{-
\frac{(\pi-\ve+2k\pi)^2}{2t}} \\
&= \frac{2e^{t/2}}{(2\pi t)^{3/2}} \frac{\ve}{\sin\ve} e^{-\frac{(\pi-\ve)^2}{2t}} e^{-\frac{2\pi \ve}{2t}} \sum_{k=0}^\infty
e^{-\frac{4k(k+1)\pi^2}{2t}} \cdot \nonumber \\
&\quad \cdot \left\{ - \cosh \frac{2(2k+1)\pi\ve}{2t} + \frac{(2k+1)^2\pi^2}{2t} \frac{\sinh\frac{2(2k+1)\pi\ve}{2t}}
{\frac{2(2k+1)\pi\ve}{2t}}\right\}\nonumber \\
&=\frac{e^{-\frac{(\pi-\ve)^2}{2t}}}{(2\pi t)^{3/2}} e^{t/2} O\left(
\frac1t + \sum_{k=1}^\infty \frac{4k^2\pi^2}{2t} e^{-\frac{4k^2\pi^2}{2t}}\right) ,\nonumber
\end{align}
where we added the terms with $k$ and $-k-1$ and then summed over $k$.
Again the final sum is $O(1)$.
Thus $P_S^{(1)}$ is smooth near $\g=\pi$ when $t>0$ and converges to $0$ when $t\to 0$, uniformly.
This justifies our claim.
Note the extra factor $t^{-1}$, a consequence of $\g=\pi$ being a conjugate point of $\g=0$, parametrized
by $S^2$.
To complete the proof we need to show that $\lim\limits_{t\to 0} \int\limits_{S^3} P_S^{(1)} = 1$.
To that end write $P_S^{(1)}=e^{t/2} (2\pi t)^{-3/2} f(\g)$ and note that \eqref{eq2.42} implies
\begin{align}
\label{eq3.52}
\int_{S^3} f(\g) &= 4\pi \sum_{k=-\infty}^\infty \int_0^\pi (\g + 2k\pi) e^{-\frac{(\g+2k\pi)^2}{2t}} \sin\g \, d\g\\
&= 4\pi \left\{ \sum_{k=0}^\infty \int_{2k\pi}^{(2k+1)\pi} + \sum_{k=1}^\infty \int_{-2k\pi}^{-(2k-1)\pi}\right\}
\g(\sin \g) e^{-\frac{\g^2}{2t}} d\g \nonumber \\
&= 4\pi \int_0^\infty \g (\sin \g) e^{-\frac{\g^2}{2t}} d\g \nonumber \\
&= (2\pi t)^{3/2} e^{-t/2} \nonumber
\end{align}
in view of Lemma 3.53.
Thus $\int\limits_{S^3} P_S^{(1)} =1$.
\end{proof}

\begin{353Lemma}
Given $\om >0$ and $\al>0$, one has
\addtocounter{equation}{1}
\begin{equation}
\label{eq3.54}
I_\al = \int_{\om-i\infty}^{\om+i\infty} e^{\la / t - \al \sqrt{2\la}} d\la = i\sqrt{2\pi} \al t^{3/2} e^{-\al^2 t/2} .
\end{equation}
\end{353Lemma}

\begin{proof}
The integrand vanishes exponentially for large $|\la|$.
We replace $\la / t$ by $\eta^2$,
\[
I_\al = t \int_{-i\infty}^{i\infty} \frac{d}{d\eta} (e^{\eta^2}) e^{-\al \sqrt{2t} \eta} d\eta 
= \al \sqrt 2 t^{3/2} e^{-\al^2 t/2} \int_{-i\infty}^{i\infty} e^{u^2} du ,
\]
after one integrates by parts and completes the square in the exponent.
\end{proof}

Thus we have completed the derivation of \eqref{eq3.39} when $n=1$.
For $n>1$ the proof will be carried out on an integral representation of the
right hand side of \eqref{eq3.39}.
Note that a residue expansion yields
\begin{equation}
\label{eq3.55}
P_S^{(1)} = \frac{e^{\frac12 t}}{(2\pi t)^{3/2}} \frac{1}{2\pi i} \int_{\omega - i\infty}^{\omega+i\infty}
\frac{e^{\la / t} d\la }{\cosh \sqrt{2\la} - \cos\g} ,\quad \omega > 0 .
\end{equation}
$P_S^{(n)}$, the heat kernel of 
$\Delta_S^{(n)}$, may be derived from $P_S^{(1)}$ as follows.
Let $V$ denote a solution of $\cL_S^{(n)} V = \p V / \p t$, that is
\begin{equation*}
\frac12 (1-x_1^2) \frac{\p^2 V}{\p x_1^2} - \left( n + \frac12\right) x_1 \frac{\p V}{\p x_1} = \frac{\p V}{\p t} .
\end{equation*}
Differentiating the two sides with respect to $x_1$ yields
\begin{equation*}
\frac12 (1-x_1^2) \frac{\p^2 V_{x_1}}{\p x_1^2} - \left( n + 1 + \frac12\right) \frac{\p V_{x_1}}{\p x_1} =
e^{-(n+\frac12)t} \frac{\p}{\p t} e^{(n+\frac12)t} V_{x_1} ,
\end{equation*}
or,
\begin{equation}
\label{eq3.56}
\cL_S^{(n+1)} \big( e^{(n+\frac12)t} V_{x_1}\big) = \frac{\p}{\p t} e^{(n+\frac12)t} V_{x_1} .
\end{equation}
This is useful in the following form:
\begin{357Lemma}
Let $u_n$ denote a solution of $\cL_S^{(n)} u_n = \p u_n /\p t$.
Then
\addtocounter{equation}{1}
\begin{equation}
\label{eq3.58}
u_{n+1} = e^{\frac{2n+1}{2}t} \frac{\p u_n}{2\pi \p x_1} = e^{\frac{2n+1}{2} t} \frac{\p u_n}{2\pi\p (\cos\g)}
\end{equation}
is a solution of $\cL_S^{(n+1)} u_{n+1} = \p u_{n+1} / \p t$.
\end{357Lemma}

\eqref{eq3.58} suggests
\begin{equation}
\label{eq3.59}
P_S^{(n)} = \prod_{\ell =1}^{n-1} e^{\frac{2\ell +1}{2} t} \left(\frac{\p}{2 \pi \p x_1}\right)^{n-1} P_S^{(1)} ,\qquad n>1,
\end{equation}
which applied to \eqref{eq3.55} implies

\begin{360Lemma}
Given $\omega>0$, one has
\addtocounter{equation}{1}
\begin{equation}
\label{eq3.61}
P_S^{(n)} = \frac{\Gamma(n) e^{\frac{n^2}{2}t}}{(2\pi)^{n-1} (2\pi t)^{3/2}}
\frac{1}{2\pi i} \int_{\omega-i\infty}^{\omega+i\infty} \frac{e^{\la/t} d\la}
{(\cosh \sqrt {2\la} - \cos\g)^n} .
\end{equation}
\end{360Lemma}

\begin{proof}
Lemma 4.44 shows that the denominator of the integrand has its zeros on the nonpositive
real axis. 
Consequently, \eqref{eq3.61} is well defined as the integrand is summable along
the path of integration, uniformly in $\g$; note that the numerator is bounded
and the denominator increases exponentially.
This allows for the exchange of the integral for $P_S^{(1)}$ in \eqref{eq3.55} with the
derivations of \eqref{eq3.59} and therefore $P_S^{(n)}$ of \eqref{eq3.61} is a
solution of the heat equation $\De_S^{(n)} P_S^{(n)} = \p P_S^{(n)} / \p t$ in view of
Lemma 3.57.
Next we note that the largest zero of the denominator in the integrand of \eqref{eq3.61}
is $-\g^2 / 2$, $0\leq \g \leq\pi$, so if $\g\ne 0$, we may move the path of integration
to $\Ree \la = -\g^2 / 4$, along which the denominator of the integrand still increases
exponentially, but the numerator is $O(e^{-\g^2/4t})$.
Consequently, $P_S^{(n)}$ of \eqref{eq3.61} vanishes exponentially as $t\to 0$, uniformly
in sets where $\g$ is bounded away from zero.

We claim that $\int\limits_{S^{2n+1}} P_S^{(n)} dS = 1$.
To this end set
\begin{equation}
\label{eq3.62}
\psi_n = \int_{\omega-i\infty}^{\omega+i\infty} \frac
{e^{\la /t} d\la} {(\cosh \sqrt{2\la} - \cos\vt_1 \cos\vp_1)^n} ,\qquad \omega > 0 .
\end{equation}
In view of \eqref{eq2.42},
\begin{equation*}
\int_{S^{2n+1}} \psi_n dS = \frac{\pi^n}{\Gamma(n)} \int_{\omega-i\infty}^{\omega+i\infty}
e^{\la /t} d\la \int_{0}^{\pi/2} \int_{-\pi}^{\pi} \frac
{\sin^{2(n-1)} \vt_1 \sin 2\vt_1 d\vp_1 d\vt_1}{(\cosh\sqrt{2\la} - \cos\vt_1\cos\vp_1)^n} .
\end{equation*}
With $u=0$, \eqref{eq4.57}, \eqref{eq4.58} yield
\begin{equation}
\label{eq3.63}
\int_{-\pi}^{\pi} \frac{d\vp_1}{(A-\cos\vt_1 \cos\vp_1)^n} =
\frac{(-1)^{n-1}}{\Gamma(n)} \frac{d^{n-1}}{dA^{n-1}} 
\frac{2\pi}{\sqrt{A^2-\cos^2 \vt_1}} ,
\end{equation}
hence
\begin{align}
\label{eq3.64}
\int_{S^{2n+1}} \psi_n dS &= \frac{2(-1)^{n-1} \pi^{n+1}}{\Gamma(n)^2}
\int_{\omega-i\infty}^{\omega+i\infty} e^{\la/t} d\la \frac{d^{n-1}}{dA^{n-1}} \int_0^1
\frac{(1-u)^{n-1} du}{\sqrt{A^2-u}}\\
&= \frac{(2\pi)^{n+1}}{n\Gamma (n)} \int_{\omega-i\infty}^{\omega+i\infty} e^{\la/t - n\sqrt{2\la}} d\la \nonumber\\
&= \frac{(2\pi)^{n+1}}{\Gamma (n)} \sqrt{2\pi} it^{3/2} e^{-\frac{n^2}{2} t} ,\nonumber
\end{align}
where we used Lemmas 4.59 and 3.53 with $A=\cosh\sqrt{2\la}$.
Then \eqref{eq3.61}, \eqref{eq3.62} and \eqref{eq3.64} justify the claim and we have proved Lemma 3.60.
\end{proof}

The denominator of the integrand in \eqref{eq3.61} has zeros at $\la_k = -(\g + 2k\pi)^2 /2$, $k\in Z$.

\begin{365Lemma}
Let $\ve_k$ denote a small circular path of positive direction, radius $\ve \ll 1$ and center $\la_k$.
Then
\addtocounter{equation}{1}
\begin{equation}
\label{eq3.66}
P_S^{(n)} = \frac{\Gamma(n)e^{\frac{n^2}{2}t}}{(2\pi)^{n-1} (2\pi t)^{3/2}} \sum_{k=-\infty}^\infty
\frac{1}{2\pi i} \oint_{\ve_k} \frac{e^{\la/t} d\la}{(\cosh\sqrt{2\la} - \cos\g)^n} ,
\end{equation}
uniformly in $\g$ when $\g$ is bounded away from $\pi$.
\end{365Lemma}

\begin{proof}
One may move the line of integration in \eqref{eq3.61} to $\omega<0$, and compensate for the difference
by adding a finite number of $\ve_k$-integrals, as usual.
Then \eqref{eq3.66} follows if we show that \eqref{eq3.61} vanishes as $\omega\to -\infty$.
To carry out this argument we use the formula
\begin{align}
\label{eq3.67}
\sqrt{\xi+i\nu} &= \left( \frac{\sqrt{\xi^2+\nu^2}+\xi}{2}\right)^{1/2} + i (\sgn \nu) \left( \frac{\sqrt{\xi^2+\nu^2}-\xi}{2}\right)^{1/2}\\
&= \sigma + i\chi .\nonumber
\end{align}
In particular, the square root maps the line $\xi+i\nu$, $-\infty < \nu<\infty$, onto the hyperbola $\sigma^2-\chi^2=\xi$.
We work with $\xi<0$, so $\chi=\sqrt{\sigma^2-\xi} = \sqrt{\sigma^2 +|\xi|}$ is an increasing function of
$|\sigma|$ which is an increasing function of $|\nu|$ with $\nu=0$ iff $\sigma=0$.
Also,
\begin{align}
\label{eq3.68}
d^2 &= |\cosh (\sigma+i\chi)-\cos\g|^2 \\
&=\cosh^2\sigma - \sin^2 \chi - 2\cosh \sigma \cos\chi \cos\g + \cos^2\g . \nonumber 
\end{align}
With $2\la = \xi+i\nu$ we choose $\omega = \frac12 \xi=- \frac12 (2m+1)^2\pi^2$ in \eqref{eq3.61}; in particular, with bounded $\sigma$ and $m$ large one
has $\sin\chi \sim 0$ and
$\cos\chi \sim -1$.
Note that
\begin{equation}
\label{eq3.69}
\nu = 2\sigma \sqrt{\sigma^2 + (2m+1)^2 \pi^2} ,
\end{equation}
\begin{equation}
\label{eq3.70}
d\nu = 2 \left( \sqrt{\s^2 + (2m+1)^2\pi^2} + \frac{\s^2}{\sqrt{\s^2 + (2m+1)^2\pi^2}}\right) d\sigma ,
\end{equation}
and we need to estimate
\begin{align}
\label{eq3.71}
&\int_{-\frac12 (2m+1)^2\pi^2-i\infty}^{-\frac12 (2m+1)^2\pi^2+i\infty} \frac{e^{\la /t} d\la}{(\cosh\sqrt{2\la}-\cos\g)^n} \\
&= ie^{-\frac{(2m+1)^2\pi^2}{2t}} \int_{-\infty}^\infty \frac{e^{i\nu/2t}d\nu}{(\cosh \sqrt{2\la}-\cos\g)^n} ;\nonumber
\end{align}
it suffices to estimate the integral on $0<\nu<\infty$, since the other half of the integral is its complex conjugate.
Let $\sinh (\s_1 / 2)=1$ define $\sigma_1$.

(i) $\s > \s_1$.  From \eqref{eq3.68},
\begin{align*}
d^2 &\geq \cosh^2 \sigma - \sin^2 \chi - 2\cosh\s + \cos^2\g - 1 + 1\\
&= (\cosh \s - 1)^2 - \sin^2 \chi - \sin^2 \g \\
&\geq 2 \left( 2 \sinh^2 \frac{\s}{2} - 1\right) \\
&\geq 2 \sinh^2 \frac{\s}{2} ,
\end{align*}
so,
\begin{equation}
\label{eq3.72}
d >\sinh \frac{\s}{2} >  \frac14 e^{\s/2} ,\qquad \s>\s_1 .
\end{equation}

(ii) $\s < \s_1$.  Here $\sigma \sim 0$, so with large $m$ and $|\g| < \pi - \delta$,
$\delta > 0$, one has
\begin{equation}
\label{eq3.73}
d^2 > | \cos\chi - \cos\g |^2
> \frac12  | 1 + \cos\g |^2
\geq \frac12 | 1 - \cos\delta |^2
\geq c(\delta)\delta^4 ,
\end{equation}
$c(\delta)>0$.
Now (i), (ii) and $\nu'(\s) \leq 2 \big( 2\s + (2m+1)\pi\big) \leq 4(2m+1)\pi (1+\s)$ imply
\begin{equation}
\label{eq3.74}
\left| \int_{-\tfrac12(2m+1)^2\pi^2-i\infty}^{-\tfrac12(2m+1)^2\pi^2+i\infty} \frac{e^{\la/t} d\la}{(\cosh\sqrt{2\la}-\cos\g)^n}\right| 
\leq
C(2m+1)\pi e^{-\frac{(2m+1)^2\pi^2}{2t}} ,
\end{equation}
which vanishes as $m\to\infty$, uniformly in $\g$, $|\g |< \pi - \de$, $\de > 0$.

To complete the proof of Lemma 3.65 we note that $d$ increases exponentially as $|\nu|\to \infty$, uniformly
in $\xi$ in a finite interval.
In particular, for $\omega > 0$, $A>0$, one has
\begin{equation}
\label{eq3.75}
\lim_{|\nu|\to\infty} \int_\omega^{-A} \frac{e^{\la/t}d\la}{(\cosh \sqrt{2\la} - \cos\g)^n} = 0 ,
\end{equation}
and we have derived the expansion \eqref{eq3.66}.
\end{proof}

\begin{proof}[End of the proof of Theorem 3.38]
Evaluating \eqref{eq3.66} yields \eqref{eq3.39} with
\begin{equation}
\label{eq3.76}
v_n = e^{\frac{\g^2}{2t}} \left( \frac{t}{2\pi} \frac{\p}{\p (\cos\g)}\right)^{n-1} e^{-\frac{\g^2}{2t}}
\frac{\g}{\sin\g} .
\end{equation}
In particular,
\begin{equation}
\label{eq3.77}
w_{n,\ell} = 0,\quad \ell \ge n .
\end{equation}
This completes the proof of Theorem 3.38.
\end{proof}

\noindent {\bf 3.78 Remark.}
\eqref{eq3.77} is a consequence of $v_{n,n}=0$, in view of \eqref{eq3.41},
which also implies $v_{n,n}=0$ if $(\cL_S^{(n)}-n^2/2) v_{n,n-1}=0$.
For example, when $n=2$,
\addtocounter{equation}{1}
\begin{equation}
\label{eq3.79}
v_{2,1} = -\frac12 \frac{\p v^2}{\g\p\g} = \frac{-1}{1-x^2} + \frac{x\cos^{-1}x}{(1-x^2)^{3/2}} ,
\end{equation}
and an elementary calculation using \eqref{eq2.41} gives $\cL_S^{(2)} v_{2,1}=2v_{2,1}$.

As expected from \eqref{eq3.2}  and \eqref{eq3.39} one has

\begin{380Lemma}
The length of a geodesic which joins $(1,0,\ldots,0)$ to $(z_1,\ldots,z_{n+1})\in S^{2n+1}$ in time $t$ is
\addtocounter{equation}{1}
\begin{equation}
\label{eq3.81}
\ell = | \cos^{-1} (\cos\vt_1 \cos\vp_1) + 2k\pi | = | \g + 2k\pi | ,\qquad k=0,\pm 1,\pm 2,\ldots
\end{equation}
\end{380Lemma}

\begin{proof}
For a geodesic, that is, the projection of a bicharacteristic on $S^{2n+1}$, $\dot\vt_j=0$, $j=2,\ldots,n$,
and $\dot\vp_j=0$, $j=2,\ldots,n+1$.
Hence
\begin{align*}
&\dz_1 = e^{i\vp_1} (-\dot\vt_1\sin\vt_1 + i\dot\vp_1\cos\vt_1),\\
&\quad \cdots \\
&\dz_k = e^{i\vp_k} (\dot\vt_1 \cos\vt_1) \sin\vt_2 \ldots \sin\vt_{k-1} \cos\vt_{k},\qquad k=2,\ldots,n ,\\
&\quad \cdots \\
&\dz_{n+1} = e^{i\vp_{n+1}} (\dot\vt_1\cos\vt_1) \sin\vt_2 \ldots \sin\vt_n ,
\end{align*}
and one has
\begin{align*}
&|z_1|^2 = \dot\vt_1^2 \sin^2\vt_1 + \dot\vp_1^2\cos^2\vt_1 = \dot\vt_1^2 \sin^2\vt_1 +\tau_1^2/\cos^2\vt_1 ,\\
&|\dz_2|^2 + \cdots + |\dz_{n+1}|^2 = \dot\vt_1^2 \cos^2\vt_1 ,\\
&|\dz_1|^2 + \cdots + |\dz_{n+1}|^2 = \dot\vt_1^2 + \tau_1^2 / \cos^2\vt_1 = 2H .
\end{align*}
Consequently, the length of the geodesic is
\begin{equation*}
\ell = \int_0^t \sqrt {2H} = \sqrt{2H} t = |Et| ,
\end{equation*}
which proves Lemma 3.80 in view of \eqref{eq3.36}.
\end{proof}

We note that the geodesic joining $1=(1,0,\ldots,0)$ and $Q=(z_1,\ldots, z_{n+1})$ is in the plane generated
by the rays $\vec 1$ and $\vec Q$.
Let $Q(s)$ denote a point of the geodesic between $1$ and $Q$, $0<s<t$.
Then replacing $\vt_1(s)$ and $\vp_1(s)$ with the appropriate function of $\g s=Es$ one obtains the following formula:
\begin{equation*}
Q(s) = (\sin\g s) (\cot\g s - \cot\g t)\vec 1 + \frac{\sin\g s}{\sin \g t} \vec Q .
\end{equation*}
Since $\vp_1(s)$ is a strictly increasing function of $s$, the geodesics are as follows:
the shortest one with length $\g$ plus $k$ added great circles of length $|\g + 2k\pi|$, $k\in Z$,
as expected.

\section{The heat kernel for $\Delta_C$}
\label{sec4}

We shall follow the line of argument in \S\ref{sec3}, and look for a solution $S$ of the Hamilton-Jacobi
equation induced by $\p/\p t - \De_C$,
\begin{align}
\label{eq4.1}
0 = &\frac{\p S}{\p t} + \frac12 \Bigg\{ \left( \frac{\p S}{\p\vt_1}\right)^2 
+ \frac{(\p S/\p\vp_1)^2}{\cos^2 \vt_1} 
+ \sum_{j=2}^{n}
\frac{(\p S / \p \vt_j)^2}{\sin^2\vt_1 \ldots \sin^2 \vt_{j-1}}  \\
&+ \sum_{j=2}^n\, \frac{(\p S /\p\vp_j)^2}{\sin^2\vt_1\ldots\sin^2
\vt_{j-1}\cos^2\vt_j} + \frac{(\p S/\p\vp_{n+1})^2}{\sin^2\vt_1\ldots \sin^2\vt_n} \nonumber \\
&- \left( \frac{\p S}{\p\vp_1} + \cdots + \frac{\p S}{\p\vp_{n+1}}\right)^2\Bigg\} ,\nonumber
\end{align}
in the form of an action integral along the bicharacteristics of $\zeta +H$; here
\begin{align}
\label{eq4.2}
H(\vt,\omega,\tau) &= \frac12 \Bigg\{ \omega_1^2 
+ \frac{\tau_1^2}{\cos^2\vt_1} 
+ \sum_{j=2}^{n} 
\frac{\omega_j^2}{\sin^2\vt_1 \ldots \sin^2 \vt_{j-1}} \\
&+ \sum_{j=2}^n \frac{\tau_j^2}{\sin^2\vt_1\ldots\sin^2\vt_{j-1}\cos^2\vt_j} +
\frac{\tau_{n+1}^2}{\sin^2\vt_1\ldots\sin^2\vt_n} \nonumber\\
&-(\tau_1+\cdots + \tau_{n+1})^2\Bigg\} .\nonumber
\end{align}
The bicharacteristic curves are solutions of \eqref{eq3.7}--\eqref{eq3.10}, when $H$ is given by \eqref{eq4.2},
and \eqref{eq3.11} still holds.
We assume \eqref{eq3.12} and this again implies \eqref{eq3.13}.
The $4n+1$ conditions \eqref{eq3.12}, \eqref{eq3.13} plus the parameter $\tau_1$ fix the
bicharacteristic.
\eqref{eq3.15} is replaced by
\begin{equation}
\label{eq4.3}
\dot\vp_j = -\tau_1 ,\qquad j=2,\ldots, n+1 .
\end{equation}
The argument \eqref{eq3.16}--\eqref{eq3.19} still applies and yields $\omega_k (s) = \omega_k (0) = 0$,
$\vt_k (s) = \vt_k$, $k=2,\ldots,n$.
Thus we have derived

\begin{44Lemma}
The boundary conditions \eqref{eq3.12}, \eqref{eq3.13} plus a choice of the parameter
$\tau_1\in\BC$ fix the following parts of the bicharacteristic curve:
\addtocounter{equation}{1}
\begin{align}
\label{eq4.5}
&\tau_1(s) = \tau_1 , \\
\label{eq4.6}
&\vt_j (s) = \vt_j ,\ \omega_j(s) = 0 ,\ j=2,\ldots, n , \\
\label{eq4.7}
&\vp_j (s) = \tau_1 (t-s) + \vp_j ,\ \tau_j(s) =0 ,\ j=2,\ldots,n+1 .
\end{align}
\end{44Lemma}
\noindent Substituting \eqref{eq4.5}--\eqref{eq4.7} into \eqref{eq4.2} yields

\begin{48Corollary}
$\vt_1 (s)$, $\vp_1(s)$ and $\omega_1 (s)$ may be obtained from the reduced Hamiltonian
\addtocounter{equation}{1}
\begin{equation}
\label{eq4.9}
H = \frac12 (\omega_1^2 + \tau_1^2 \tan^2 \vt_1 ) ,
\end{equation}
as solutions of
\begin{equation}
\label{eq4.10}
\dot\vt_1 = \omega_1 = \pm \sqrt{2H - \tau_1^2 \tan^2 \vt_1}\, ,
\end{equation}
\begin{equation}
\label{eq4.11}
\dot\vp_1 = \tau_1 \tan^2\vt_1 .
\end{equation}
\end{48Corollary}
\noindent We note that \eqref{eq4.9} represents the Hamiltonian of $\cL_C$ of \eqref{eq1.17}.

\begin{412Lemma}
Set $2 H = E^2$, $\Omega^2 = E^2 + \tau_1^2$.
Then
\addtocounter{equation}{1}
\begin{equation}
\label{eq4.13}
\sin^2\vt_1 (s) = \frac{E^2}{\Omega^2} \, \sin^2\Omega s .
\end{equation}
\end{412Lemma}

\begin{proof}
\eqref{eq4.10} yields
\begin{align*}
\pm\frac{d\vt_1}{ds} &= \sqrt{2H + \tau_1^2 - \frac{2\tau_1^2}{1+\cos 2\vt_1}} \\
&= \frac{|\tau_1|}{\sin 2\vt_1} \frac{k^2}{\sqrt{k^2+1}} \sqrt{1 - \left( \frac{k^2+1}{k^2} \,\cos 2\vt_1-
\frac{1}{k^2}\right)^2}\, ,
\end{align*}
with $k^2 = 2H / \tau_1^2$.
Separating the variables and integrating one obtains
\begin{equation}
\label{eq4.14}
\frac{k^2+1}{k^2} \,\cos 2\vt_1 - \frac{1}{k^2} =
\sin \big( \pm 2 |\tau_1| \sqrt{k^2+1} (s-s_0)\big) ,
\end{equation}
with some $s_0$ still to be determined.
At $s=0$, \eqref{eq4.14} is $1=\sin (\pm 2 |\tau_1| \sqrt{k^2+1} s_0 )$, and expanding the sine
function in \eqref{eq4.14} via the addition formula gives
\begin{equation*}
\frac{k^2+1}{k^2} \, \cos 2\vt_1 - \frac{1}{k^2} = \cos \big( 2|\tau_1| \sqrt{k^2+1}\, s\big) ,
\end{equation*}
and subtracting $1$ from each side leads to \eqref{eq4.13}.
\end{proof}

\begin{415Lemma}
One has
\addtocounter{equation}{1}
\begin{equation}
\label{eq4.16}
\vp_1 (s) - \vp_1(0) = \tan^{-1} \Big( \frac{\tau_1}{\Omega} \, \tan\Omega s\Big) - \tau_1 s .
\end{equation}
\end{415Lemma}

\begin{proof}
In view of
\begin{equation}
\label{eq4.17}
1-\sin^2\vt_1 (s) = 1 - \frac{E^2}{\Omega^2} \, \sin^2\Omega s =
\frac{1 + \frac{\tau_1^2}{\Omega^2} \, \tan^2 \Omega s}{1+\tan^2 \Omega s} \, ,
\end{equation}
and \eqref{eq4.11}, one has
\begin{equation}
\label{eq4.18}
d\vp_1 = -\tau_1 ds + \frac{\tau_1 ds}{\cos^2 \vt_1 (s)} = -\tau_1 ds +
\frac{d\big( \frac{\tau_1}{\Omega} \,\tan\Omega s\big)}{1 + \big( \frac{\tau_1}{\Omega}\,\tan\Omega s\big)^2} ,
\end{equation}
which implies \eqref{eq4.16}.
\end{proof}

To find $S = Ht = (E^2 / 2)t$ we need $E$.
Assume $\vp_1 (0) = 0$.
Then at $s=t$,
\begin{equation}
\label{eq4.19}
\tan(\vp_1 + \tau_1 t) = \frac{\tau_1}{\Omega} \, \tan \Omega t ,
\end{equation}
and comparing with \eqref{eq4.13},
\begin{align*}
&\frac{\sin^2 \vt_1}{\sin^2 \Omega t} = 1 - \frac{\tau_1^2}{\Omega^2} = 1 -
\frac{\tan^2 (\vp_1 + \tau_1 t)}{\tan^2 \Omega t} , \\
&\sin^2\vt_1 = \sin^2\Omega t - \cos^2\Omega t\, \tan^2(\vp_1 +\tau_1 t) = 1 -
\frac{\cos^2\Omega t}{\cos^2 (\vp_1+\tau_1 t)} ,\\
&\cos\Omega t  =\pm \cos\vt_1 \cos (\vp_1 + \tau_1 t) ,\\
&\Omega t = \cos^{-1}\big( \cos\vt_1 \cos (\vp_1 + \tau_1 t)\big) + 2k\pi ,\qquad k\in\BZ .
\end{align*}
Squaring both sides yields
\begin{equation*}
(Et)^2 = -(\tau_1 t)^2 + \big( \cos^{-1} \big( \cos\vt_1 \cos (\vp_1 + \tau_1 t)\big) + 2k\pi \big)^2 ,
\end{equation*}
hence
\begin{equation}
\label{eq4.20}
S = Ht = \frac{-(\tau_1 t)^2 + \big(\cos^{-1} \big( \cos\vt_1 \cos(\vp_1+\tau_1 t)\big) +2k\pi\big)^2}{2t} .
\end{equation}
$\tau_1$ does not occur in the final formula of $P_C$ so we sum over $u=\tau_1 t$; 
$\tau_1$ is treated as a running parameter, even though $\vp_1(0)=0$ in \eqref{eq4.19} 
is supposed to fix it.
Set
\begin{equation*}
e(u) = -u^2 + \big( \cos^{-1} \big( \cos\vt_1 \cos (\vp_1 + u)\big)\big)^2 .
\end{equation*}
$S$ will be positive on the imaginary axis so we set $f(u)=e(-iu)$,
\begin{align}
\label{eq4.21}
f(u) &= u^2 + \big( \cos^{-1}\big( \cos\vt_1 \cos(\vp_1-iu)\big)\big)^2 = u^2 + \tgamma^2 \\
&= u^2 - \big(\cosh^{-1} \big( \cos\vt_1 \cosh (u+i\vp_1)\big)\big)^2 = u^2 -\kappa^2 ,\nonumber
\end{align}
where we set
\begin{equation}
\label{eq4.22}
\cos\tgamma = \cos\vt_1 \cos (\vp_1 - iu) = \cosh \kappa ,\qquad \kappa=i\tgamma .
\end{equation}

\begin{423Theorem}
The heat kernel $P_C$ of $\De_C$ on $S^{2n+1}$ is given by
\addtocounter{equation}{1}
\begin{equation}
\label{eq4.24}
P_C = \frac{e^{\frac{n^2}{2}t}}{(2\pi t)^{n+1}}\, \sum_{k\in\BZ}\, \int_{-\infty}^\infty e^
{-\frac{u^2-(\kappa + i 2k\pi)^2}{2t} }
\sum_{\ell=0}^{n-1} \, V_{n,\ell} (\kappa + i2k\pi) t^\ell du ,
\end{equation}
when $(\vt_1,\vp_1) \ne (0,  \pi)$, where
\begin{equation}
\label{eq4.25}
V_{n,0} = V^n = \left( \frac{\kappa}{\sinh\kappa}\right)^n , \qquad V_{n,\ell} = V_{n,0} W_{n,\ell},\ \ W_{n,0} = 1,
\end{equation}
and $W_{n,\ell}$ is found by iteration,
\begin{equation}
\label{eq4.26}
W_{n,\ell} = \frac{1}{\kappa^\ell} \int_0^\kappa \, V^{-n} \g^{\ell-1} \left( \cL_S - \frac{n^2}{2}\right) V_{n,\ell-1} d\g .
\end{equation}
\end{423Theorem}
\noindent This is the principal result of \S\ref{sec4} and
its proof occupies the rest of this section.

\begin{427Lemma}
\eqref{eq4.25} and \eqref{eq4.26} are consequences of the expansion \eqref{eq4.24}.
\end{427Lemma}

\begin{proof}
The action of $\cL_C$ on the integrands of \eqref{eq4.24} is given by
\addtocounter{equation}{1}
\begin{equation}
\label{eq4.28}
\cL_C = - \left( \frac12 \frac{\p^2}{\p\kappa^2} + n(\coth\kappa)\frac{\p}{\p\kappa}\right) - \frac12 \frac{\p^2}{\p\vp_1^2} 
= \cL_S - \frac12 \frac{\p^2}{\p\vp_1^2} ,
\end{equation}
with a slight abuse of notation when compared to \eqref{eq1.17} and \eqref{eq3.43}.
Let $V' = \p V / \p \kappa$.
Then
\begin{equation}
\label{eq4.29}
\cL_S e^{\frac{\kappa^2}{2t}} V = e^{\frac{\kappa^2}{2t}}
\left\{ - \frac{\kappa^2 V}{2t^2} - \frac{2\kappa V' + (2n\kappa\coth\kappa +1)V}{2t} + \cL_S V\right\} .
\end{equation}
Next $\p h (u + i \vp_1 )/\p\vp_1 = i\p h (u + i\vp_1)/\p u$ implies
\begin{align}
\label{eq4.30}
-\frac{\p^2}{\p\vp_1^2} \, e^{-\frac{u^2}{2t}} h
&=e^{-\frac{u^2}{2t}} \frac{\p^2 h}{\p u^2} \\
&=\frac{\p^2}{\p u^2} \left( e^{-\frac{u^2}{2t}}\right) h 
-\frac{\p}{\p u} \left( \frac{\p}{\p u} \left( e^{-\frac{u^2}{2t}}\right) h - e^{-\frac{u^2}{2t}}
\frac{\p h}{\p u}\right) ,\nonumber
\end{align}
hence
\begin{align}
\label{eq4.31}
&\cL_C \left( e^{-\frac{f}{2t}} V(\kappa)\right) \\
&= e^{-\frac{f}{2t}} \left\{ \frac{fV}{2t^2} - \frac{\kappa V' + (n\kappa \coth\kappa + 1)V}{t} + \cL_S V \right\} +
O\left( \frac{\p}{\p u}\right) .\nonumber
\end{align}
Consequently,
\begin{align}
\label{eq4.32}
&\left( \cL_C - \frac{\p}{\p t}\right) \frac{e^{At}}{t^\al} \int_{-\infty}^\infty e^{-\frac{f}{2 t}} V \, du \\
&= - \frac{e^{At}}{t^\al} \int_{-\infty}^\infty e^{-\frac{f}{2 t}} \Bigg\{ 
\frac{\kappa V' + (n\kappa \coth\kappa + 1-\al)V}{t} \nonumber\\
&\qquad\qquad\qquad\qquad\qquad - (\cL_S - A) V + \frac{\p V}{\p t}\Bigg\}\, du ,\nonumber
\end{align}
and this vanishes if the curly bracket in the integrand vanishes.
Set $V=V_n=V_{n,0}(\kappa) + V_{n,1} (\kappa) t + \cdots + V_{n,\ell} (\kappa) t^\ell + \cdots$, expand the
curly bracket in powers of $t$ and set the coefficient of each power of $t$ equal to zero.
The vanishing of the coefficient of the lowest power of $t$ requires
\begin{equation}
\label{eq4.33}
V'_{n,0} = \left( \frac{\al-1}{\kappa} - n\coth\kappa\right) V_{n,0},\qquad \mathrm{i.e.} \quad V_{n,0} = 
\frac{c\kappa^{\al-1}}{\sinh^n\kappa} .
\end{equation}
At the base point, $(\vt_1,\vp_1)=(0,0)$, $V_{n,0}=(-iu)^{\al-1} (-i\sinh u)^{-n}$ which is finite
and nonzero at $u=0$ if and only if $\al=n+1$.
Thus we choose
\begin{equation}
\label{eq4.34}
V_{n,0} = \left( \frac{\ka}{\sinh\ka}\right)^n = V^n ,
\end{equation}
in agreement with \eqref{eq4.25}.
The coefficient of $t^{\ell -1}$ vanishes if
\begin{equation}
\label{eq4.35}
\ka V'_{n,\ell} + n(\ka\coth\ka-1)V_{n,\ell} + \ell V_{n,\ell} = (\cL_S - A) V_{n,\ell-1} .
\end{equation}
With $V_{n,\ell}=V^n W_{n,\ell}$, \eqref{eq4.35} is reduced to
\begin{equation*}
\ka W'_{n,\ell} + \ell W_{n,\ell} = V^{-n} (\cL_S -A) V_{n,\ell-1} ,
\end{equation*}
in view of \eqref{eq4.33}, and this implies \eqref{eq4.26}.
\end{proof}

The derivation of \eqref{eq3.48} also yields
\begin{equation}
\label{eq4.36}
W_{n,1} = \binom{n}{2} \frac{d\log V}{d(\ka^2 /2)} ,\qquad \text{if }A=\frac{n^2}{2} .
\end{equation}

Consequently, replacing $\gamma$ by $\tilde\gamma$ in \eqref{eq3.39} yields the following part
of the integrand of \eqref{eq4.24}:
\begin{equation}
\label{eq4.37}
\frac{e^{\frac{n^2}{2}t}}{(2\pi t)^{n+1}} \sum_{k=-\infty}^{\infty} e^{-\frac{(\tgamma^2+2k\pi)^2}{2t}}
\sum_{\ell=0}^{n-1} v_{n,\ell} (\tgamma + 2k\pi)t^\ell .
\end{equation}
In particular, with $\kappa = i\tgamma$, $V_n (\kappa) = v_n (\tgamma)$ and \eqref{eq3.76} translates as
\begin{equation}
\label{eq4.38}
V_n(\kappa , t) = e^{-\frac{\ka^2}{2t}} \left( \frac{t}{2\pi} \frac{\p}{\p (\cosh \ka)}\right)^{n-1}
e^{\frac{\ka^2}{2t}} \frac{\ka}{\sinh\ka} ,
\end{equation}
which implies $V_{n,\ell}=0$, $\ell \ge n$, and the derivation of
\eqref{eq3.55} implies

\begin{439Proposition}
Let $\al>0$.
Then
\addtocounter{equation}{1}
\begin{equation}
\label{eq4.40}
P_C = \frac{\G (n) e^{\frac{n^2}{2} t}}{2(2\pi)^{n-1}(2\pi t)^2} \int_{-\infty}^\infty du \int_{\al- i\infty}^{\al+i\infty}
\frac{e^{\la/2t} d\la / (2\pi i)}{(\cosh \sqrt{\la + u^2} - \cosh \ka )^n} .
\end{equation}
\end{439Proposition}

The proof is somewhat lengthy and we shall break it up
into several steps.
Note that \eqref{eq4.24} is the residue expansion of \eqref{eq4.40} and therefore Theorem 4.23
is a consequence of Proposition 4.39 and
one also has \eqref{eq4.38}.

\begin{441Lemma}
$|\cosh \sqrt{\xi+i\nu}|$ is an increasing function of $|\nu|$ with
\addtocounter{equation}{1}
\begin{equation}
\label{eq4.42}
\min_\nu |\cosh \sqrt{\xi+i\nu}| = \cosh \sqrt \xi ,\qquad \xi,\nu \in \BR .
\end{equation}
\end{441Lemma}

\begin{proof}
Let $\sqrt{\xi + i\nu}=\s + i\chi$.
In view of \eqref{eq3.67} the square root maps the line $\xi+i\nu$, $-\infty < \nu < \infty$, onto the hyperbola $\s^2 -\xi^2 =\chi$.

1) $\xi>0$: $\s = \sqrt{\xi+\chi^2} > 0$ is an increasing function of $|\chi|$ which is an increasing function of
$|\nu|$; $\chi=0$ if and only if $\nu=0$.  Also,
\begin{equation}
\label{eq4.43}
|\cosh (\s + i\chi)|^2 = \cosh^2\s -\sin^2\chi =\cosh^2 \sqrt{\xi+\chi^2} - \sin^2\chi .
\end{equation}
Using the variables $(\xi,\chi)$ we find
\begin{equation*}
\frac{\p}{\p|\chi|} |\cosh (\s + i\chi)|^2 = 2|\chi|
\left(
\frac{\sinh 2\sqrt{\xi+\chi^2}}{2\sqrt{\xi+\chi^2}} - \frac{\sin 2\chi}{2\chi}\right) \geq 0 ,
\end{equation*}
and vanishes only when $\chi=0$, or, equivalently, when $\nu=0$.
Hence $|\cosh \sqrt{\xi+i\nu}|$ is an increasing function of $|\chi|$, and therefore an increasing function of $|\nu|$
with $\min\limits_\nu |\cosh \sqrt{\xi+i\nu}|=\cosh \sqrt{\xi}$.

2) $\xi<0$.  We switch to $(\xi,\s)$ coordinates with $\chi=\sqrt{\s^2-\xi}$ an increasing function of $|\s|$, which
is an increasing function of $|\nu|$; $\s =0$ if and only if $\nu=0$.
Now
\begin{align*}
|\cosh \sqrt{\xi+i\nu}|^2 &= \cosh^2 \s - \cos^2 \sqrt{\s^2-\xi} ,\\
\frac{\p}{\p |\s|} |\cosh \sqrt{\xi+i\nu}|^2 &= 2|\s| \left(
\frac{\sinh 2\s}{2\s} - \frac{\sin 2\sqrt{\s^2-\xi}}{2\sqrt{\s^2 -\xi}}\right) \geq 0 ,
\end{align*}
and vanishes only when $\s=0$, that is, when $\nu=0$.
Again, $|\cosh \sqrt{\xi+i\nu}|$ is an increasing function of $|\nu|$ with the minimum at $\nu=0$,
\begin{equation*}
\min_\nu |\cosh \sqrt{\xi+i\nu}| = \cos\sqrt{|\xi|} .
\qedhere
\end{equation*}
\end{proof}

\begin{444Lemma}
Set $x=\cos\vt_1 \cosh (u+i\vp_1)$.
Then
\addtocounter{equation}{1}
\begin{equation}
\label{eq4.45}
d = \cosh \sqrt{\la + u^2} - x = 0 \ \Rightarrow \ \Ree\la \leq 0 .
\end{equation}
\end{444Lemma}

\begin{proof}
$\la = \al + i\nu$. Then \eqref{eq4.42} implies
\begin{equation*}
\min_\nu  \big|\cosh \sqrt{\la + u^2}\big| = \cosh \sqrt{\al + u^2} .
\end{equation*}
In view of \eqref{eq4.43},
\begin{equation}
\label{eq4.46}
\max_{\vp_1} \big |\cosh (u+i\vp_1)\big| = \cosh u ,
\end{equation}
and therefore
\begin{equation*}
|d| \geq \cosh \sqrt{\al +  u^2} - \cosh u \geq \int_{|u|}^{\sqrt{u^2+\al}} \s \, d\s = \frac{\al}{2} .
\qedhere
\end{equation*}
\end{proof}

\begin{447Lemma}
The $P_C$-integrand \eqref{eq4.40} is summable; it and all its $u$, $\la$, $\vt_1$ and $\vp_1$ derivatives vanish
exponentially for large $|u|$ and $|\nu|$.
\end{447Lemma}

\begin{proof}
Define $\rho$ by $\nu = \rho (\al + u^2)$.
We shall make use of \eqref{eq4.43}, \eqref{eq3.67} and \eqref{eq4.46} to estimate $|d|$ from below:
\addtocounter{equation}{1}
\begin{align}
\label{eq4.48}
|d| &\geq \sinh \left( \frac{\sqrt{(\al+u^2)^2 +\nu^2} + \al + u^2}{2}\right)^{1/2} - \cosh u \\
&> \frac12 \left( e^{\sqrt{\al +u^2} (\frac{\sqrt{1+\rho^2}+1}{2})^{1/2}} - 1 - (e^{|u|} +1)\right) \nonumber \\
&= \frac12 e^{\sqrt{\al +u^2}(\frac{\sqrt{1+\rho^2}+1}{2})^{1/2}} \Bigg\{ 1 - e^
{-(\sqrt{\al+u^2}(\frac{\sqrt{1+\rho^2}+1}{2})^{1/2} - |u|)} \nonumber\\
&\qquad\qquad\qquad\qquad\qquad\quad - 2 e^{-\sqrt{\al+u^2}(\frac{\sqrt{1+\rho^2}+1}{2})^{1/2}}
\Bigg\} . \nonumber
\end{align}
The last two terms in the curly bracket vanish exponentially as $\rho\to \infty$, uniformly in $|u|$, hence
the curly bracket is bounded from below by $1/2$
when $\rho > M > 0$, $M$ sufficiently large.
When $\rho < M$, we note that
\begin{align}
\label{eq4.49}
&\sqrt{\al +u^2} \left( \frac{\sqrt{1+\rho^2}+1}{2}\right)^{1/2} - |u| \\
&\geq \sqrt{\al+u^2} - \sqrt{u^2} = \int_0^{\al} \frac{ds}{2\sqrt{u^2+s}} \geq \frac{1}{4|u|} \int_0^\al
ds = \frac{\al}{4|u|} ,\nonumber
\end{align}
and
\begin{align}
\label{eq4.50}
&1-e^{-(\sqrt{\al+u^2}(\frac{\sqrt{1+\rho^2}+1}{2})^{1/2} - |u|)}\\
&\leq 1 - e^{-\frac{\al}{4|u|}} = \int_0^{\frac{\al}{4|u|}} e^{-s} ds \le \frac{\al}{4|u|} .\nonumber
\end{align}
The third term in the curly bracket is negligible as $|u|\to\infty$.
Consequently, when either $|u|>M$ or $\rho>M$, $M>0$ sufficiently large, one has
\begin{equation}
\label{eq4.51}
|d| \geq A e^{\sqrt{\al+u^2}(\frac{\sqrt{1+\rho^2}+1}{2})^{1/2}} (1+|u|)^{-1} ,
\end{equation}
$A>0$.
$d$ does not vanish on the domain of integration, so on the bounded domain $0<\rho<M$, $0<|u|<M$
we can always estimate $|d|$ from below by \eqref{eq4.51} with a judicious choice of $A$.
Hence \eqref{eq4.51} holds for all $\rho$ and $u$, and returning to $(u,\nu)$ we have
\begin{align}
\label{eq4.52}
|d|^{-n} &\leq A^{-n} e^{-n (\frac{\sqrt{(\al+u^2)^2+\nu^2}+\al+u^2}{2})^{1/2}} (1+|u|)^n \\
&\leq A^{-n} e^{-\frac{n}{2} (\sqrt{|\nu |} + |u|)} (1+|u|)^n ,\nonumber
\end{align}
$(u,\nu)\in\BR^2$, hence $|d|^{-n}$ is summable.
Differentiating the integrand in $u$, $\la$, $\vt_1$ or $\vp_1$ may add a $\sinh$ or $\cosh$ term to the
numerator, but this added increase is cancelled by the extra power of $d$ in the denominator.
\end{proof}

\begin{453Lemma}
$P_C$ of \eqref{eq4.40} is a solution of the heat equation $\De_C P_C = \p P_C / \p t$.
\end{453Lemma}

\begin{proof}
It suffices to show that $\cL_C P_C = \p P_C / \p t$.
Write
\begin{equation*}
P_C = \frac{Be^{\frac{n^2}{2}t}}{t^2} \int_{\al-i\infty}^{\al+i\infty}
e^{\la/2t} k(\la , \vt_1 , \vp_1)d\la ,
\end{equation*}
\begin{equation*}
k = \int_{-\infty}^{\infty} \frac{du}{h(\la, u,x)^n} ,\qquad h =\cosh \sqrt{\la + u^2} - x .
\end{equation*}
Then
\begin{align*}
\left( \cL_C - \frac{\p}{\p t}\right)\, \frac{e^{\frac{n^2}{2}t} e^{\la/2t} k}{t^2} 
&= \frac{e^{\frac{n^2}{2}t} e^{\la/2t}}{t^2} \left\{ 2\la \frac{\p^2 k}{\p\la ^2} + \left(\cL_C
-\frac{n^2}{2}\right)k\right\} \\
&\quad + \frac{e^{\frac{n^2}{2}t}}{t^2} \frac{\p}{\p\la} e^{\la/2t} \left\{ \frac{\la k}{t} + 2
\left( k - \la \frac{\p k}{\p \la}\right)\right\} ,
\end{align*}
where we exchanged the $t$-singularities with $\la$-derivatives.
In view of Lemma 4.47, the $\la$-integral of the last line vanishes, so
\begin{align*}
&\cL_C P_C - \frac{\p P_C}{\p t} \\
&=\frac{Be^{\frac{n^2}{2}t}}{t^2} \int_{\al-i\infty}^{\al+i\infty}
e^{\la/2t} \left\{ 2\la \frac{\p^2 k}{\p\la^2} + \left( \cL_S - \frac{n^2}{2}\right)k - \frac12
\frac{\p^2 k}{\p\vp_1^2}\right\} d\la .
\end{align*}
With $x=\cos\vt_1 \cosh (u + i\vp_1 )$ we integrate by parts and obtain
\begin{align*}
\frac{\p k}{\p\vp_1}
&= \int_{-\infty}^\infty \frac{\p h^{-n}}{\p x} \frac{\p x}{\p\vp_1} du = i \int_{-\infty}^\infty
\frac{\p h^{-n}}{\p x} \frac{\p x}{\p u} du \\
&= i \int_{-\infty}^\infty \left\{ \left[ \frac{\p h^{-n}}{\p x} \frac{\p x}{\p u} + \frac{\p h^{-n}}{\p u}\right] -
\frac{\p h^{-n}}{\p u}\right\} du 
= -i \int_{-\infty}^{\infty} \frac{\p h^{-n}}{\p u} d u ,
\end{align*}
since the square bracket is the total derivative of $h^{-n}$ with respect to $u$ and its integral vanishes.
A second derivative in $\vp_1$ yields the negative of the second derivative in $u$ under the integral,
so
\begin{align*}
&\left(  \frac{Be^{\frac{n^2}{2}t}}{t^2} \right)^{-1}
\left( \cL_C P_C - \frac{\p P_C}{\p t}\right) \\
&= \int_{\al-i\infty}^{\al+i\infty}
e^{\la/2t} d\la \int_{-\infty}^\infty \left\{ 2\la \frac{\p^2 h^{-n}}{\p \la ^2} + \left( \cL_S - \frac{n^2}{2}\right) h^{-n} +
\frac12 \frac{\p^2 h^{-n}}{\p u^2}\right\} du ,
\end{align*}
and it suffices to show that the curly bracket vanishes.
With $\la = \mu^2$,
\begin{align*}
\Big\{\cdots\Big\} &=\frac12 \left( \frac{\p^2}{\p\mu^2} - \frac{1}{\mu} \frac{\p}{\p\mu} + 2\cL_S - n^2 + \frac{\p^2}{\p u^2}\right) h^{-n}\\
&= \frac{1}{2h^n} \Bigg\{ 
\frac{n(n+1)}{h^2} 
\left[ \left( \frac{\p h}{\p\mu}\right)^2 + (1-x^2) 
\left(\frac{\p h}{\p x}\right)^2 
+ \left( \frac{\p h}{\p u}\right)^2 \right] \\
&\qquad\qquad - \frac{n}{h} \left[ \frac{\p^2 h}{\p\mu^2} - \frac{1}{\mu} \frac{\p h}{\p \mu} + 2\cL_S h +
\frac{\p^2 h}{\p u^2}\right] - n^2\Bigg\} \\
&=\frac{1}{2h^n} \Bigg\{ \frac{n(n+1)}{h^2} \left( \sinh^2 \sqrt{\mu^2 +u^2} +1-x^2\right) \\
&\qquad \qquad - \frac{n}{h} \big( \cosh \sqrt{\mu^2+ u^2} + x + 2nx\big) - n^2\Bigg\} \\
&= \frac{n^2}{2h^n} \left\{ \frac{1}{h} \big( \cosh \sqrt{\mu^2 + u^2} + x\big) - \frac{2 x}{h} - 1\right\} \\
&= 0 .
\qedhere
\end{align*}
\end{proof}

\begin{454Lemma}
One has
\addtocounter{equation}{2}
\begin{equation}
\label{eq4.55}
\int_{S^{2n+1}} P_C dS = 1 .
\end{equation}
\end{454Lemma}

\begin{proof}
Integrating functions of $(\vt_1,\vp_1)$ only on $S^{2n+1}$ the volume element is reduced to
\begin{equation}
\label{eq4.56}
dS = \frac{\pi^n}{\G (n)} \sin^{2(n-1)} \vt_1 \sin 2\vt_1 \, d\vt_1 d\vp_1 ,
\end{equation}
see \eqref{eq2.42}.
Therefore,
\begin{align*}
&\int_{S^{2n+1}} P_C dS \\
&= 
\frac{(-1)^{n-1} e^{\frac{n^2} {2}t}}{2^{n+2} \G(n) \pi t^2} \int_{-\infty}^\infty e^{-\frac{u^2}{2t}} du 
\frac{1}{2\pi i} \int_{\alpha-i\infty}^{\alpha+i\infty} d\la e^{\la / 2t}
\frac{d^{n-1}}{d(\cosh \sqrt\la)^{n-1}} \ \cdot \\
&\qquad\cdot \int_0^{\pi/2} \sin^{2(n-1)} \vt_1 \sin 2\vt_1 d\vt_1 \int_{-\pi}^\pi
\frac{d\vp_1}{\cosh\sqrt\la - \cos\vt_1 \cosh (u+ i\vp_1)} .
\end{align*}
With $A=(\cosh \sqrt\la )/\cos\vt_1$ the last integral becomes
\begin{equation}
\label{eq4.57}
\frac{1}{\cos\vt_1} \int_{-\pi}^{\pi} \frac
{d\vp_1} {A-\cosh u \cos\vp_1 - i\sinh u \sin\vp_1} = \frac{I}{\cos\vt_1} ,
\end{equation}
and the customary substitution $v=\tan (\vp_1 / 2)$ yields
\begin{equation}
\label{eq4.58}
I = \frac{2}{A+\cosh u} \int_{-\infty}^\infty
\frac{dv}{\Big( v - i\frac{\sinh u}{A+\cosh u}\Big)^2+\frac{A^2-1}{(A+\cosh u)^2}} = \frac
{2\pi}{\sqrt{A^2-1}} ,
\end{equation}
since $A^2\ne 1$.
Thus
\begin{equation*}
2\pi \int_0^{\pi/2}
\frac{\sin^{2(n-1)} \vt_1 \sin 2\vt_1 d\vt_1}{\sqrt{\cosh^2 \sqrt\la - \cos^2\vt_1}} = 2\pi
\int_0^1 \frac{(1-u)^{n-1} du}{\sqrt{\cosh^2 \sqrt\la - u}} ,
\end{equation*}
and in view of Lemmas 4.59 and 3.53 one has
\[
\int_{S^{2n+1}} P_C dS = \frac{e^{\frac{n^2}{2} t}}{4\pi int^2} \int_{-\infty}^\infty
e^{-\frac{u^2}{2t}} du \int_{\alpha-i\infty}^{\alpha+i\infty} e^{\la/2t-n\sqrt\la} d\la 
= 1 .
\qedhere
\]
\end{proof}

\begin{459Lemma}
$A\not\in [0,1]$.
Then
\addtocounter{equation}{1}
\begin{equation}
\label{eq4.60}
\frac{d^k}{dA^k} \int_0^1 \frac{(1-u)^k du}{\sqrt{A^2-u}} = \frac{(-1)^k 2^{k+1} k!}{k+1}
\big( A -\sqrt{A^2-1}\big)^{k+1} .
\end{equation}
\end{459Lemma}

\begin{proof}
We argue by induction on $k$.  The result is elementary when $k=0$.
For $k>0$ an integration by parts yields
\begin{align}
\label{eq4.61}
q_k(A)  &\equiv \int_0^1 \frac{(1-u)^k du}{\sqrt{A^2-u}} = -2 \int_0^1 (1-u)^k d\sqrt{A^2-u} \\
&= 2A - 2k (A^2-1) q_{k-1} (A) - 2k q_k (A) ,\nonumber
\end{align}
so,
\begin{equation}
\label{eq4.62}
q_k (A) = \frac{2A}{2k+1} - \frac{2k}{2k+1} (A^2-1) q_{k-1} (A) .
\end{equation}
Differentiating \eqref{eq4.61} and then integrating by parts one has
\begin{equation}
\label{eq4.63}
q_k^{(1)} (A) = -2A \int_0^1 (1-u)^k d(A^2-u)^{-1/2} = 2-2k A q_{k-1} (A) .
\end{equation}
In particular, \eqref{eq4.63} gives the result for $k=1$ and we may assume that $k\geq 2$.
Our induction hypothesis is that $q_j^{(j)} = C_j \big( A-\sqrt{A^2-1})^{j+1}$, $j<k$.
Differentiating \eqref{eq4.63} and \eqref{eq4.62} one obtains
\begin{align}
\label{eq4.64}
q_k^{(k)} (A)
&= -2k \Big\{ Aq_{k-1}^{(k-1)} (A) + (k-1) q_{k-1}^{(k-2)} (A)\Big\} ,\\
q_k^{(k)}(A) &= -\frac{2k}{2k+1} \Big\{ (A^2-1)q_{k-1}^{(k)} (A) + 2k A q_{k-1}^{(k-1)} (A)\nonumber \\
&\qquad\qquad\qquad + k (k-1) q_{k-1}^{(k-2)} (A)\Big\} .\nonumber
\end{align}
Equating the two right hand sides we introduce the induction hypothesis for $q_{k-1}^{(k-1)}$ and find its
derivative $q_{k-1}^{(k)}$.
This yields
\begin{equation*}
q_{k-1}^{(k-2)} (A) = \frac{-C_{k-1}}{(k+1)(k-1)} \big( A - \sqrt{A^2-1}\big)^k 
\big( A+k \sqrt{A^2-1}\big) ,
\end{equation*}
which we substitute into \eqref{eq4.64} together with the induction hypothesis for $q_{k-1}^{(k-1)}$ and
find
\begin{equation*}
q_k^{(k)} (A) = -\frac{2k^2}{k+1} C_{k-1} \big( A - \sqrt{A^2-1}\big)^{k+1} .
\end{equation*}
A simple calculation leads to the explicit $C_k$ of \eqref{eq4.60}.
\end{proof}

To complete the proof of Propositon 4.39 we still need to prove that $\lim\limits_{t\to 0} P_C (t)=
\delta (1-z)$, $1=(1,\ldots,0)$, $z=(z_1,\ldots, z_{n+1})$.
In view of Lemma 4.54 it suffices to show that $P_C (t)$ converges to zero as $t\to 0$ uniformly
on compact subsets of $S^{2n+1}$ which do not contain $(1,0,\ldots,0)$; in other words, on compact
subsets of $S^{2n+1}$ on which $\cos\vt_1 \cos\vp_1 < 1$.
We shall prove this uniform convergence first on sets where $\vt_1$ is bounded away from zero, then on
sets where $\vp_1$ is bounded away from zero.

\begin{465Lemma}
$\lim\limits_{t\to 0} P_C (t) = 0$ uniformly on any compact subset $\Omega$ of $S^{2n+1}$ where
$\vt = \min\limits_\Omega \vt_1 >0$.
\end{465Lemma}

\begin{proof}
Staying with the notation $\la = \al +i\nu$, $\al > 0$,
\begin{equation*}
P_C = \frac{Be^{\frac{n^2} {2}t} e^{\frac{\al}{2t}}}{t^2}
\int_{-\infty}^{\infty} e^{i\frac{\nu}{2t}} d\nu \int_{-\infty}^\infty \frac{du}{d^n} .
\end{equation*}
We intend to show that with our hypothesis $\al$ may be chosen negative, while $|d|^{-n}$ remains
integrable on $\BR^2$, uniformly with respect to $\vt_1 \geq \vt > 0$ and $\vp_1$.
Hence the result.
It suffices to do this on the first quadrant $u>0$, $\nu>0$.
With the $\vt>0$ of the hypothesis,
\addtocounter{equation}{1}
\begin{equation}
\label{eq4.66}
|d| \geq \cosh \sqrt{\al + u^2} - \cos\vt \cosh u > 0 ,\qquad \al > 0 ,
\end{equation}
see \eqref{eq4.42}, \eqref{eq4.46}.
Set $\cos\vt = e^{-\mu}$, $\mu>0$.
Then $\cos\vt\cosh u \sim \cosh (u-\mu)$, $u\sim\infty$, and this suggests that \eqref{eq4.66}
holds with some $\al < 0$ when $u\sim \infty$.
To this end we define $u_\vt$ by $\cos\vt \cosh u_\vt = 1$, and let $\tilde u = \tilde u (u)$ be given by
\begin{equation}
\label{eq4.67}
\cos\vt \cosh u = \cosh \tilde u (u) ,\qquad u > u_\vt .
\end{equation}
Then $u=\tilde u (u) + \ve(u)$, $\ve(u) > 0$, $\ve(u_\vt)=u_\vt$, $\tilde u (u_\vt)=0$ and $\ve(\infty)=\mu$.
The $u$-derivative of \eqref{eq4.67} yields
\begin{equation}
\label{eq4.68}
\cos\vt \sinh u = \big( 1- \ve'(u)\big) \sinh \tilde u ,
\end{equation}
so $\ve' (u) < 1$ and $\tilde u\vphantom{u}'(u) = 1-\ve' (u) > 0$.
Subtracting the square of \eqref{eq4.68} from the square of \eqref{eq4.67} gives 
$-\sin^2\vt = \ve' (2-\ve')\sinh^2 \tilde u$, hence $\ve'(u) < 0$ and we have derived

\begin{469Lemma}
$u-\tilde u (u)$ decreases from $u_\vt$ at $u=u_\vt$ to $\mu$ at $u=\infty$.
\end{469Lemma}
\renewcommand{\qedsymbol}{}
\end{proof}

\begin{proof}[Proof of Lemma 4.65 continued]
1) $u_\vt < M < u$.
Choose a negative $\al$, $-\mu^2 / 2 < \al < 0$.
Then $\al+u^2=\al+(u-\mu+\mu)^2 \geq \al + \mu^2 + (u-\mu)^2 > \mu^2/2 + (u-\mu)^2$, and
\addtocounter{equation}{1}
\begin{align}
\label{eq4.70}
|d| &> \frac12 \big( e^{((\sqrt{(\al+u^2)^2+\nu^2}+\al+u^2)/2)^{1/2}} - 1 - (e^{u-\mu} +1)\big) \\
&> \frac12   e^{((\sqrt{(\al+u^2)^2+\nu^2}+\al+u^2)/2)^{1/2}}\ \cdot \nonumber \\
&\cdot \Big\{ 1 - e^{-(\sqrt{\mu^2/2+(u-\mu)^2}-(u-\mu))} 
- 2 e^{-\sqrt{\al+ u^2}}\Big\} .\nonumber
\end{align} 
The last term in the curly bracket vanishes exponentially
as $u\to\infty$, and the difference of the first two terms is bounded from below by $\mu^2/8(u-\mu)$,
see \eqref{eq4.50}.
Therefore,
\[
|d| > Ae^{\sqrt{\frac{|\nu| + (u-\mu)^2}{2}}} (1+u)^{-1} ,
\]
and
\begin{equation*}
|d|^{-n} < A^{-n}  e^{-\frac{n}{2}(\sqrt\nu + |u-\mu|)} (1+u)^n ,
\end{equation*}
which is integrable on $u>M$, $0 < |\nu| < \infty$.

2) $0 < u < M$. With $\al < 0$, $\al+u^2$ becomes negative for small $u^2$.
To compensate we let $\nu$ be large.
Then \eqref{eq4.48} yields
\begin{align*}
|d| &\geq \sinh \left( \frac{\sqrt{(\al+u^2)^2+\nu^2}+\al+u^2}{2}\right)^{1/2} -\cosh M \\
&\geq \sinh \sqrt{\frac{\nu-|\al|}{2}} - \cosh M \\
&= \frac12 e^{\frac{1}{\sqrt 2} \sqrt{\nu-|\al|}} \Big( 1 - (1+2\cosh M)e^{-\frac{1}{\sqrt 2}\sqrt{\nu-|\al|}}\Big) \\
&\geq \frac14 e^{\frac{1}{\sqrt 2}\sqrt{\nu - |\al|}} , \qquad \nu >M' > M ,
\end{align*}
when $M'$ is sufficiently large.
Consequently, one has
\begin{equation*}
|d|^{-n} \leq 4^n e^{-\frac{n}{\sqrt 2} \sqrt{\nu-|\al|}}
\end{equation*}
which is integrable on $\{ 0 < u < M , \ M'< \nu < \infty \}$.

3) $0 < u < M ,\ 0 < \nu < M'$.
In view of Lemma 4.41 and formulas \eqref{eq4.45} and \eqref{eq4.46}, one has 
\begin{equation*}
|d(\al = 0)| \geq \cosh u (1-\cos\vt)  > 0 ,
\end{equation*}
and by continuity $d\ne 0$ for some $\al < 0$ and all $\vp_1$ and $\vt_1 > \vt$.
\end{proof}

When $\vt_1 \sim 0$ we need $\vp_1$ bounded away from zero.
The work with nonzero $\vp_1$ requires the extension of \eqref{eq4.67} from $u \in \BR$ to
$u + i\vp_1 \in \BC$:
\begin{equation}
\label{eq4.71}
\cos\vt_1 \cosh (u+i\vp_1) = \cosh (\tilde u + i\tilde\vp_1) ,
\end{equation}
where $\tilde u = \tilde u (\vt_1 , \vp_1 , u)$ and $\tilde\vp_1 = \tilde\vp_1 (\vt_1 , \vp_1 , u)$ represent
the solutions of \eqref{eq4.71} in $0 < \tilde u < \infty$, $-\pi < \tilde\vp_1 \leq \pi$ when $0<u<\infty$,
$-\pi < \vp_1 \leq\pi$; to simplify matters we recall that $\cosh (u+i\vp_1)$ maps
$\big\{ 0 < u < \infty ,\ 0\leq \vp_1 \leq \pi/2\big\}$ one-to-one and onto the first quadrant of the
complex plane.
One separates \eqref{eq4.71} into its real and imaginary parts,
\begin{equation}
\label{eq4.72}
\cos\vt_1 \cosh u \cos\vp_1 = \cosh \tilde u \cos \tilde \vp_1 ,
\end{equation}
\begin{equation}
\label{eq4.73}
\cos\vt_1 \sinh u \sin\vp_1 = \sinh \tilde u \sin \tilde\vp_1 .
\end{equation}
Note that $\tilde\vp_1 (\vp_1=0)=0$, $\tilde\vp_1 (\vp_1 = \pi/2) = \pi/2$, and by continuity
$0 < \vp_1 < \pi/2 \Rightarrow 0 < \tilde\vp_1 < \pi/2$.
Also,
\begin{equation}
\label{eq4.74}
\tilde\vp_1 (\vp_1) = \begin{cases}
\pi- (\pi-\vp_1)^\sim ,&\pi/2 < \vp_1 \leq \pi ,\\
- (-\vp_1)^\sim, &-\pi < \vp_1 < 0, 
\end{cases}
\end{equation}
\begin{equation}
\label{eq4.75}
\tilde u (\vp_1) = \begin{cases}
\tilde u (\pi-\vp_1) ,&\pi/2 < \vp_1 \leq \pi ,\\
\tilde u (-\vp_1), &-\pi < \vp_1 < 0, 
\end{cases}
\end{equation}
and we may restrict our attention to $0 <  \vp_1 < \pi/2$.

\begin{476Lemma}
$0 < \vp_1 < \pi/2$, $\cos\vt_1 = e^{-\mu}$, $\mu >0$.
Then
\medskip	

\noindent {\rm (i)} $\tilde u < u$,
\medskip
	
\noindent {\rm (ii)} $\vp_1 < \tilde\vp_1$,
\medskip

\noindent {\rm (iii)} $\displaystyle{\frac{\p\tilde u}{\p\vt_1} < 0}$, $\displaystyle{\frac{\p\tilde\vp_1}{\p \vt_1} > 0}$,
\medskip

\noindent {\rm (iv)} $\displaystyle{\frac{\p \tilde u}{\p u} > 0}$, $\displaystyle{\frac{\p\tilde\vp_1}{\p u} < 0}$,
\medskip

\noindent {\rm (v)} $\tilde u (u) \sim u - \mu$, $\tilde\vp_1 (u) \sim \vp_1$, $u\to \infty$.
\end{476Lemma}

\begin{proof}
(i) Assume that $\tilde u \geq u$.
Then \eqref{eq4.72}, \eqref{eq4.73} yield
\begin{equation*}
\cos\vt_1 \cosh u < \cosh u \leq \cosh \tilde u  \ \Rightarrow\ \vp_1 < \tilde\vp_1 ,
\end{equation*}
\begin{equation*}
\cos\vt_1 \sinh u < \sinh u \leq \sinh \tilde u \ \Rightarrow \ \vp_1 > \tilde\vp_1 .
\end{equation*}
Contradiction, hence $\tilde u < u$.
\medskip

(ii) Again, from \eqref{eq4.72}, \eqref{eq4.73}
\begin{equation*}
\frac{\tanh u}{\tanh\tilde u} = \frac{\tan\tilde\vp_1}{\tan \vp_1} ,
\end{equation*}
and in view of (i) this implies (ii).
\medskip

(iii) We differentiate \eqref{eq4.72} and \eqref{eq4.73} with respect to $\vt_1$, then eliminating
$\p \tilde\vp_1 / \p \vt_1$ yields $\p\tilde u / \p \vt_1 < 0$, and eliminating $\p \tilde u / \p \vt_1$ gives
$\p \tilde\vp_1 / \p\vt_1 > 0$.
\medskip

(iv) Same argument as in (iii).
\medskip

(v) Note that
\begin{equation*}
\cosh (\tilde u + i\tilde\vp_1 ) \sim \frac12 e^{u-\mu} e^{i\vp_1} ,\qquad u\to\infty ,
\end{equation*}
hence $\tilde u (u) \sim u-\mu$ and $\tilde\vp_1 (u)\sim \vp_1$ when $u\sim\infty$.
\end{proof}
Of course, $\vt_1=0$ implies that $\tilde u + i\tilde\vp_1 = u+i\vp_1$.
Also, $\tilde\vp_1 (\vp_1 = \pm\pi)=\pm\pi$ and \eqref{eq4.73} yields
\addtocounter{equation}{1}
\begin{equation}
\label{eq4.77}
\vp_1 \ne 0 \ \Rightarrow \ \tilde u (u=0) = 0 ,\qquad
\cos\tilde\vp_1 (u=0)=\cos\vt_1 \cos\vp_1 ;
\end{equation}
note the difference from $\vp_1 = 0$, when $\tilde u (u_{\vt_1}) = 0$, and
$u_{\vt_1}=0$  if and only if $\vt_1=0$.

We shall show that $P_C\to 0$ as $t\to 0$, uniformly on sets where 
(i) $0 < \delta < |\vp_1| < \pi - \delta$,  and (ii) $\vp_1 $ is near $\pm\pi$.

\begin{478Lemma}
Let $0<\de \leq |\vp_1|<\pi-\de'$, $\de'  >0$, and $m=1,2,\ldots,$ sufficiently large.
Then
\addtocounter{equation}{1}
\begin{align}
\label{eq4.79}
&\frac{\G(n)t^{n-1}}{2} \int_{\al-i\infty}^{\al+i\infty}
\frac{e^{\la / 2t}}{(\cosh \sqrt\la - \cos\vt_1 \cosh (u+i\vp_1))^n} 
\frac{d\la}{2\pi i} \\
&= \sum_{k=-m}^m e^{\frac{(\kappa + i2k\pi)^2}{2t}} V_n (\kappa + i2k\pi, t)+R_m,\nonumber
\end{align}
where $\al > u^2$, $\kappa = \cosh^{-1} \big( \cos\vt_1 \cosh (u+i\vp_1)\big)=\tilde u + i\tilde\vp_1$, $V_n$ is given
by \eqref{eq4.38}, and
\begin{equation}
\label{eq4.80}
|R_m| \leq C(\de') (2m+1) \pi t^{n-1} e^{\frac{\tilde u^2 - (2m+1)^2\pi^2}{2t}} (1+ |\tilde u|)^{n+\frac32}
e^{-\frac{n|\tilde u|}{2}} .
\end{equation}
\end{478Lemma}

\begin{proof}
The idea is to move the line of integration to the left and estimate the new integral.
$d=\cosh \sqrt\la -\cos\vt_1 \cosh (u+i\vp_1)=\cosh \sqrt\la - \cosh (\tilde u + i\tilde\vp_1)=0$
if $\la =\la_k =\big(\tilde u +i(\tilde\vp_1 + 2k\pi)\big)^2 = \tilde u^2 - (\tilde \vp_1 + 2k\pi)^2 +
i 2\tilde u (\tilde\vp_1 + 2k\pi)$, $k\in\BZ$; recall that $0<\tilde u < u$ and $|\tilde\vp_1|<\pi - \de'$.
To avoid $\la_k$, $k\in\BZ$, we move the line integration to $\Ree\la =\tilde u^2 - (2m+1)^2\pi^2=\al_m$ and set
\begin{equation}
\label{eq4.81}
R_m = \frac{\G (n) t^{n-1}}{2} \int_{\Ree\la =\al_m} \frac{e^{\la/2t}}{d^n} \frac{d\la}{2\pi i} .
\end{equation}
Again $\sqrt\la = \s + i\chi$, which yields
\[
d = 2\sinh \frac{\s+\tilde u + i(\chi+\tilde\vp_1)}{2} \sinh\frac{\s-\tilde u +i(\chi-\tilde\vp_1)}{2} ,
\]
\[
|d|^2 = 4 \left( \sinh^2 \frac{\s +\tilde u}{2} +\sin^2 \frac{\chi+\tilde\vp_1}{2}\right)
\left( \sinh^2 \frac{\s-\tilde u}{2} + \sin^2 \frac{\chi-\tilde\vp_1}{2}\right) .
\]
We shall switch the variable of integration from $\nu$ to $\s$; it suffices to work with $\s > 0$,
$\tilde u > 0$ and $\nu>0$.

1) $|\s^2 - \tilde u^2| \geq 2$. Here
\begin{align*}
|d|^2 
\geq 4\sinh^2 \frac{\s +\tilde u}{2}\sinh^2 \frac{\s-\tilde u}{2} &\geq (\s-\tilde u)^2 \sinh^2\frac{\s+\tilde u}{2} \\
&\geq \left( \sinh^2 \frac{\s+\tilde u}{2}\right) \Big/ \left( \frac{\s+\tilde u}{2}\right)^2 ,
\end{align*}
and $\sinh\g /\g > e^\g / (1+2\g) > e^\g / 2(1+\g)$ implies
\[
\int_{|\s^2-\tilde u^2|>2} |d|^{-n} d\nu < 2^n \int_{|\s^2-\tilde u^2|>2} 
\left( 1 + \frac{\s+\tilde u}{2}\right)^n e^{-n\frac{\s+\tilde u}{2}} d\nu .
\]
In view of \eqref{eq3.67}, $\nu = 2\s \chi = 2\s \sqrt{\s^2 - \al_m}$, and
\[
\frac{d\nu}{d\s} = 2 \left( \sqrt{\s^2-\al_m} + \frac{\s^2}{\sqrt{\s^2-\al_m}} \right);
\]
$0<\nu < \infty$ goes into (i) $0<\s < \infty$ when $\al_m<0$, or (ii) $\sqrt{\al_m}<\s<\infty$ when $\al_m>0$.

(i) $\al_m < 0$: $\s^2 < \s^2 - \al_m < \s^2 + (2m+1)^2 \pi^2$, so
\[
\frac{d\nu}{d\s} < 2 \big( \s + (2m+1) \pi +\s\big) < 2(2m+1) \pi (1+\s) ,
\]
\[
\int_{\tilde u^2 + 2<\s^2} |d|^{-n} d\nu < C(n) (2m+1)\pi e^{-\frac{n\tilde u}{2}} \int_{\tilde u < \s} (1+\s)^{n+1}
e^{-n\s/2} d\s ,
\]
\[
\int_{\s^2 < [\tilde u^2 -2]_+} |d|^{-n} d\nu < C(n) (2m+1)\pi (1+\tilde u)^{n+1} e^{-\frac{n\tilde u}{2}}
\int_{\s < \tilde u} e^{-\frac{n\s}{2}} d\s ,
\]
which yield
\begin{equation}
\label{eq4.82}
\int_{|\s^2-\tilde u^2|>2} |d|^{-n} d\nu < C(n) (2m+1) \pi (1+\tilde u)^{n+\frac32} e^{-\frac{n\tilde u}{2}} .
\end{equation}

(ii) $\al_m\geq 0$: $\sqrt{\s} \sqrt{\s - \sqrt{\al_m}} \leq \sqrt{\s^2 - \al_m} \leq \s $, since $\s^2>\al_m$
according to \eqref{eq3.67}, so
\[
\frac{d\nu}{d\s} < 2 \left( \s + \frac{\s^{3/2}}{\sqrt{\s-\sqrt{\al_m}}}\right) < \frac{4(1+\s)^{3/2}}{\sqrt{\s-\sqrt{\al_m}}} ,
\]
\[
\int_{\tilde u^2 +2<\s^2} |d|^{-n} d\nu < C(n) e^{-\frac{n\tilde u}{2}} \int_{\tilde u < \s}
(1+\s)^{n+\frac32}  \frac{e^{-n\s /2}d\s}{\sqrt{\s-\sqrt{\al_m}}} ,
\]
\[
\int_{\s^2 < [\tilde u^2 -2]_+} |d|^{-n} d\nu < C(n) (1+\tilde u)^{n+\frac32} e^{-\frac{n\tilde u}{2}}
\int_{\sqrt{\al_m}}^{\tilde u} \frac{e^{-n\sigma / 2} d\sigma}{\sqrt{\s - \sqrt{\al_m}}} ,
\]
and \eqref{eq4.82} holds for all $\al_m$.

2) $|\s^2-\tilde u^2|<2$.
We shall estimate
\[
|d|^2 \geq 4\left( \sinh^2 \frac{\s+\tilde u}{2} + \sin^2 \frac{\chi+\tilde\vp_1}{2}\right) \sin^2
\frac{\chi-\tilde\vp_1}{2} .
\]
Here
\begin{align*}
\chi &= \sqrt{\s^2-\al_m} = \sqrt{\s^2 - \tilde u^2 + (2m+1)^2 \pi^2} \\
&= (2m+1)\pi + (2m+1)\pi \left( \sqrt{1 + \frac{\s^2 - \tilde u^2}{(2m+1)^2\pi^2}} - 1 \right) ,
\end{align*}
\[
\sin \frac{\chi-\tilde\vp_1}{2} = \pm \sin \left( \frac{\pi - \tilde\vp_1}{2} + \frac{(2m+1)\pi}{2} \left\{
\sqrt{ 1 + \frac{\s^2-\tilde u^2}{(2m+1)^2\pi^2}} - 1 \right\}\right) ,
\]
and
\[
\left| \sqrt{ 1 + \frac{\s^2-\tilde u^2}{(2m+1)^2\pi^2}} - 1 \right| \leq
\sqrt{ 1 + \frac{2}{(2m+1)^2\pi^2}} - 1 \leq \frac{1}{(2m+1)^2\pi^2} ,
\]
implies that
\[
\left| \sin\frac{\chi-\tilde\vp_1}{2}\right| > \sin \left( \frac{\de'}{2} - \frac{1/2}{(2m+1)\pi}\right) >\sin\frac{\de'}{4}
\]
for large $m$; by hypothesis $\delta' < \pi - |\tilde\vp_1 |$.
Also,
\[
\sin\frac{\chi+\tilde\vp_1}{2} = \pm \sin \left( \frac{-\pi+\tilde\vp_1}{2} + \frac{(2m+1)\pi}{2}
\left\{ \sqrt{1 + \frac{\s^2-\tilde u^2}{(2m+1)^2\pi^2}} - 1 \right\} \right) ,
\]
so
\[
\left| \sin \frac{\chi+\tilde\vp_1}{2}\right| > \sin \frac{\de'}{4} ,
\]
and then
\[
\sinh^2 \frac{\s + \tilde u}{2} + \sin^2 \frac{\chi+\tilde\vp_1}{2} \geq C(\de') \cosh^2 \frac{\s+\tilde u}{2} ,\qquad
C(\de')>0,
\]
yields
\[
|d|^2 > C(\de') \cosh^2 \frac{\s+\tilde u}{2} > C(\de') e^{\s+\tilde u} .
\]
Here, $-2<\s^2 -\tilde u^2 < 2$, hence $(2m+1) \pi/2 < \sqrt{\s^2 -\al_m} < 2(2m+1)\pi$, and
$d\nu / d\s < 8 (2m+1)\pi$.
The length of the interval of the $\s$-integration for $\s < \tilde u$ and also for $\s > \tilde u$ is less
than 2 units, therefore one has
\begin{equation}
\label{eq4.83}
\int_{|\s^2-\tilde u^2|<2} |d|^{-n} d\nu < C(\de') (2m+1) \pi e^{-\frac{n\tilde u}{2}} .
\end{equation}
Thus \eqref{eq4.81}, \eqref{eq4.82} and \eqref{eq4.83} imply \eqref{eq4.80}; the term 
$e^{\la / 2t}$ in the integrand adds the factor $\exp \big(\tilde u^2 / 2t - (2m+1)^2\pi^2 / 2t\big)$ when
we integrate on $\la = \al_m + i\nu$.

To obtain the finite sum of residues in \eqref{eq4.79}, one integrates on the boundary of a square
in the $\la$-plane with sides parallel to the real and imaginary $\la$-axes, $\tilde u^2 - (2m+1)^2\pi^2<
\Ree \la < \tilde u^2 + \ve$, and $-\nu < \Imm \la < \nu$.
For fixed $u$ and $m$, $\Ree \la$ is in a finite interval, where $|d|$ increases exponentially as $|\nu|\to\infty$,
which makes the integrals on $\Imm\la =\pm \nu$ vanish at $|\nu|=\infty$.
Therefore,
\begin{align*}
&\frac{\G (n)t^{n-1}}{2} \int_{\Ree \la = \tilde u^2 +\ve} 
\frac{e^{\la / 2t} d\la / 2\pi i}{(\cosh \sqrt \la - \cos\vt_1 \cosh (u+i\vp_1))^n} \\
&= \sum_{k=-m}^m \frac{\G (n)t^{n-1}}{2} \oint_{\la_k}
\frac{e^{\la /2t} d\la / 2\pi i}{(\cosh \sqrt\la - \cos\vt_1 \cosh (u+i\vp_1))^n} + R_m ,
\end{align*}
where the integrals are on small circles about $\la_k$.
Now,
\[
\frac{\G(n)t^{n-1}}{2} \oint_{\la_k} \frac{e^{\la / 2t} d\la /2 \pi i }{(\cosh \sqrt\la - x)^n}= \frac{t^n\p^n}{\p x^n} e^
{\frac{(\cosh^{-1} x + i2k\pi)^2}{2t}}
= \left( \frac{t\p}{\sinh\kappa \p \kappa}\right)^n e^{\frac{(\kappa+i2k\pi)^2}{2t}} ,
\]
$x = \cosh \kappa$,
and we have justified \eqref{eq4.38} and derived \eqref{eq4.79}, \eqref{eq4.80}.
\end{proof}

\begin{484Corollary}
Assuming $0<\de\leq|\vp_1| \leq\pi - \de'$, $\de' >0$, $P_C$ of \eqref{eq4.40}
can be put in the following form
\addtocounter{equation}{1}
\begin{align}
\label{eq4.85}
P_C = &\frac{e^{\frac{n^2}{2} t}} {(2\pi t)^{n+1}} \sum_{k=-m}^{m} \int_{-\infty}^\infty 
e^{-\frac{u^2 -  (\kappa +i2k\pi)^2}{2t} }
V_n (\kappa + i2k\pi , t) du \\
&+ O\Big( C_{\de'} (2m+1)\pi t^{-2} e^{-\frac{(2m+1)^2\pi^2}{2t}}\Big) ,\nonumber
\end{align}
$t<1$, where $m=1,2,\ldots$ is large.
In particular, $P_C \to 0$ as $t\to 0$, uniformly in $(\vt_1 , \vp_1)$, $0\leq\vt_1 \leq\pi/2$ and $\de < |\vp_1| < \pi - \de'$.
\end{484Corollary}

\begin{proof}
\eqref{eq4.40}, \eqref{eq4.79} and \eqref{eq4.80} imply \eqref{eq4.85}.
As for the last statement, it suffices to show that each integral term vanishes as $t\to 0$.
Note that $V_{n,j}(z)$ is analytic and can be estimated as follows,
\begin{equation}
\label{eq4.86}
V_{n,j} (z) = O \left( \frac
{|z|^n}{(1+|z|^j) |\sinh z|^n} \right) .
\end{equation}
Consequently, with $\kappa = \tilde u + i \tilde \vp_1$, each integral term is dominated by
\[
e^{-\frac{\tilde\vp_1^2}{2t}} \int_{-\infty}^\infty \frac
{(1+ |\tilde u|)^n}{\sin^n \de'}
e^{-\frac{u^2-\tilde u^2}{2t}} du ,
\]
where $0 < \de < |\tilde\vp_1| < \pi - \de'$.
The integral is finite since $|\tilde u|\leq|u|$, $u^2 - \tilde u^2 \ge \mu |u|$, and the factor
$e^{-\tilde\vp_1^2 / 2t} < e^{-\de^2 / 2t}$ in front vanishes as $t\to 0$.
Hence Corollary 4.84.
\end{proof}

The integrands of the individual summands in \eqref{eq4.24} are singular when
$\cosh^2 \kappa = \cos^2 \vt_1 \cosh^2 (u+i\vp_1)=1$ which happens
at $\vp_1 = 0$, $\pi$ and $\cosh u = (\cos\vt_1)^{-1}$.
The singularity is integrable in the $u$-variable when $\vt_1 \ne 0$, but behaves
like $u^{-1}$, hence nonintegrable, when $\vt_1=0$.
Of course, the sum $P_C$ has no singularities at $\vt_1=0$, $\vp_1=0, \pi$; more later.

We still need to show that $P_C$ of \eqref{eq4.40} vanishes when $t\to 0$ at $\vp_1\sim \pm \pi$.
To this end one notes that the argument of Lemma 4.78 redone with the line of integration
$\Ree\la = \tilde u^2 - (2m\pi)^2 = \al_m$ yield the following analogue of Corollary 4.84:

\begin{487Lemma}
Assuming $0<\de \leq |\vp_1 | \leq \pi - \de'$, $\de'>0$, $P_C$ of \eqref{eq4.40} can be put in the following
form
\addtocounter{equation}{1}
\begin{align}
\label{eq4.88}
P_C = &\frac{e^{\frac{n^2}{2}t}}{(2\pi t)^{m+1}} \sum_{k=-m-1}^m \int_{-\infty}^\infty 
e^{-\frac{u^2-(\kappa + i2k\pi)^2}{2t}} V_n (\kappa + i2k\pi, t) du \\
&+ O \Big( C_\de 2m\pi t^{-2} e^{-\frac{(2m\pi)^2}{2t}}\Big) ,\nonumber
\end{align}
$t < 1$, and $m=1,2,\ldots$, large.
\end{487Lemma}

\begin{489Lemma}
$P_C$ of \eqref{eq4.40} converges to zero as $t\to 0$, uniformly in $\vt_1\sim 0$, $\vp_1\sim\pm\pi$.
\end{489Lemma}

\addtocounter{equation}{1}

\begin{proof}
On $S^{2n+1}$ one writes \eqref{eq4.88} as
\begin{equation}
\label{eq4.90}
P_C = \frac{e^{\frac{n^2}{2}t}}{(2\pi t)^{n+1}} \sum_{k=-m-1}^{m} \int_{-\infty}^\infty I_k^{(n)} + R_m .
\end{equation}
$R_m$ and the estimate on $R_m$ extend to $\vp_1 \sim \pm\pi$ and $R_m$ vanishes as $t\to 0$,
uniformly in $\vp_1 \sim \pm\pi$.
The rest of the argument amounts to showing that $I_k^{(n)} + I_{-k-1}^{(n)}$, $k=0,1,2,\ldots$ extend
to $\vp_1 \sim \pm\pi$, and vanish as $t\to 0$, uniformly in $\vp_1$, $\vp_1\sim\pm\pi$.
We start with $n=1$.
Set $\vp_1=\pi-\ve$, $\tilde\vp_1= \pi - \tilde\ve$.
Then
\begin{align}
\label{eq4.91}
&I_k^{(1)} + I_{-k-1}^{(1)} = -2 e^{-\frac{u^2-(\tilde u - i\tilde\ve)^2}{2t}} \frac{\tilde u - i\tilde\ve}{\sinh
(\tilde u - i\tilde\ve)} e^{-\frac{(2k+1)^2\pi^2}{2t} }\cdot \\
&\cdot \left\{ \cosh \frac{2(2k+1)\pi (\tilde \ve + i \tilde u)}{2t} -
\frac{2(2k+1)^2\pi^2}{2t}
\frac{\sinh \frac{2(2k+1)\pi(\tilde\ve + i\tilde u)}{2t}}{\frac{2(2k+1)\pi (\tilde\ve + i\tilde u)}{2t}}
\right\} ,\nonumber
\end{align}
and $I_k^{(1)} + I_{-k-1}^{(1)}$ is well defined at $\tilde\ve = 0$,
that is, at $\vp_1 = \pi$.
When $\vt_1\sim 0$ one may assume that $\tilde \ve < \pi / 8$.
Since $|a+ib|/|\sinh (a+ib)|$ is an increasing function of $|b|$, one has
\begin{align*}
&\left| \frac{\tilde u - i\tilde \ve}{\sinh (\tilde u - i\tilde\ve)}\right|^2 =
\frac{\tilde u^2 + \tilde\ve ^2}{\sinh^2 \tilde u + \sin^2 \tilde\ve} \leq
\frac{\tilde u^2 + \pi^2 / 64}{\sinh^2 \tilde u^2 + \sin^2 (\pi/8)} , \\
&\left| \frac{\sinh \frac{2(2k+1)\pi(\tilde\ve +i\tilde u)}{2t}} {\frac{2(2k+1)\pi(\tilde\ve + i\tilde u)}{2t}}
\right|^2 \leq
\frac{\sinh^2 \frac{2(2k+1)\pi\tilde\ve}{2t}} {\frac{2(2k+1)\pi\tilde\ve}{2t}}
\leq \cosh^2 \frac{2(2k+1)\pi\tilde\ve}{2t} \leq e^{\frac{(2k+1)\pi^2/2}{2t}} .
\end{align*}
Consequently,
\[
|I_k^{(1)} + I_{-k-1}^{(1)}| \leq C \frac{\tilde u}{\sinh\tilde u}
\left( 1 + \frac{2(2k+1)^2 \pi^2}{2t}\right) e^{-\frac{(2k+1)^2\pi^2/2}{2t}} ,
\]
which converges to $0$ as $t\to 0$.
Note that modulo $\exp \big( -u^2 / 2t - (2k+1)^2 \pi^2 / 2t\big)$, \eqref{eq4.91} has the following
form:
\begin{equation}
\label{eq4.92}
e^{\frac{(\tilde u - i\tilde\ve)^2}{2t}}
\left( \frac{\tilde u - i\tilde\ve}{\sinh (\tilde u - i\tilde\ve)}\right)
\psi_1 (\tilde\ve + i\tilde u) ,
\end{equation}
where each factor is an analytic function of $(\tilde u - i\tilde \ve)^2$.
To find $I_k^{(n)} + I_{-k-1}^{(n)}$ we use \eqref{eq4.38} on \eqref{eq4.92} and
note that $\p / \p \big( \cosh (\tilde u - i\tilde \ve)\big)\sim \p/\p (\tilde u - i\tilde\ve)^2$.
Thus such derivatives still yield an analytic function of $(\tilde u - i\tilde\ve)^2$.
Derivatives of the first factor lead to the first factor times the second factor.
Derivatives of the second factor give powers of the second factor times a polynomial
in $\coth (\tilde u - i\tilde \ve)$ which is bounded for large $u$, and the third factor
and its derivatives are dominated by
$\cosh \frac{2(2k+1)\pi\tilde\ve}{2t}$ and $\sinh \frac{2(2k+1)\pi\tilde\ve}{2t}$ with
$\tilde\ve < \pi/8$.
Therefore one finds
\begin{equation}
\label{eq4.93}
|I_{k}^{(n)} + I_{-k-1}^{(n)}| \leq C \left( \frac{\tilde u}{\sinh \tilde u}\right)^n
\left( 1 + \frac{2(2k+1)\pi}{2t}\right)^{n+1} e^{-\frac{(2k+1)^2\pi^2/2}{2t}} ,
\end{equation}
which is integrable on $-\infty < u < \infty$, and the integral vanishes as $t\to 0$.
\end{proof}

Individual terms in \eqref{eq4.24} do not extend to $\vt_1=0$, $\vp_1=0$, but the sum of
the $k$-th and $(-k)$-th term does: with $\kappa= \tilde u +i\tilde\vp_1$, one has
\begin{align}
\label{eq4.94}
I_k^{(1)} + I_{-k}^{(1)} =
&2e^{-\frac{u^2-\kappa^2}{2t}}
\frac{\kappa}{\sinh \kappa}
\sum_{k=0}^\infty
e^{-\frac{4k^2\pi^2}{2t}} \cdot \\
&\cdot\left\{
\cosh \frac{4k\pi (-i\kappa)}{2t} - \frac{2k^2\pi^2}{2t}
\frac{\sinh \frac{4k\pi (-i\kappa)}{2t}} {\frac{4k\pi (-i\kappa)}{2t}}
\right\} ;\nonumber
\end{align}
compare with \eqref{eq3.50}.
As above, one can extend this to $I_k^{(n)} + I_{-k}^{(n)}$.
This yields

\begin{495Theorem}
$t>0$.

{\rm 1) }
\addtocounter{equation}{1}
\begin{align}
\label{eq4.96}
P_C = \frac{e^{\frac{n^2}{2}t}} {(2\pi t)^{n+1}} \Bigg\{ &\int_{-\infty}^\infty
e^{-\frac{u^2-\kappa^2}{2t}} V_n (\kappa,t) du \\
& +\sum_{k=1}^\infty \int_{-\infty}^\infty \sum_{j=\pm k} e^{-\frac{u^2-(\kappa+i 2j\pi)^2}{2t}} V_n
(\kappa + i2j\pi,t)du\Bigg\} \nonumber
\end{align}
converges absolutely and uniformly in $0 \leq \vt_1 \leq \pi/2$, $0 \leq
|\vp_1 | \leq \pi - \de'$, $\de' > 0$.

{\rm 2)} 
\begin{equation}
\label{eq4.97}
P_C = \frac{e^{\frac{n^2}{2}t}}{(2\pi t)^{n+1}} \sum_{k=0}^\infty \int_{-\infty}^\infty
\sum_{k=j,-j-1} 
e^{-\frac{u^2-(\kappa+i2j\pi)^2}{2t}} V_n (\kappa + i2j\pi , t) du\nonumber
\end{equation}
converges absolutely and uniformly in  $0\leq \vt_1 \leq \pi/2$,
$0 < \de < |\vp_1 | \leq \pi$.
\end{495Theorem}

\section{CR geodesics on $S^{2n+1}\subset \BC^{n+1}$}
\label{sec5}

The parameters $\tau_1$ and $\Omega$ are defined by \eqref{eq4.13} and \eqref{eq4.16}.
Write \eqref{eq4.13} in the form
\begin{equation}
\label{eq5.1}
\frac{\sin^2\vt_1}{\sin^2\Om t} = 1 - \frac{\tau_1^2}{\Om^2} = \al^2 \leq 1 .
\end{equation}
This requires that $\vt_1 + k\pi \leq \Om t \leq (k+1)\pi - \vt_1$, $k=0,1,2,\ldots$, at least
when $\Om t>0$. From \eqref{eq5.1} one has $\tau_1=\tau_1(\Om)=\Om\sqrt{1-\al^2}$.
We set $\eta=\Om t$.
In view of \eqref{eq4.16} one may write
\begin{equation}
\label{eq5.2}
\vp_1(s;t,\Om,\tau_1) = \vp_1 \left( \frac{s}{t} ; 1 , \eta, \eta \sqrt{1-\al^2}\right) ,
\end{equation}
which vanishes when $s=0$ if and only if
\begin{equation}
\label{eq5.3}
\vp_1 \big( 0; 1, \eta, \eta\sqrt{1-\al^2}\big) = 0 .
\end{equation}
Given $\vp_1$, $0<\vp_1 \leq \pi$, let $\vp_1(1;1,\eta,\eta\sqrt{1-\al^2})=\vp_1$.
Then Lemma 4.15 implies that \eqref{eq5.3} is equivalent to
\begin{equation}
\label{eq5.4}
\vp_1 = -\eta \sqrt{1-\al^2} + T(\eta) ,\qquad \vt_1 > 0 ,
\end{equation}
where we set
\begin{equation}
\label{eq5.5}
T(\eta) = \int_0^\eta \frac {d( \sqrt{1-\al^2}\tan\s)}{1+ (\sqrt{1-\al^2}\tan\s)^2} .
\end{equation}
When $\vt_1=0$, one has $\eta = k\pi$, $k\in\BZ$,  and $\tau_1$ is given by
\begin{equation}
\label{eq5.6}
\vp_1 = -\tau_1 t +  T(\eta) =  -\tau_1 t + k\pi .
\end{equation}
Let \eqref{eq5.4} define $\vp_1(\eta)$.

\begin{57Lemma}
$\vp_1(\eta)$ is an increasing function of $\eta$ when $\vt_1\leq \eta \leq \pi/2$.
At the endpoints one has $\vp_1(\vt_1)=0$, $\vp_1(\pi/2)=\pi(1-\cos\vt_1)/2$.
\end{57Lemma}

\begin{proof}
Here $T(\eta)=\tan^{-1} (\sqrt{1-\al^2}\tan\eta)$.
With $\al'(\eta) = -\al\cot \eta$, we have
\addtocounter{equation}{1}
\begin{equation}
\label{eq5.8}
\frac{d}{d\eta} \tan^{-1} (\sqrt{1-\al^2}\tan \eta) = \frac{1}{\sqrt{1-\al^2}} ,
\end{equation}
hence
\begin{equation}
\label{eq5.9}
\vp'_1 (\eta) = \frac{\al^2}{\sqrt{1-\al^2}} (1-\eta \cot \eta) > 0 .
\qedhere
\end{equation}
\end{proof}

We note that
\begin{equation}
\label{eq5.10}
\frac{d}{d\eta} \eta\cot\eta = \frac{\sin 2\eta - 2\eta}{2\sin^2 \eta} \leq 0 ,
\end{equation}
and $\eta\cot\eta$ is a decreasing function of $\eta$ on $0<\eta<\pi/2$ with maximum $1$
at $\eta =0$.

\begin{511Corollary}
$\vp_1 (\eta)$ is an increasing function of $\eta$ on $\vt_1\leq \eta \leq \pi-\vt_1$ with
$\vp_1 (\pi-\vt_1)=\pi$
\end{511Corollary}

\begin{proof}
It suffices to assume that $\pi/2\leq\eta \leq \pi-\vt_1$. 
Set $\eta=\pi-\eta'$, so
$\vt_1 \leq\eta'\leq \pi/2$.
Then
\addtocounter{equation}{1}
\begin{align}
\label{eq5.12}
T(\eta) &= \frac{\pi}{2} + \int_{\pi/2}^{\pi-\eta'} \frac{d(\sqrt{1-\al^2}\tan\s)}{1+(\sqrt{1-\al^2}\tan\s)^2}\\
&=\frac{\pi}{2}+\int_{-\pi/2}^{-\eta'}\frac{d(\sqrt{1-\al^2}\tan\s)}{1+(\sqrt{1-\al^2}\tan\s)^2}\nonumber\\
&=\pi - \tan^{-1} (\sqrt{1-\al^2} \tan\eta') ,\nonumber
\end{align}
and
\begin{align}
\label{eq5.13}
\vp_1 (\eta) &= -(\pi-\eta') \sqrt{1-\al^2} + \pi-\tan^{-1} (\sqrt{1-\al^2} \tan \eta') \\
&=\pi (1-\sqrt{1-\al^2}) - \vp_1 (\eta') ,\nonumber
\end{align}
which is a decreasing function of $\eta'$ and therefore an increasing function of $\eta$.
\end{proof}

\begin{514Lemma}
Let $\eta_k$ denote the unique solution of $\tan\eta=\eta$ in the interval $k\pi\leq \eta \leq k\pi
+\pi/2$ and set $\eta_k = \eta'_k + k\pi$.

{\rm 1)} $\vt_1 < \eta'_k$.
On the interval $\vt_1 + k\pi \leq \eta \leq (k+1)\pi-\vt_1$ the curve $\vp_1(n)$ decreases from
$\vp_1(\vt_1 +k\pi)=k\pi$ to $\vp_1(\eta_k)\geq 0$ and then increases to $\vp_1
\big( (k+1)\pi-\vt_1\big)=(k+1)\pi$.

{\rm 2)} $\vt_1\geq \eta'_k$.
$\vp_1(\eta)$ increases from $\vp_1(\vt_1+k\pi)=k\pi$ to 
$\vp_1\big( (k+1)\pi-\vt_1\big) = (k+1)\pi$.
\end{514Lemma}

\begin{proof}
1) $\vt_1 < \eta'_k$. (i) $\vt_1+k\pi\leq \eta \leq \pi/2 + k\pi$.
With $\eta=\eta'+k\pi$,
\addtocounter{equation}{1}
\begin{align}
\label{eq5.15}
\vp_1(\eta) &=-\eta \sqrt{1-\al^2} + T(\eta) \\
&= -\eta \sqrt{1-\al^2} + k\pi + \tan^{-1} (\sqrt{1-\al^2} \tan\eta')\nonumber \\
&= \vp_1 (\eta') + k\pi ( 1 - \sqrt{1-\al^2}) \nonumber \\
&\geq 0 ,\nonumber
\end{align}
and
\begin{equation}
\label{eq5.16}
\vp'_1 (\eta) = \frac{\al^2}{\sqrt{1-\al^2}} \frac{\tan\eta-\eta}{\tan\eta} .
\end{equation}
$\tan\eta-\eta$ is an increasing function on $k\pi <\eta < k\pi + \pi/2$, it increases from
$-k\pi$ at $\eta=k\pi$ to $\infty$ at $\eta=k\pi + \pi/2$ with $\tan\eta_k - \eta_k =0$.
In particular,
\begin{equation}
\label{eq5.17}
\vp'_1 (\eta) = \begin{cases}
<0 , &\vt_1+k\pi < \eta < \eta_k,\\
>0 , &\eta_k < \eta < k\pi + \pi/2 .
\end{cases}
\end{equation}

(ii) Set $\eta = (k+1)\pi -\eta'$ on $k\pi + \pi/2 < \eta < k\pi +\pi$.
Then the argument of \eqref{eq5.12} yields
\begin{align}
\label{eq5.18}
T(\eta) &= k\pi + \frac{\pi}{2} + \int_{k\pi+\frac{\pi}{2}}^{-\eta'+(k+1)\pi}
\frac{d(\sqrt{1-\al^2}\tan\s)}{1+(\sqrt{1-\al^2}\tan\s)^2} \\
&= (k+1) \pi - \tan^{-1} (\sqrt{1-\al^2} \tan\eta') ,\nonumber
\end{align}
hence
\begin{equation}
\label{eq5.19}
\vp_1 (\eta) = (k+1)\pi ( 1 - \sqrt{1-\al^2}) - \vp_1 (\eta') .
\end{equation} 
Thus $\vp_1(\eta)$ is a decreasing function $\eta'$ hence an increasing function of $\eta$.

2) follows from \eqref{eq5.16}.
\end{proof}

Note that \eqref{eq5.18} yields $T(k\pi)=k\pi$,
see also \eqref{eq5.6}.
Furthermore
\begin{align}
\label{eq5.20}
&\vp_1 (\vt_1 + k\pi) = k\pi ,\\
\label{eq5.21}
&\vp_1 \left( k\pi + \frac{\pi}{2}\right) = \left( k + \frac12\right) \pi ( 1 - \cos\vt_1) ,\\
\label{eq5.22}
&\vp_1 \big( (k+1) \pi - \vt_1\big) = (k+1)\pi .
\end{align}
In particular, \eqref{eq5.15} implies
\begin{equation}
\label{eq5.23}
\lim_{\vt_1\to 0} \vp_1 (\eta_k) = 0.
\end{equation}

\begin{524Lemma}
On $\vt_1+k\pi \leq \eta \leq (k+1)\pi - \vt_1$ one has
\addtocounter{equation}{1}
\begin{equation}
\label{eq5.25}
\min_\eta \vp_1 (\eta+\pi) \geq \min_\eta \vp_1 (\eta) + \pi(1-\cos\vt_1) .
\end{equation}
\end{524Lemma}

\begin{proof}
\eqref{eq5.15} implies
\begin{align}
\label{eq5.26}
\vp_1(\eta+\pi) &= - (\eta+\pi) \sqrt{1-\al^2} + T(\eta+\pi) \\
&= - (\eta+\pi) \sqrt{1-\al^2} + \pi + T (\eta) \nonumber \\
&= \vp_1 (\eta) + \pi ( 1 - \sqrt{1-\al^2}) ,\nonumber
\end{align}
hence
\begin{equation}
\label{eq5.27}
\min_\eta \vp_1 (\eta+\pi) \geq \min_\eta \vp_1 (\eta) + \min_\eta \pi (1-\sqrt{1-\al^2})
\end{equation}
which is \eqref{eq5.25}.
\end{proof}

\begin{528Corollary}
When $\vt_1<\eta'_k$, \eqref{eq5.25} implies
\addtocounter{equation}{1}
\begin{equation}
\label{eq5.29}
\vp_1 (\eta_{k+1}) \geq \vp_1 (\eta_k) + \pi (1-\cos\vt_1) .
\end{equation}
In particular,
\begin{equation}
\label{eq5.30}
\lim_{k\to\infty} \vp_1 (\eta_k) = \infty ,\qquad \vt_1 > 0 .
\end{equation}
\end{528Corollary}

We have derived

\begin{531Proposition}
Let $0\leq \vp_1 \leq \pi$ and $\vt_1 \leq \eta' \leq \pi-\vt_1$.  
Then

{\rm (i)} $\vp_1 (\eta')=\vp_1$ has a unique solution.

{\rm (ii)} When $\vt_1>0$, $\vp_1 (\eta'+k\pi)=\vp_1$ has solutions for at most a finite number of $k$'s among
$k=1,2,\ldots$
If it has a solution for a particular $k$ then it has exactly two solutions with that $k$,
$\eta_{k,1} \le \eta_k \le \eta_{k,2}$; as a solution $\eta_k$ is
counted twice.

{\rm (iii)} When $\vt_1=0$, $\vp_1(\eta'+k\pi)=\vp_1$ has two solutions,
$\eta' = 0,\pi$, for each $k=1,2,\ldots$; see \eqref{eq5.6}.
\end{531Proposition}

Note that $\vp_1(\eta)$ can be given by a single expression on
$\vt_1 \leq \eta \leq \pi-\vt_1$ if we use
\addtocounter{equation}{1}
\begin{equation}
\label{eq5.32}
\cot^{-1} \left( \frac{\cot\eta}{\sqrt{1-\al^2}}\right) = \begin{cases}
\tan^{-1} (\sqrt{1-\al^2} \tan \eta),&\eta < \frac{\pi}{2} ,\\
\tan^{-1} (\sqrt{1-\al^2}\tan\eta) + \pi,&\eta>\frac{\pi}{2} ,
\end{cases}
\end{equation}
since $\cot^{-1} \left( \frac{\cot\eta}{\sqrt{1-\al^2}}\right)$ is continuous on $\vt_1<\eta<\pi-\vt_1$.

\begin{533Lemma}
Let $\vt_1+k\pi < \eta < (k+1)\pi - \vt_1$.
Then
\addtocounter{equation}{1}
\begin{equation}
\label{eq5.34}
\vp_1(\eta) 
=\vp_1 (\eta-k\pi) + k\pi (1-\sqrt{1-\al^2}) .
\end{equation}
\end{533Lemma}

\begin{proof}
In view of \eqref{eq5.15} we only need to derive \eqref{eq5.34} for $\eta$ in $k\pi + \pi/2 < \eta < (k+1)\pi -\vt_1$.
Here \eqref{eq5.19} and \eqref{eq5.12} imply
\begin{align*}
\vp_1(\eta) &= (k+1)\pi (1-\sqrtonealtwo) - \vp_1 \big( (k+1)\pi - \eta\big) ,\\
&= (k+1)\pi (1-\sqrtonealtwo) + \big( (k+1)\pi - \eta\big) \sqrtonealtwo \\
&\quad - \tan^{-1} \big( \sqrtonealtwo \tan \big( \pi - (\eta - k\pi) \big) \\
&= -\eta\sqrtonealtwo + k\pi + T (\eta - k\pi) ,
\end{align*}
which is \eqref{eq5.34}.
\end{proof}

We are interested in the number of CR-geodesics, projections of CR-bicharacteristics on $S^{2n+1}$,
between two points and their lengths; rotation invariance allows us to assume that the geodesic
starts at $1=(1,0,\ldots,0)$ and ends at $Q=(z_1,z_2,\ldots, z_{n+1})$ with $\vt_1(Q)=\vt_1$,
$\vp_1(Q)=\vp_1$.
This requires parameters $\tau_1$ and $\eta=\Om t$ which are solutions of \eqref{eq5.1} and \eqref{eq5.4}.
Since $z_2,\ldots, z_{n+1}$ depend on $\sin\vt_1$, $\eta$ must belong to intervals with even $k$,
$\vt_1+2k\pi \leq \eta \leq (2k+1)\pi -\vt_1$.
Then $\vp_1(0)\equiv 0$ (mod $2k\pi$) and this leaves the starting point $1$ undisturbed.
The following result is a translation of Proposition 5.31 into the language of geodesics.

\begin{535Proposition}
{\rm (i)} Given $Q\in S^{2n+1}$, there are an odd number of CR-geodesics when $\vt_1(Q)>0$,
at least one, and an infinite number of CR-geodesics when $\vt_1(Q)=0$ which connect $1$ and $Q$.

{\rm (ii)} When $\vt_1(Q)\geq \eta_1$, there is exactly one CR-geodesic between $1$ and $Q$.
\end{535Proposition}

\begin{proof}
(ii) is a consequence of Lemmas 5.14 and 5.24.
\end{proof}

The geodesics in question are local geodesics in the following sense.

\begin{536Lemma}
$\vp_1 (Q)(s)$ is an increasing function of $s$ along a CR-geodesic of Proposition 5.35 between $1$ and $Q$.
\end{536Lemma}

\begin{proof}
\eqref{eq4.16} and its periodic extension yields
\addtocounter{equation}{2}
\begin{equation}
\label{eq5.37}
\vp'_1(s) = \frac{\tau_1 \left( 1 - \frac{\tau_1^2}{\Om^2}\right) \tan^2 \Om s }{1+\left(\frac{\tau_1}{\Om}\tan\Om s\right)^2}
> 0 .\qedhere
\end{equation}
\end{proof}

The number of CR-geodesics joining $1$ and $Q$ can be parametrized by a finite $(\vt_1>0)$ or an
infinite $(\vt_1=0)$ sequence of $\eta' s$.
The first $\eta$ is $\eta_{0,2}$, $\vt_1\leq \eta_{0,2} \leq \pi-\vt_1$.
If there is a solution to \eqref{eq5.4} in $\vt_1+2\pi \leq \eta \leq 3\pi-\vt_1$, then there are two solutions
$\vt_1+2\pi \leq \eta_{2,1} \leq \eta_2 \leq \eta_{2,2} \leq 3\pi - \vt_1$, etc.
So the sequence parametrizing the CR-geodesics is $\eta_{0,2} < \eta_{2,1} \leq \eta_{2,2} < \cdots < \eta_{2k,1}\leq
\eta_{2k,2}< \cdots$; there is no $\eta_{0,1}$.
We shall show that the lengths of the CR-geodesics which connect $1$ and $Q$ is an increasing function of
the attached parameter.
To this end we begin with the definition of CR-curves, which have lengths, and the demonstration that
CR-geodesics are CR-curves.

\begin{538Definition}
The Cauchy-Riemann (CR) tangent space of $S^{2n+1}$ is the complex orthogonal complement
of $N$ in $\BC^{n+1}$ restricted to $S^{2n+1}$.
A CR-curve on $S^{2n+1}$ has all its tangents in the CR-tangent space.
\end{538Definition}

In CR-geometry only CR-curves have length, which is then the usual Euclidean length.
According to the definition a curve $\big( z_1(s), \ldots, z_{n+1}(s)\big) \subset S^{2n+1}$ is a CR-curve if
and only if
\addtocounter{equation}{1}
\begin{equation}
\label{eq5.39}
\dz_1 (s) \oz_1(s) + \cdots + \dz_{n+1} (s) \oz_{n+1}(s) = 0 .
\end{equation}

\begin{540Lemma}
$\big( \vt(s),\vp(s)\big)\subset S^{2n+1}$ is a CR-curve if and only if
\addtocounter{equation}{1}
\begin{align}
\label{eq5.41}
&\dvp_1\cos^2\vt_1 + \sin^2\vt_1 \big( \dvp_2 \cos^2 \vt_2 + \sin^2 \vt_2 (\dvp_3 \cos^2 \vt_3 + \cdots \\
&\quad \cdots + \sin^2 \vt_{n-1} (\dvp_n \cos^2 \vt_n + \dvp_{n+1} \sin^2 \vt_n) \cdots\big) = 0 .\nonumber
\end{align}
\end{540Lemma}

\begin{proof}
We write \eqref{eq5.39} in spherical coordinates and find that the coefficient of $\dvt_{n}(s)$ vanishes, then we notice
that the coefficient of $\dvt_{n-1}(s)$ vanishes, and, continuing in this manner, we see that the coefficients
of $\dvt_n(s),\dvt_{n-1}(s) , \ldots, \dvt_1(s)$ all vanish.
This yields \eqref{eq5.41}.
\end{proof}

\begin{542Proposition}
A CR-geodesic $z(s) \subset S^{2n+1}$, $0\leq s \leq t$, is a CR-curve.
Its length is $\sqrt{2H}t = Et$, where $H$ is given by \eqref{eq4.9}.
\end{542Proposition}

\begin{proof}
We may assume that $z(0)=(1,0,\ldots,0)$, set $z(1)\cdot\oz(0)=\cos\vt_1 e^{i\vp_1}$ and extend
$(\vt_1,\vp_1)$ to a set of spherical coordinates.
Substituting \eqref{eq4.3} and \eqref{eq4.11} into \eqref{eq5.41} the sum on the left telescopes to $0$
and the geodesic is a CR-curve.
Its length is calculated using the Euclidean metric.
According to the CR-version of \eqref{eq3.19}, $\dvt_j = 0$, $j=2,\ldots, n$, and this implies
\begin{align*}
&\dz_1 = e^{i\vp_1} (-\dvt_1\sin\vt_1 + i \dvp_1 \cos\vt_1 ) ,\\
&\qquad\qquad ....\\
&\dz_k = e^{i\vp_k} (\dvt_1\cos\vt_1 + i\dvp_k \sin\vt_1) \sin\vt_2 \ldots \sin\vt_{k-1}\cos\vt_k,\quad
k=2,\ldots,n , \\
&\qquad\qquad ....\cr
&\dz_{n+1} = e^{i\vp_{n+1}} (\dvt_1 \cos\vt_1 + i\dvp_{n+1} \sin\vt_1) \sin\vt_2 \ldots \sin\vt_n .
\end{align*}
Since $\dvp_j = -\tau_1$, $j=2,\ldots, n+1$, see \eqref{eq4.3}, one has
\begin{equation*}
|\dz_k|^2 + \cdots + |\dz_{n+1}|^2 = (\dvt_1^2 \cos^2 \vt_1 + \tau_1^2 \sin^2 \vt_1 ) \sin^2\vt_2 \ldots
\sin^2\vt_{k-1} ,
\end{equation*}
$k=2,\ldots, n$, which yields
\begin{align*}
|\dz_1|^2 + \cdots + |\dz_{n+1}|^2 &=\dvt_1^2 \sin^2 \vt_1 + \dvp_1^2 \cos^2 \vt_1
+ \dvt_1^2 \cos^2 \vt_1 + \tau_1^2 \sin^2\vt_1 \\
&= \dvt_1^2 +\tau_1^2 \tan^2 \vt_1 \\
&=2H ,
\end{align*}
as $\dvp_1 = \tau_1 \tan^2 \vt_1$, see \eqref{eq4.11}.
$H$ is constant along $z(s)$, hence
\addtocounter{equation}{1}
\begin{equation}
\label{eq5.43}
\ell = \int_0^t \sqrt{2H} ds = \sqrt{2H} t = Et .
\qedhere
\end{equation}
\end{proof}

By definition
\begin{equation}
\label{eq5.44}
(Et)^2 = (\Om t)^2 - (\tau_1 t)^2 = \eta^2 - (\tau_1 t)^2 .
\end{equation}
$Et$ as a function of $\eta$ is defined in all the intervals $\vt_1+k\pi\leq \eta \leq (k+1)\pi-\vt_1$,
$k=0,1,2,\ldots$
Let $E_{0,2}, E_{1,1}, E_{1,2}, \ldots$ denote the value of $Et$ when $\eta$ is $\eta_{0,2},
\eta_{1,1}, \eta_{1,2},\ldots$

(i) $\vt_1=0$.  Here
\begin{equation}
\label{eq5.45}
\eta_{k,1} = k\pi ,\qquad k=1,2,\ldots ,
\end{equation}
\begin{equation}
\label{eq5.46}
\eta_{k,2} = (k+1)\pi ,\qquad k=0,1,2,\ldots ,
\end{equation}
so $\eta'_{k,j}=0$.
$\al$ is not defined, but \eqref{eq5.6} yields
\begin{equation}
\label{eq5.47}
\tau_1 t = \begin{cases}
k\pi - \vp_1 ,&\eta=\eta_{k,1} ,\\
(k+1)\pi - \vp_1 ,&\eta=\eta_{k,2}.
\end{cases}
\end{equation}
Consequently, one has
\begin{equation}
\label{eq5.48}
E_{k,2}^2 = E_{k+1,1}^2 = (k+1)^2 \pi^2 - \big( (k+1)\pi - \vp_1\big)^2 ,
\end{equation}
$k=0,1,2,\ldots$, and we have derived

\begin{549Lemma}
$\vt_1 =0$. When $0<\vp_1 \leq \pi$,
\addtocounter{equation}{1}
\begin{equation}
\label{eq5.50}
E_{k,2}^2 = E_{k+1,1}^2 = \big( 2(k+1) \pi -\vp_1\big)\vp_1 ,\qquad k=0,1,2,\ldots
\end{equation}
\end{549Lemma}

(ii)
Let $Q\in S^{2n+1}$ with $\vt_1 (Q) > 0$.
There is a finite number of solutions $\eta_{0,2},\eta_{1,1},\eta_{1,2}$, \ldots of $\vp_1(\eta)=\vp_1$.
To each $\eta_{k,j}$ we attach $E_{k,j}$ by 
\begin{equation}
\label{eq5.51}
E_{k,j}^2 = \eta_{k,j}^2 - \eta_{k,j}^2 (1-\al_{k,j}^2) = \sin^2\vt_1 \frac{\eta_{k,j}^2}{\sin^2\eta_{k,j}} ,
\end{equation}
and we shall show that $E_{k,j}$ increases with increasing $\eta_{k,j}$.
Start with

\begin{552Lemma}
Let $\vt_1 > 0$, $0<\vp_1\leq \pi$ and $\vp_1=\vp_1 (\eta_{k,1})=\vp_1(\eta_{k,2})$
for some $k$.
Then
\addtocounter{equation}{1}
\begin{equation}
\label{eq5.53}
E_{k,1}^2 < E_{k,2}^2 .
\end{equation}
\end{552Lemma}

\begin{proof}
Suppose $\vp_1 (\eta_-)=\vp_1 (\eta_+)$, $\vt_1+k\pi \leq \eta_- \leq \eta_k \leq \eta_+ \leq
(k+1)\pi - \vt_1$, and set $\eta=\eta' + k\pi$.
Then \eqref{eq5.34} yields
\begin{equation}
\label{eq5.54}
\vp_1 (\eta'_-) + k\pi \Big(1- \sqrtonealminustwo\Big) =\vp_1(\eta'_+) + k\pi \Big(1- \sqrtonealplustwo\Big) ,
\end{equation}
hence
\begin{equation*}
0 < \vp_1 (\eta'_+) - \vp_1 (\eta'_-) = k\pi \Big(\sqrtonealplustwo - \sqrtonealminustwo \Big) ,
\end{equation*}
and
\begin{equation}
\label{eq5.55}
\sqrtonealminustwo < \sqrtonealplustwo .
\end{equation}
Set
\begin{equation}
\label{eq5.56}
\psi(\eta) = \sin^2 \vt_1 \frac{\eta^2}{\sin^2\eta} ,
\end{equation}
and note that
\begin{equation}
\label{eq5.57}
\psi' (\eta) = \sin^2 \vt_1 \frac{2\eta}{\sin^2\eta} (1-\eta \cot\eta ) = 2\eta \sqrtonealtwo \vp'_1(\eta) .
\end{equation}
In particular, $\psi(\eta)$ decreases from $\psi(\eta_{k,1})$ to $\psi(\eta_k)$ then increases
to $\psi(\eta_{k,2})$.
We claim that $\psi(\eta_{k,1}) < \psi (\eta_{k,2})$.
Indeed,
\begin{align*}
\psi (\eta_{k,1}) - \psi(\eta_k) &= \int_{\eta_k}^{\eta_{k,1}} \psi' (\eta_-) d\eta_- \\
&= \int_{\eta_k}^{\eta_{k,1}} 2\eta_- \sqrtonealminustwo \vp'_1 (\eta_-)d\eta_- \\
&= \int_{\vp_1 (\eta_k)}^{\vp_1} 2\eta_- \sqrtonealminustwo d\vp_1 ,
\end{align*}
and similarly,
 \begin{equation*}
\psi(\eta_{k,2}) - \psi (\eta_k) = \int_{\vp_1 (\eta_k)}^{\vp_1} 2\eta_+ \sqrt{1-\al_+^2} d\vp_1 .
\end{equation*}
Consequently,
\[
\psi(\eta_{k,2}) - \psi(\eta_{k,1}) = \int_{\vp_1(\eta_k)}^{\vp_1}
\Big( 2 \eta_+ \sqrtonealplustwo - 2\eta_- \sqrt{1 - \al_-^2}\Big) d\vp_1 > 0
\]
since \eqref{eq5.55} and $\eta_- <\eta_k <\eta_+$ imply that the integrand is positive.
\end{proof}

\begin{558Lemma}
Fix $\vp_1$, $0<\vp_1 \leq \pi$, and choose $\eta$ in
$k\pi + \pi/2 \leq \eta \leq (k+1)\pi$.
Let \eqref{eq4.16} and \eqref{eq4.13}, with $s=t$ and $\vp_1 (0) = 0$, define
$\tau_1=\tau_1 (\eta)$ and $\vt_1 = \vt_1 (\eta)$  and set $S^2=S^2\big( \vt_1 (\eta) ,\eta\big)=\eta^2-
(\tau_1 t)^2$.
Then
\addtocounter{equation}{1}
\begin{equation}
\label{eq5.59}
\frac{dS^2}{d\eta} > 0.
\end{equation}
\end{558Lemma}

\begin{proof}
Let $\g = \g(\eta)=\sqrtonealtwo$.
Then \eqref{eq5.19} implies 
\begin{equation}
\label{eq5.60}
\eta\g = (k+1)\pi - \vp_1 + \tan^{-1} (\g \tan \eta' ) ,
\end{equation}
or,
\begin{equation}
\label{eq5.61}
\g\tan\eta' =
\g\tan \eta = \tan (\vp_1 + \eta\g ) ,
\end{equation}
and the $\eta$-derivative of \eqref{eq5.61} yields
\begin{equation}
\label{eq5.62}
\g' = \frac{\g(1-\g^2)\tan^2\eta}{\eta - \tan\eta +\eta\g^2 \tan^2\eta} <
\frac{1-\g^2}{\eta\g} ,
\end{equation}
if we drop $\eta-\tan\eta > 0$ from the denominator.
Now
\begin{equation}
\label{eq5.63}
S^2 = \eta^2 - (\tau_1 t)^2 = \eta^2 (1-\g^2) ,
\end{equation}
and
\begin{equation*}
\frac{d S^2}{d\eta} = 2\eta (1-\g^2 - \eta \g \g') > 2\eta \left( 1 - \g^2 - \eta\g
\frac{1-\g^2}{\eta\g}\right) = 0 .\qedhere
\end{equation*}
\end{proof}

\begin{564Corollary}
$\vp_1$ is fixed, $0 < \vp_1 \leq \pi$.
Then
\addtocounter{equation}{1}
\begin{equation}
\label{eq5.65}
\lim_{\eta\to (k+1)\pi}\, S(\eta)^2 = \big( 2 (k+1)\pi - \vp_1\big) \vp_1 .
\end{equation}
\end{564Corollary}

\begin{proof}
$(1-\al^2)\tan^2\eta = (\sin^2 \eta - \sin^2\vt_1) / \cos^2 \eta \to 0 $ as $\eta\to (k+1)\pi$,
since $\vt_1 \leq (k+1)\pi - \eta\to 0$.
Consequently, in view of \eqref{eq5.60}
\begin{align*}
\lim_{\eta\to (k+1)\pi} \tau_1t &= (k+1)\pi - \vp_1 + \lim_{\eta\to (k+1)\pi} \tan^{-1}
\Big( \sqrtonealtwo \tan\eta' \Big) \\
&= (k+1)\pi - \vp_1 ,
\end{align*}
since $\tau_1 t = \eta\g$, and
\begin{align*}
\lim_{\eta\to (k+1)\pi} S(\eta)^2 &= \lim_{\eta\to (k+1)\pi} \big( \eta^2 - (\tau_1 t)^2\big) \\
&= (k+1)^2 \pi^2 - \big( (k+1)\pi - \vp_1\big)^2 .\qedhere
\end{align*}
\end{proof}

Note the agreement with Lemma 5.49.
\begin{566Lemma}
$\vp_1$ is fixed, $0< \vp_1 \leq \pi$, and $k\pi < \eta < k\pi + \pi/2$.
Then
\addtocounter{equation}{1}
\begin{equation}
\label{eq5.67}
\frac{dS(\eta)^2}{d\eta} = \begin{cases}
>0,&k\pi < \eta <\eta_k,\\
<0,&\eta_k<\eta<k\pi + \frac{\pi}{2} .
\end{cases}
\end{equation}
\end{566Lemma}

\begin{proof}
With $\g=\sqrtonealtwo$, \eqref{eq5.60} implies \eqref{eq5.61}
and the equality in \eqref{eq5.62}.
Since
\begin{equation}
\label{eq5.68}
\eta - \tan\eta = \begin{cases}
>0,&k\pi<\eta<\eta_k ,\\
<0,&\eta_k<\eta<k\pi+\frac{\pi}{2} ,
\end{cases}
\end{equation}
dropping $\eta-\tan\eta$ in \eqref{eq5.62} gives us
\begin{equation}
\label{eq5.69}
\eta\g' = \begin{cases}
> \frac{1-\g^2}{\g},&k\pi < \eta < \eta_k ,\\
< \frac{1-\g^2}{\g},&\eta_k < \eta < k\pi + \frac{\pi}{2} .
\end{cases}
\end{equation}
Substituting \eqref{eq5.69} into
\begin{equation}
\label{eq5.70}
\frac{dS^2}{d\eta} = 2\eta \big( 1 - \g (\g + \eta\g')\big)
\end{equation}
justifies \eqref{eq5.67}.
\end{proof}

Lemmas 5.52, 5.58 and 5.66 and Corollary 5.64 show that for a given $\vt_1$ and $\vp_1$, $S(\eta_{k,1})<S(\eta_{k,2})$, and
when $k\pi + \pi/2 < \eta_{k,2}$ one also has $S(\eta_{k,2})<S(\eta_{k+1,1})$.
We still need to show that $S(\eta_{k,2}) < S(\eta_{k+1,1})$ when $\eta_k < \eta_{k,2}<k\pi + \pi/2$.
Set
\begin{equation}
\label{eq5.71}
\eta_{k,2}^* = \begin{cases}
k\pi + \frac{\pi}{2} ,&\vp_1 \left( k\pi + \frac{\pi}{2}\right) \leq \pi ,\\
\max\limits_{\vp_1 (\eta_{k,2})\leq\pi} \eta_{k,2},&\vp_1 \left( k\pi + \frac{\pi}{2}\right) > \pi .
\end{cases}
\end{equation}

\begin{572Lemma}
\addtocounter{equation}{1}
$\vt_1 >0$.
Let $\vp_1 (\eta_{k+1})\leq \vp_1 (\eta_{k,2}) = \vp_1 (\eta_{k+1,1}) < \vp_1 (\eta_{k,2}^*)$.
Then
\begin{equation}
\label{eq5.73}
S(\eta_{k,2}) < S(\eta_{k+1,1}) .
\end{equation}
\end{572Lemma}

\begin{proof}
When $\vp_1 (\eta_{k+1,1})=\vp_1 (\eta_{k+1,2})$, one has $\eta_{k+1,2} < (k+1)\pi + \pi/2$
and $\vp_1 (\eta_{k+1,2} -\pi) < \vp_1 (\eta_{k+1,2})$
in view of \eqref{eq5.26}.
Since $\eta_k < \eta_{k,2} < \eta_{k,2}^*\leq k\pi + \pi/2$, \eqref{eq5.26} yields
\begin{equation*}
k\pi + \frac{\pi}2 - \eta_{k,2} < (k+1)\pi + \frac{\pi}{2} -\eta_{k+1,2} < (k+1)\pi +\frac{\pi}{2} -\eta_{k+1,1} ,
\end{equation*}
so $\sin^2 \eta_{k,2} > \sin^2 \eta_{k+1,1}$ and
\[
S(\eta_{k+1,1})^2 = \sin^2 \vt_1 \frac{\eta_{k+1,1}^2}{\sin^2\eta_{k+1,1}} 
> \sin^2 \vt_1 \frac{\eta_{k,2}^2}{\sin^2\eta_{k,2}} = S(\eta_{k,2})^2 .\qedhere
\]
\end{proof}

\begin{574Lemma}
Let $\vp_1(\eta_{0,2})=\vp_1 (\eta_{1,1})$, $0<\eta_{0,2}<\pi<\eta_{1,1}<\eta_1$.
Then
\addtocounter{equation}{1}
\begin{equation}
\label{eq5.75}
S(\eta_{0,2}) < S(\eta_{1,1}) .
\end{equation}
\end{574Lemma}

\begin{proof}
Lemmas 5.58 and 5.66 demonstrate the truth of \eqref{eq5.75} when $\pi/2 \leq \eta_{0,2} < \pi$,
hence a proof is required only when $0<\eta_{0,2}<\pi/2$.
In that case $\vp_1 (\eta_{0,2})=\vp_1 (\eta_{1,1})=\vp_1 (\eta'_{1,1}) +\pi (1-\sqrtonealtwo)$, therefore
$\eta'_{1,1}<\eta_{0,2} < \pi/2$, and $\sin^2 \eta_{1,1} = \sin^2 \eta'_{1,1} < \sin^2 \eta_{0,2}$.
Consequently,
\begin{equation*}
S(\eta_{1,1})^2 = \sin^2 \vt_1 \frac{\eta_{1,1}^2}{\sin^2\eta_{1,1}} > \sin^2 \vt_1
\frac{\eta_{0,2}^2}{\sin^2\eta_{0,2}} = S(\eta_{0,2})^2 .\qedhere
\end{equation*}
\end{proof}

We have derived

\begin{576Proposition}
Given $\vp_1$, $0\leq \vp_1 \leq \pi$, let $0\leq \eta_{0,2}\leq \eta_{1,1} \leq \eta_{1,2}\leq \cdots$
denote the solutions of $\vp_1(\eta) =\vp_1$.
Then $S(\eta_{k,j}) = E_{k,j}$, $j=1,2$, $k=0,1,2,\ldots$, is an increasing function of $\eta_{k,j}$.
\end{576Proposition}

\begin{577Corollary}
Let $k$ denote an even positive integer and let $S(\eta_{k,j})$ denote the length of the CR-geodesic
connecting $1$ and $Q=(z_1,\ldots, z_{n+1})\in S^{2n+1}$ which is parametrized by $\eta_{k,j}$.
Then $S(\eta_{k,j})$ is an increasing function of $\eta_{k,j}$.
In particular the shortest geodesic between $1$ and $Q$ has length $S(\eta_{0,2})$, usually referred to as the
Carnot-Caratheodory distance between $1$ and $Q$.
\end{577Corollary}

We still need to consider points $Q\in S^{2n+1}$ whose $\vp_1$-component is negative,
$-\pi<\vp_1(Q)<0$.
To this end we note that $\vp_1(\eta)$ is skew-symmetric,
\addtocounter{equation}{2}
\begin{equation}
\label{eq5.78}
\vp_1(-\eta) = - \vp_1(\eta) .
\end{equation}
Consequently, all the results of \S\ref{sec5} derived for $0<\vp_1\leq\pi$ hold when $-\pi<\vp_1<0$ if we
change $k$ to $-k$, $\eta$ to $-\eta$ and $\vp_1$ to $-\vp_1$.

\section{On the small time asymptotic of $P_C$}
\label{sec6}

We shall derive a formula for $\lim_{t\to 0} P_C$ in terms of the Carnot-Caratheodory distance which
is the analogue of the classical formula for $\lim_{t\to 0} P_S$ given in terms of the usual Riemannian
distance.
Since $P_C$ is given by integrals of exponentials we shall employ the stationary phase technique.
The calculations are based on the following two observations:

1) the values of the action function $f$ at its imaginary critical points represent geodesic lengths, and

2) all critical points of $f$ are imaginary.

1) and 2) are invariant under translations of $f(w)$, $w=u+iv$, along the imaginary axis and it is
convenient to work with $2g(w) = f(w-i\vp_1)$.
Let
\begin{equation}
\label{eq6.1}
g_k (w) = \frac12 (w-i\vp_1)^2 - \frac12 \big( \cosh^{-1} (\cos\vt_1 \cosh w)+i2k\pi\big)^2 ,
\end{equation}
$w\in \BC$, $k\in Z$; we shall usually suppress the subscript in $g_0$ and write $g_0=g$.
Recall that the projection of a bicharacteristic on $S^{2n+1}$ is a geodesic if $\vp_1(0)=0$.
Lemmas 4.12 and 4.15 yield
\begin{equation}
\label{eq6.2}
\vp_1 (0) = \vp_1 - \tan^{-1} \left( \sqrt{1 - \frac{\sin^2\vt_1}{\sin^2\Omega t}}\, \tan\Omega t\ \right) 
+ \Omega t \sqrt{1 - \frac{\sin^2\vt_1}{\sin^2\Omega t}}\  ,
\end{equation}
$\vt_1 \leq \Omega t \leq \pi - \vt_1$.
With
\begin{equation}
\label{eq6.3}
\tan^{-1} x = \cos^{-1} \frac{1}{\sqrt{1+x^2}} ,
\end{equation}
\begin{equation}
\label{eq6.4}
x = \sqrt{1 - \frac{\sin^2 \vt_1}{\sin^2 \Omega t}} \, \tan\Omega t ,\qquad 1+x^2 = \frac
{\cos^2\vt_1}{\cos^2 \Omega t} ,
\end{equation}
one has
\begin{align}
\label{eq6.5}
\vp_1 (0) &= \vp_1 - \cos^{-1}\frac{\cos\Omega t}{\cos\vt_1} + \Omega t \sqrt{1- \frac{\sin^2\vt_1}
{\sin^2 \Omega t}} \\
&= - (\psi - \vp_1) + \frac{\cos^{-1} (\cos\vt_1 \cos\psi)}{\sqrt{1-\cos^2\vt_1 \cos^2\psi}}\, \cos\vt_1\sin\psi ,\nonumber
\end{align}
where we set
\begin{equation}
\label{eq6.6}
\cos\vt_1 \cos\psi = \cos\Omega t ,\qquad 0\leq \psi \leq \pi .
\end{equation}
Differentiating $g =g_0$ of \eqref{eq6.1} with respect to $w$ one obtains
\begin{equation}
\label{eq6.7}
ig_w (i\psi) = -(\psi-\vp_1) + \frac{\cos^{-1}(\cos\vt_1 \cos\psi)}{\sqrt{1-\cos^2\vt_1 \cos^2\psi}}\, \cos\vt_1\sin\psi ,
\end{equation}
and a comparison of \eqref{eq6.5} and \eqref{eq6.7} yields

\begin{68Theorem}
Let $\Omega$ induce $\psi$ by \eqref{eq6.6},
$\vt_1 \leq \Omega t \leq \pi - \vt_1$, $0\leq\psi\leq \pi$.
Then
\addtocounter{equation}{1}
\begin{equation}
\label{eq6.9}
ig_w (i\psi) = \vp_1 (0; t,\Omega,\tau_1),\qquad
\tau_1 = \Omega \sqrt{1 - \frac{\sin^2 \vt_1}{\sin^2 \Omega t}}\ ,
\end{equation}
see \eqref{eq4.13}.
In particular, a geodesic which connects $(1,0,\ldots, 0)$ and $(\vt,\vp)$, $0<\vp_1 < \pi$, induces
a unique imaginary critical point $i\psi$ of $g$ via \eqref{eq6.6}, and vice versa.
\end{68Theorem}

A direct proof of the uniqueness of the critical point $i\psi\in i\BR$ when $\psi$ is in the interval
$0\leq \psi \leq \pi$ follows from
\begin{equation}
\label{eq6.10}
\frac{d}{d\psi} \left( \psi - \frac{\cos^{-1}(\cos\vt_1 \cos\psi)}{\sqrt{1-\cos^2\vt_1 \cos^2\psi}}\, \cos\vt_1
\sin\psi\right)
= \al^2 (1-\eta \cot\eta )
> 0 ,
\end{equation}
in the terminology of \S\ref{sec5}; compare with \eqref{eq5.9}.

According to \eqref{eq5.44} and Proposition 5.76 the minimum of geodesic lengths,
or Carnot-Caratheodory distance $d_c$ between $(\vt,\vp)$ and $(1,0,\ldots,0)$, is
\begin{equation}
\label{eq6.11}
d_c^2 = E_{0,2}^2 = \eta_{0,2}^2 - \eta_{0,2}^2 \left( 1 - \frac{\sin^2 \vt_1}{\sin^2 \eta_{0,2}}\right),\qquad
\vt_1 \ne 0 ,
\end{equation}
where $\eta_{0,2} = \Omega_{0,2} t$ is the smallest positive solution of $\vp_1 = \vp_1 (\eta)$, assuming
$\vp_1 >0$.
\eqref{eq6.5} implies that if $i\psi$ is the first imaginary critical point of $g$ and $\cos\Omega t = \cos\vt_1
\cos\psi$ then $\Omega = \Omega_{0,2}$,
\begin{equation}
\label{eq6.12}
\Omega t \sqrt{ 1 - \frac{\sin^2 \vt_1}{\sin^2 \Omega t}} \ =\ \cos^{-1} 
\frac{\cos\Omega t}{\cos\vt_1} - \vp_1 ,
\end{equation}
and then \eqref{eq6.11} and \eqref{eq6.6} imply
\begin{align}
\label{eq6.13}
d_c^2 &= \big( \cos^{-1} (\cos\vt_1 \cos\psi)\big)^2 - \left( \cos^{-1} \frac{\cos\Omega t}{\cos\vt_1} - \vp_1\right)^2 \\
&= - (\psi-\vp_1)^2 + \big( \cos^{-1} (\cos\vt_1 \cos\psi)\big)^2 .\nonumber 
\end{align}
In particular we have derived

\begin{614Theorem}
Suppose $\vt_1(Q)\ne 0$.
Let $i\psi (Q)$ denote the critical point of $g$ with the smallest
modulus.
Then
\addtocounter{equation}{1}
\begin{equation}
\label{eq6.15}
g(i\psi) = \frac12 d_c (1,Q)^2 .
\end{equation}
\end{614Theorem}

One may rewrite $d_c$.  From \eqref{eq6.12}
\[
\frac{\Omega t}{\sin\Omega t} \sqrt{\sin^2 \Omega t - \sin^2 \vt_1} = \cos^{-1}
\frac{\cos\Omega t}{\cos\vt_1} - \vp_1 ,
\]
$\Omega = \Omega_{0,2}$, and then \eqref{eq6.11} yields
\begin{align*}
d_c^2 
= \left( \frac{\Omega t}{\sin\Omega t}\right)^2 \sin^2 \vt_1 &=
\frac
{\Big( \cos^{-1} \frac{\cos\Omega t}{\cos\vt_1} -\vp_1\Big)^2}
{1 - \Big( \frac{\cos\Omega t}{\cos\vt_1}\Big)^2}\ 
\tan^2\vt_1 \\
&= \frac{(\psi-\vp_1)^2}{\sin^2\psi} \, \tan^2\vt_1 ,\qquad \vt_1 \ne 0 .
\end{align*}
We note that
\[
\vp_1 (\vt_1) = 
\psi - \frac{\cos^{-1} (\cos\vt_1 \cos\psi)}{\sqrt{1-\cos^2 \vt_1 \cos^2\psi}}\, \cos\vt_1 \sin\psi
\]
increases from $0$ to $\psi$.
Hence $0\leq \vp_1 < \psi \leq \pi $, see \eqref{eq6.5}.

\begin{616Theorem}
Let $0 \leq v < \pi$.
With $g=g_1 + ig_2 $, $g_1,g_2\in\BR$, one has
\addtocounter{equation}{1}
\begin{equation}
\label{eq6.17}
\frac{\partial g_1}{\partial u} \ \begin{cases}
\ >0,&\text{$u\ne 0$,}\\
\ =0,&\text{$u=0$,}
\end{cases}
\end{equation}
and all of $g$'s critical points are imaginary.
\end{616Theorem}

\begin{proof}
One has
\begin{equation}
\label{eq6.18}
\frac{\partial g}{\partial u} = w -i \vp_1 - G(\vt_1 , w) ,
\end{equation}
where we set
\begin{equation}
\label{eq6.19}
G(\vt_1 , w) = \frac{\cosh^{-1} (\cos\vt_1 \cosh w)}{\sqrt{\cos^2 \vt_1 \cosh^2 w - 1}}\, \cos\vt_1 \sinh w .
\end{equation}

1) $u=0$. Here,
\[
\frac{\p g}{\p u}\Big|_{u=0} = i (v-\vp_1) - \frac
{\cosh^{-1} (\cos\vt_1 \cos v)}{\sqrt{1-\cos^2 \vt_1 \cos^2 v}}\, 
\cos\vt_1 \sin v .
\]
If $\cos\vt_1 \cos v=x$, then $\cosh^{-1} x = i\cos^{-1} x = i\gamma$,
\[
\frac{\p g}{\p u}\Big|_{u=0} = i\Big( v - \vp_1 - \frac{\g}{\sin\g} \sqrt{\cos^2 \vt_1 - \cos^2 \g}\Big) ,
\]
and
\begin{equation}
\label{eq6.20}
\frac{\p g_1}{\p u} \Big|_{u=0} = 0 .
\end{equation}

2) $u>0$. 
We set $G=G_1+iG_2$, $G_1,G_2\in\BR$, and note that
\begin{equation}
\label{eq6.21}
\frac{\p g_1}{\p u} \Big|_{\vt_1 = 0} = 0 ,\  \text{or}\ G_1 (0,w) = u .
\end{equation}
It is our intention to show that $G_1$ is a decreasing function of $\vt_1$ which implies Theorem 6.16;
we shall discuss the $v=0$ and $v>0$ cases separately.

(i) $v=0$. $(\al)$ When $\cos\g = \cos\vt_1 \cosh u <1$, one has
\begin{equation}
\label{eq6.22}
G_1(\vt_1,u) = G(\vt_1,u) = \frac{\cos^{-1} (\cos\vt_1 \cosh u)}{\sqrt{1-\cos^2 \vt_1 \cosh^2 u}}\,
\cos\vt_1 \sinh u ,
\end{equation}
and
\begin{equation}
\label{eq6.23}
\frac{\p G_1}{\p \vt_1} = -\sin\vt_1 \sinh u \frac{2\g - \sin 2\g}{2\sin^3 \g} < 0 , 
\quad \text{if}\ \vt_1 >0 .
\end{equation}

$(\beta)$ $\cosh \kappa = \cos\vt_1 \cosh u > 1$.  Here
\begin{align}
\label{eq6.24}
\frac{\p G_1}{\p \vt_1} &= \frac{\p}{\p\vt_1} \frac
{\cosh^{-1} (\cos\vt_1 \cosh u)}{\sqrt{\cos^2 \vt_1 \cosh^2 u-1}}\, \cos\vt_1 \sinh u \\
&= -\tan\vt_1 \sqrt{\cosh^2 \kappa - \cos^2 \vt_1} \, \frac{\sinh 2\kappa - 2\kappa}{2\sinh^3 \kappa} \nonumber \\
&< 0, \quad \vt_1 > 0 . \nonumber
\end{align}

(ii) $v>0$. Set $\cos\vt_1 \cosh w = \cosh z$, $z=\tilde u + i\tilde v$ and $\beta^2 = \sin^2 \vt_1$;
to simplify the notation we shall write $\tilde u + i \tilde v = x+iy$.
Then
\begin{equation}
\label{eq6.25}
\frac{\p G}{\p\vt_1} = \tan\vt_1 \sqrt{\sinh^2 z + \beta^2} \, \frac{z-\sinh z \cosh z}{\sinh^3 z} ,
\end{equation}
\begin{align}
\label{eq6.26}
&|\sinh z|^6 \cot\vt_1 \frac{\p G}{\p \vt_1}  \\
&= \sqrt{\sinh^2 z + \beta^2}  \cdot \nonumber \\
&\qquad \cdot \big( z \, \overline{\sinh^3 z} - |\sinh z|^2 |\cosh z|^2 \, \overline{\cosh z} + |\sinh z|^2 \cosh z\big) 
\nonumber \\
&= \sqrt{ \sinh^2 z + \beta^2} (A - iB) ,\nonumber
\end{align}
with
\begin{align}
\label{eq6.27}
A = &x\sinh x \cos y (\sinh^2 x \cos^2 y - 3\cosh^2 x \sin^2 y) \\
&+ y \cosh x \sin y (3\sinh^2 x \cos^2 y - \cosh^2 x \sin^2 y ) \nonumber \\
&- \cosh x \cos y (\sinh^2  x + \sin^2 y) (\sinh^2 x - \sin^2 y) ,\nonumber
\end{align}
\begin{align}
\label{eq6.28}
B = &x\cosh x \sin y (3\sinh^2 x \cos^2 y - \cosh^2 x \sin^2 y) \\
&- y \sinh x \cos y (\sinh^2 x \cos^2 y - 3\cosh^2 x \sin^2 y ) \nonumber \\
&- \sinh x \sin y (\sinh^2  x + \sin^2 y) (\cosh^2 x + \cos^2 y) ;\nonumber
\end{align}
here we used
\begin{align}
\label{eq6.29}
&|\sinh (x+iy)|^2 = \cosh^2 x - \cos^2 y = \sinh^2 x + \sin^2 y ,\\
\label{eq6.30}
&|\cosh (x+iy)|^2 = \cosh^2 x - \sin^2 y = \sinh^2 x + \cos^2 y .
\end{align}
With $\sinh^2 z + \beta^2 = \rho + i\delta = \omega$, \eqref{eq3.67} yields
\begin{equation}
\label{eq6.31}
\sqrt{\sinh^2 z + \beta^2} = \sqrt{\frac{|\omega|+\rho}{2}} + i\sqrt{\frac{|\omega|-\rho}{2}} ,\qquad \de > 0 ,
\end{equation}
so,
\begin{equation}
\label{eq6.32}
(\cot\vt_1) |\sinh z|^6 \,\frac{\p G_1}{\p\vt_1} = \sqrt{\frac{|\omega|+\rho}{2}} \,A + \sqrt{\frac{|\omega|-\rho}{2}}\, B ;
\end{equation}
$\delta > 0$ when $0 < y < \pi/2$ which is equivalent to $0 < v < \pi/2$.
We note that
\begin{equation}
\label{eq6.33}
A < 0 \ \text{when}\ 0 < y < \frac{\pi}{2} \ \text{and}\  x> 0 .
\end{equation}
Indeed,
\begin{align}
\label{eq6.34}
A &= x\sinh^3 x \cos^3 y - 3x \sinh x \cosh^2 x \sin^2 y \cos y\\
    &\qquad\qquad\qquad\qquad +3 \sinh^2 x \cosh x \sin y (y \cos y) \cos y \nonumber \\
    &\qquad\qquad\qquad\qquad - y \cosh^3 x \sin^3 y \nonumber \\
    &\quad\  - \sinh^4 x \cosh x \cos y + \cosh x \sin^4 y \cos y \nonumber \\
    &= -\frac12 \sinh^3 x \cos y (\sinh 2x - 2x) - x\sinh^3 x \sin^2 y \cos y \nonumber \\
    &\quad\  - 3 \sinh x \cosh x \sin y \cos y \big( x \cosh x \sin y - \sinh x (y\cos y)\big) \nonumber \\
    &\quad\  - \frac12 \cosh x \sin^3 y (2y - \sin 2y) - y\sinh^2 x \cosh x \sin^3 y \nonumber \\
    &< 0 .\nonumber
\end{align}
As for $B$, $B(y=x) = x^4 + O(x^6) >  0$, $B(y=-x) = -x^4 + O(x^6)<0$ when $x$ is small.
Next,
\begin{equation}
\label{eq6.35}
2 \sqrt{\frac{|\omega|+\rho}{2}} \, (\cot \vt_1 ) |\sinh z|^6 \frac{\p G_1}{\p \vt_1} =
|\omega| A + \rho A + \de B < 0 , 
\end{equation}
$0 \leq y \leq \frac{\pi}{2}$.
We argue as follows.
With
\begin{equation}
\label{eq6.36}
\omega = \beta^2 + \sinh^2 (x+iy) 
= \beta^2 -\sin^2 y + (\sinh x ) \sinh (x+i2y) ,
\end{equation}
one has
\begin{equation}
\label{eq6.37}
|\omega| \geq |\beta^2 - \sin^2 y | + \sinh^2 x .
\end{equation}
Indeed, when $\beta^2 < \sin^2 y$,
\begin{align*}
&|\omega |^2 - (\sin^2 y - \beta^2 + \sinh^2 x)^2 \\
&= (\beta^2 - \sin^2 y + \sinh^2 x \cos 2y)^2 + \sinh^2 x \cosh^2 x \sin^2 2y \\
&\qquad - (\sin^2  y - \beta^2 + \sinh^2 x)^2 \\
&= 4\beta^2 \sinh^2 x \cos^2 y  .
\end{align*}
When $\beta^2 >\sin^2 y$, an analogous calculation yields
\[
|\omega|^2 - (\beta^2 - \sin^2 y + \sinh^2 x)^2 = 4(1-\beta^2)\sinh^2 x \sin^2 y ,
\]
and we have derived \eqref{eq6.37}.
In view of this one rewrites the right hand side of \eqref{eq6.35}:
\begin{align}
\label{eq6.38}
&|\omega|A + \rho A + \de B \\
&= \big( |\omega| - (\sin^2 y - \beta^2 + \sinh^2 x)\big) A + 2(\sinh^2 x \cos^2 y) A + \de B .\nonumber
\end{align}
The first term is negative, as follows from \eqref{eq6.37} and \eqref{eq6.33}, and
the sum of the second and third terms,
\begin{align}
\label{eq6.39}
&2 (\sinh^2 x \cos^2 y ) A + \de B \\
&= 2 \sinh x \cos y \big( (\sinh x \cos y) A + (\cosh x \sin y)B\big) \nonumber \\
&= 2 \sinh x \cos y (\sinh^2 x + \sin^2 y ) C ,\nonumber
\end{align}
with
\begin{align}
\label{eq6.40}
C &= x (\sinh^2 x \cos^2 y - \cosh^2 x \sin^2 y) \\
&\quad\ + 2y \sinh x \cosh x \sin y \cos y - \sinh x \cosh x (\sinh^2 x + \sin^2 y) \nonumber \\
&= \frac12 \sinh^2 x (2x \cos^2 y - \sinh 2x ) \nonumber \\
&\quad\ - \cosh x \sin y \big( x \cosh x \sin y - \sinh x (y\cos y)\big) \nonumber \\
&\quad\ + \sinh x \cosh x \sin y (y\cos y - \sin y ) \nonumber \\
&<0 \nonumber
\end{align}
also negative
and we have derived \eqref{eq6.17} when $0\leq y = \tilde v \leq \pi/2$.
In particular we have shown that
\begin{align}
\label{eq6.41}
\Ree g (u) &= \frac12 \big( u^2 - (v-\vp_1)^2\big) - \frac12 (\tilde u^2 - \tilde v^2 ) \\
&= \frac12 (u^2 -\tilde u^2 ) - \frac12 (v-\vp_1)^2 + \frac12 \tilde v^2 \nonumber
\end{align}
is an increasing function of $|u|$ when $0\leq v\leq \pi/2$.
$\tilde v (u)$ decreases to $v$ as $|u|\to \infty$, so $u^2 - \tilde u (u)^2$ is an increasing function of $|u|$.

\begin{sloppypar}
When $\pi / 2 < v < \pi$, $\tilde v(u)$ increases to $v$,
and $u^2 - \tilde u (u,v)^2 = u^2 - {\tilde u (u,\pi-v)^2}$ is also an increasing function of $|u|$ since
$0<\pi -v \leq \pi/2$.
Consequently, \eqref{eq6.41} is an increasing function of $|u|$ on $\pi/2 < v < \pi$, and we
have completed the proof of Theorem 6.16. \qedhere
\end{sloppypar}
\end{proof}

We are ready to derive the small time asymptotic of $P_C$.
Write
\begin{equation}
\label{eq6.42}
P_C = \frac{e^{\frac{n^2}{2}t}}{(2\pi t)^{n+1}}\, \sum_{k=-\infty}^\infty\, S_k ,
\end{equation}
see \eqref{eq4.24}.
The first step is the derivation of the small time asymptotic of $S_0$ and
we start by
replacing $u+i\vp_1$ with $u$ in the integrand,
\begin{equation}
\label{eq6.43}
S_0 = \int_{-\infty}^\infty e^{-\frac{g(u)}{t}} 
V_n \big( \cosh^{-1} (\cos\vt_1 \cosh u),t\big) du ,
\end{equation}
see \eqref{eq6.1}.
This is the change of variable formula if the path of integration is $(-\infty + i\vp_1 , \infty+i\vp_1)$, and
we claim that it may be moved to the real axis.
Indeed, let $u=v+i\ve$, $0\leq \ve \leq \pi - \de$, $\de > 0$, and $\cos\vt_1 \cosh (v + i\ve)=\cosh
(\tilde v + i\tilde \ve)$.
Then $2\Ree g = v^2 - \tilde v^2 - (\ve - \vp_1)^2$, and $e^{- \Ree g / t} < e^{\pi^2 / 2t}$.
Also, the estimates \eqref{eq4.86} imply that $|V_n (\tilde v + i\tilde\epsilon , t)|^2 \leq
C(t) (\tilde v^2 + \pi^2)^n / (\sinh^2 \tilde v + \sin^2 \de)^n$, so the integrand vanishes exponentially
as $|u|\to \infty$, uniformly in $\ve$, $0\leq \ve \leq \pi -\de$, hence \eqref{eq6.43}.

The work on \eqref{eq6.43} follows the proof of Theorem 2.4.2 of \cite{3}.
$g(u+i\psi)$ attains its global minimum at $u=0$;
according to Theorem 6.16 it is a strictly increasing function of $|u|$
with a global minimum at $u=0$.
Set
\begin{equation}
\label{eq6.44}
\Phi (u) = g(u+i\psi) - g(i\psi) ,\\
\end{equation}
where $i\psi$ is the first critical point of $g$ with $\psi > 0$.
Then
\begin{align}
\label{eq6.45}
\Phi' (0) &=g_u (i\psi) = 0 ,\\
\label{eq6.46}
\Phi'' (0)&= g_{u^2} (i\psi) \\
&= \frac{\sin^2 \vt_1}{1-\cos^2\vt_1 \cos^2\psi}
\left( 1 - \frac{\cos\vt_1 \cos\psi \cos^{-1} (\cos\vt_1 \cos\psi)}{\sqrt{1-\cos^2\vt_1 \cos^2 \psi}}\right)\nonumber \\
&> 0 ,\qquad \vt_1 > 0 ;\nonumber
\end{align}
compare with \eqref{eq6.10}.
Therefore
\begin{equation}
\label{eq6.47}
\Ree \Phi (u) \geq Cu^2  ,\qquad C>0 , \ |u| < 1 .
\end{equation}
We move the path of integration in $S_0$ of \eqref{eq6.43} to $(-\infty + i\psi , \infty + i\psi)$,
\begin{align}
\label{eq6.48}
S_0 &= e^{-\frac{g(i\psi)}{t}} \left\{ \int_{-\de}^\de + \int_{|u|>\de}\right\} e^{-\frac{\Phi (u)}{t}} \cdot \\
&\quad\ \cdot V_n \big( \cosh^{-1} (\cos\vt_1 \cosh (u+i\psi)), t\big) du \nonumber \\
&= S_\de + S'_\de ,\qquad 0 < \de < 1 ,\nonumber
\end{align}
and note that the error term $S'_\de$ is bounded by
\[
|S'_\de| \leq e^{-\frac{g_1 (\de + i\psi)}{t}} \int_{-\infty}^\infty \big | V_n \big( 
\tilde u (u) + i \tilde \psi (u) , t\big) \big| du \leq Ce^{-\frac{g_1 (\de  + i\psi)}{t}} ;
\]
in view of \eqref{eq4.86}.
As for $S_\de$, one introduces a new variable $z$ by
\begin{equation}
\label{eq6.49}
\Phi (u) = \Phi'' (0) u^2 \big( 1 + O (|u|)\big) = \Phi'' (0) z^2 ,
\end{equation}
which leads to
\begin{align}
\label{eq6.50}
e^{\frac{g(i\psi)}{t}} S_\de &= \int_{z (-\de)}^{z(\de)} \, e^{-\frac{\Phi'' (0) z^2}{t}} 
V_n \big( \cosh^{-1} \big( \cos\vt_1 \cosh (u(z)+i\psi)\big),t\big) du \\
&= \int_{-\de}^\de \, e^{-\frac{\Phi''(0) z^2}{t}} V_n \big( \cosh^{-1} \big( \cos\vt_1 \cosh ( u(z)+i\psi)\big),t\big) du 
\nonumber \\
&\quad \ + O \left( e^{-\frac{c}{t}}\right) ,\qquad c>0 ;\nonumber
\end{align}
the integration contour in \eqref{eq6.50} may be complex, but the integrand is holomorphic in $z$, so by
moving the contour to the real axis the error committed is $\sim \exp \big( - \Phi'' (0) \de^2 / 2t\big)$.
Hence
\begin{align}
\label{eq6.51}
e^{\frac{g(i\psi)}{t}} S_\de
&= \int_{-\de}^\de \, e^{-\frac{\Phi''(0)z^2}{t}} v^n\big( \cos^{-1} (\cos\vt_1 \cos\psi)\big) dz + O(t) \\
&= \int_{-\infty}^\infty \, e^{-\frac{\Phi''(0)z^2}{t}} v^n \big( \cos^{-1} (\cos\vt_1 \cos\psi)\big)dz +O(t) \nonumber \\
&= \sqrt{\frac{\pi t}{\Phi''(0)}}\, v^n \big( \cos^{-1} (\cos\vt_1 \cos\psi)\big) + O (t) ,\nonumber
\end{align}
and substituting for $\Phi''(0)$ and $v$ one has

\begin{652Theorem}
Let $0<\vt_1 < \pi/2$, $0\leq \vp_1 < \pi$.
Then
\addtocounter{equation}{1}
\begin{equation}
\label{eq6.53}
P_C (t,1,Q) = \frac{\sin (\eta_{0,2}) v (\eta_{0,2})^n e^{-\frac{d_c (1,Q(\vt,\vp))^2}{2t}} (1+O(t))}
{\sqrt 2 (2\pi t)^{n+\frac12} \sin\vt_1 \sqrt{1-\eta_{0,2}\cot\eta_{0,2}}}\, , 
\end{equation}
where $\eta_{0,2} =\Omega_{0,2} t$.
\end{652Theorem}

\begin{proof}
\eqref{eq5.51} implies
\begin{equation}
\label{eq6.54}
d_c = \eta_{0,2} \frac{\sin\vt_1}{\sin\eta_{0,2}} < \eta_{0,2} < \pi ,
\end{equation}
so, in view of \eqref{eq4.85}, it suffices to show that individual $S_k$'s are dominated by $\exp(-\pi^2/2t)$
when $k\ne 0$.

(i) $\vp_1 > 0$.
Here
\begin{align*}
|S_k| &\leq \left| \int_{-\infty}^\infty \, e^{-\frac{u^2-(\tilde u + i\tilde\vp_1 + i2k\pi)^2}{2t}}
V_n (\tilde u + i\tilde\vp_1 + i2k\pi , t)du\right| \\
&\leq e^{-\frac{(\tilde\vp_1 + 2k\pi)^2}{2t}} \int_{-\infty}^\infty \, \big| V_n \big( \tilde u + i (\tilde \vp_1 +
2k\pi),t\big) \big| du \\
&\leq Ce^{-\frac{\pi^2}{2t}} ,\qquad k\ne 0 ,
\end{align*}
see \eqref{eq4.86}.

(ii) When $\vp_1 =0$, \eqref{eq4.94} and its extension to $I_k^{(n)} + I_{-k}^{(n)}$  yield
\[
| S_k + S_{-k} | \leq Ce^{-\frac{4k^2 \pi^2}{2t}} , \qquad k = 1,2,\ldots .\qedhere
\]
\end{proof}

\begin{655Theorem}
$\vt_1 = 0$, $0< |\vp_1| < \pi$.
Then
\addtocounter{equation}{1}
\begin{equation}
\label{eq6.56}
P_C = \frac{|\vp_1|^{n-1}}{2^n \G (n) t^{2n}} e^{-\frac{d_c ( 1, Q (\vt,\vp))^2}{2t}} \big( 1 + O(t)\big) .
\end{equation}
\end{655Theorem}

\begin{proof}
We shall work with $\vp_1>0$; when $\vp_1 < 0$ the derivation is analogous.
\begin{align*}
P_C &= \frac{e^{\frac{n^2}{2}t}}{(2\pi t)^{n+1}}\, \sum_{k=-\infty}^{\infty} 
\int_{-\infty}^\infty e^
{-\frac{u^2 - (u+i\vp_1 +i2k\pi)^2}{2t}}
V_n (u+i\vp_1 + i2k\pi , t) du \\
&= \frac{2\pi i}{(2\pi t)^{n+1}} \left\{ \sum_{k=-\infty}^\infty \, e^{-\frac{(\vp_1 +2k\pi)^2}{2t}} J_k \right\}
\big( 1 + O (t)\big) ,
\end{align*}
where we set
\[
J_k = \frac{1}{2\pi i} \int_{-\infty}^{\infty} \,
\frac{(z+i\vp_1 + i2k\pi)^n}{\sinh^n (z+i\vp_1)} e^{\frac{iz(\vp_1 +2k\pi)}{t}} dz .
\]
The singularities of the integrand are $z_\ell = i(\ell \pi - \vp_1)$, $\ell\in\BZ$, and
$\sinh (z+i\vp_1) = (-1)^\ell (z-z_\ell) (1+\cdots )$,
$z_\ell + i\vp_1 + i2k\pi = i(\ell+2k)\pi$, hence
\begin{equation}
\label{eq6.57}
J_k \! =\!  e^{\frac{iz_\ell (\vp_1 + 2k\pi)}{t}}\, \frac{1}{2\pi i} \! \int_{-\infty}^\infty\,
\frac{(z-z_\ell + i(\ell +2k)\pi)^n}{(-1)^{\ell n} (z-z_\ell)^n(1+\cdots )}
e^{\frac{i(z-z_\ell)(\vp_1 + 2k\pi)}{t}} \! dz .
\end{equation}
One expands the numerator of the integrand in powers of $z-z_\ell$, collects the
coefficients of $(z-z_\ell)^{n-1}$ and keeps the terms with the highest negative power
of $t$ only.
This yields
\[
J_k = e^{\frac{iz_\ell (\vp_1 + 2k\pi)}{t}} 
\left\{ \sum_\ell \, \frac{(i(\ell+2k)\pi)^n}{(-1)^{\ell n}} \, 
\frac{i^{n-1} (\vp_1 + 2k\pi)^{n-1}}{\G(n) t^{n-1}}\right\} \, \big( 1 + O(t)\big) ,
\]
hence one has
\begin{align}
\label{eq6.58}
P_C = \frac{(-1)^{n(\ell-1)}}{2^n \G (n) t^{2n}}\,\Big\{ &\sum_{k,\ell}\, (\ell+2k)^n (\vp_1 + 2k\pi)^{n-1}\cdot \\
&\cdot e^{ \frac{\vp_1^2 - 4k(k+\ell)\pi^2 - 2\ell\pi\vp_1}{2t}} \Big\}\, \big( 1 + O(t)\big) .\nonumber
\end{align}
The exponent of the integrand in \eqref{eq6.57} must have a negative real part.
This allows for two sums and we set $P_C = P_{C,1} + P_{C,2}$.

1) $P_{C,1}$ is the sum over $k=0,1,2,\ldots$ and $\ell=1,2,\ldots$.
In this summation the largest exponent in
\eqref{eq6.58} occurs when $k=0$, $\ell=1$, so
\[
P_{C,1} = \frac{\vp_1^{n-1}}{2^n\G (n) t^{2n}} \, e^{-\frac{(2\pi - \vp_1)\vp_1}{2t}} \big( 1 + O(t)\big)  ,
\]
which agrees with \eqref{eq6.56}; see also \eqref{eq5.50}.

2) $P_{C,2}$ is the sum over $k=-1,-2,\ldots$ and $\ell=0,-1,-2,\ldots$.
The largest exponent occurs when $k=-1$ and $\ell=0$ which gives the following
exponential factor for $P_{C,2}$:
\begin{equation}
\label{eq6.59}
e^{-\frac{4\pi^2-\vp_1^2}{2t}} .
\end{equation}
This is smaller than the exponential  factor in $P_{C,1}$ hence $P_{C,2}$ does not contribute
to the principal part of the small time asymptotic of $P_C$.
\end{proof}

\noindent {\bf 6.60 Remark.}
Replacing $\tan^{-1} (\sqrt{1-\al^2}\,\tan\Omega t)$ by its integral form \eqref{eq5.5} one may extend
Theorem 6.8 to a correspondence between $\Omega$ and $\psi$, where $2k\pi + \vt_1 \leq \Omega t \leq
(2k+1)\pi - \vt_1$ and $2k\pi \leq \psi \leq (2k+1)\pi $.

\bibliographystyle{amsalpha}
\bibliography{}

\bibliographystyle{amsalpha}

\providecommand{\bysame}{\leavemode\hbox to3em{\hrulefill}\thinspace}

\end{document}